*The Birth of Number Theory (Book VII of Euclid's* Elements*)*
*from the Arithmetization of Pythagorean Music*

Stelios Negrepontis, Vassiliki Farmaki, Angeliki Pisimisi

***Abstract***
Book VII of Euclid's *Elements* represents a remarkable mathematical achievement, marking the birth of number theory; although this book falls short of explicitly stating the Fundamental Theorem of Arithmetic (that every natural number can be written uniquely as the product of primes), it nevertheless contains all the tools necessary for its proof: namely, the Principle of the Least – mathematically equivalent to the modern Principle of Mathematical Induction – and the method of arithmetical *anthyphairesis* for finding the greatest common divisor of two natural numbers, commonly known as the Euclidean algorithm.
Our work presents novel arguments in support of the thesis that Book VII evolved directly from early Pythagorean arithmetized music. The dominant ancient accounts, which attributed this arithmetization to acoustical experiments devised by Pythagoras himself, have been definitively shown to be fictitious by Vincenzo Galilei. In contrast, the alternative acoustical experiments and the construction of Hippasus' 4–chord – consisting of four bronze cylinders with identical bases and heights 6, 8, 9, 12 – are validated as physically correct by Euler's Law for pipes. This correct experimental foundation led not only to the arithmetization of the fundamental musical intervals but, even more significantly, as we argue, to the discovery of the fundamental concept of *anthyphairesis* (initially in its musical, multiplicative form).
By applying Aristotle's *Topics* 158b24–29 Principle regarding the dynamic interaction between Postulates and Definitions, we reconstruct how musical (multiplicative) *anthyphairesis* was transferred – through the generation of an upward inductive ladder of ratios, starting with multiple and epimoric ratios, as recounted by Theon of Smyrna – to its arithmetical (additive) counterpart. Furthermore, we show that the numerous mathematical peculiarities in the definitions, statements of propositions, and proofs in Book VII, when compared with a modern treatment of the Fundamental Theorem of Arithmetic and the theory of positive rational numbers, possess a convincing explanation only through their musical origin. Ultimately, Hippasus emerges as the pivotal figure who discovered Pythagorean musical *anthyphairesis* (which later evolved into Philolaus' *Fragment 6*), and who, as we also know, transferred arithmetical *anthyphairesis* to its geometric counterpart, which was crucial both for Pythagorean incommensurability and for the two principles of the Infinite and the Finite in Pythagorean philosophy.

*Introduction*

Book VII of Euclid's *Elements* represents an impressive mathematical achievement. By employing the Principle of the Least – mathematically equivalent to the modern Principle of Mathematical Induction – introducing the fundamental algorithm of *anthyphairesis* of two natural numbers for finding their greatest common divisor, and culminating in Propositions VII.30 and VII.31, it essentially establishes the Fundamental Theorem of Arithmetic. The purpose of this paper is to present novel arguments in support of the thesis that Book VII of the *Elements* evolved directly from early Pythagorean arithmetized music.

The arithmetization of early empirical Pythagorean music was traditionally credited to Pythagoras himself (in Nicomachus' account); however, the experiments described there were conclusively shown to be fictitious by Vincenzo Galilei (Sections 1, 2).

Alternative acoustical experiments involving vessels, described by Theon and attributed to Lasus and Hippasus, are of uncertain physical validity. In contrast, the construction of a 4–chord in the form of four bronze cylinders with identical bases and heights 6, 8, 9, 12 – described in an *Anonymous Scholion* based on the authority of Aristoxenus and attributed to Hippasus – proves that the experimental foundation on which Hippasus based his construction of the 4–chord is entirely valid, as is conclusively demonstrated by Euler's Law for pipes.

However, Hippasus' construction of the 4–chord is of even greater importance because, in relation to these validated acoustical experiments, it incorporates the birth of the fundamental concept of *anthyphairesis* in its musical, multiplicative form. The anthyphairetic nature of Hippasus' 4–chord is strongly confirmed by the later Pythagorean Philolaus' *Fragment 6*, which extends this musical anthyphairesis by two more steps to generate both the *diesis* and the *Pythagorean comma* (cf. Boethius, 1989) (Sections 3, 4, 5).

Based on the properties embodied in this 4–chord and in these earlier accounts, we reconstruct the basic properties of the arithmetized Pythagorean theory of music (Section 6).

A backward application of Aristotle's *Topics* 158b24–29 Principle regarding the dynamic relation between Postulates and Definitions clarifies the sweeping Pythagorean conjecture that generalized the properties of the four musical intervals of the 4–chord to all arithmetical ratios (Sections 8, 9, and 10). The musical *anthyphairesis* was transferred to its arithmetical counterpart – the Euclidean algorithm – by gradually constructing multiple, epimoric, and epimeric ratios with an *anthyphairesis* of length one, two, three, and so forth, following an upward ladder as suggested by Theon of Smyrna (Section 11).

However, the contents of Book VII, when compared with a modern text in *number theory*, present a long array of striking peculiarities: numbers are represented as straight line segments; addition of ratios is absent; multiplication of ratios is restricted to a special form; the definition of proportion $a/b = c/d$ is given only for unequal numbers, encumbered by the unnecessary distinction between "part" and "parts"; and finally, the proof of the commutativity of multiplication is highly unexpected (Section 7). Ultimately, all these peculiarities of Book VII find a systematic and convincing explanation in their musical origin (Section 12).

It must be noted that Szabó (1978) had the key insight that arithmetical *anthyphairesis* originated in Pythagorean music, but he failed to ground this insight on the essential arguments developed in our work (Section 13).

In conclusion, Hippasus emerges as the pivotal figure of Pythagorean mathematics and philosophy who:

(a) shaped the early Pythagorean arithmetized theory of music (followed by Philolaus) on the basis of musical (multiplicative) *anthyphairesis*;

(b) employed infinite geometrical *anthyphairesis* in order to prove the incommensurability of the diameter to the side of a square; and

(c) formulated an early version of the two Pythagorean philosophical Principles of the *Infinite* and the *Finite*, which possess a clear anthyphairetic content and were inspired by the incommensurability proof.

We therefore conclude that it is highly improbable that the Pythagorean who accomplished the intermediate step of transferring musical *anthyphairesis* to arithmetic – thereby discovering arithmetical *anthyphairesis* (the Euclidean algorithm) – could be anyone other than Hippasus (Section 15).

## **1.** *Ancient Empirical Music*

Nicomachus' problematic account of Pythagoras' arithmetization of ancient empirical music is nevertheless instructive in showing how this empirical music was conceived by an experienced practical musician. Each of the four fundamental musical intervals had the role of an empirical equivalence class, whose elements were the dichords producing this musical interval. The class of dichords in each of the musical intervals satisfied the properties of empirical transitivity, empirical most equality, and possibly the arithmetical property: if a dichord ($a$, $b$) produces a musical interval, then the dichord ($ka$, $kb$) belongs to the same class for any natural number $k$. In addition, there is an empirical operation of composition on 3–chords, with corresponding equivalence classes, subject to the rule of commutativity.

***1.1.*** *Each of the four fundamental musical intervals, the octave, the fifth, the fourth, and the tone, is an empirical equivalence class of dichords, satisfying the transitivity property*

Since ancient, certainly pre–Greek times, professional musicians were able to recognize and isolate certain fundamental combinations of sounds as pleasing to the ear. These musical intervals were produced in various ways by musical instruments, but principally by plucking two unequal strings/chords of a dichord simultaneously. In early Pythagorean music, there were just four basic musical intervals out of which music was composed: namely, the octave, the fifth, the fourth, and the tone. Thus, many centuries before the appearance of mathematics, professional musicians, who had no knowledge of or interest in mathematics, were able to conceive of music in a pre–mathematical manner as follows: a musical interval is an empirical equivalence class, and the class of all dichords, in their various forms, that produce this musical interval, constitutes the elements or representatives of the equivalence class, representing the mutually equivalent dichords. Thus, again, a musical interval is the empirical equivalence class of all dichords that produce it.
This relation between dichords and musical intervals is clearly described in Nicomachus' account of the hypothetical discovery of the arithmetization of musical intervals by Pythagoras himself in Chapter 6 of his work *Harmonicum enchiridion*, in which Pythagoras, depicted as passing by a blacksmith's workshop, recognized (*epeginoske*) the four basic musical intervals that he heard (*epekouse*) being produced by hammers beating iron, the equivalent of dichords plucked simultaneously by a lyre player.

> παρὰ τι χαλκοτυπεῖον περιπατῶν, ἔκ τινος δαιμονίου συντυχίας
> ἐπήκουσε ῥαιστήρων σίδηρον ἐπ' ἄκμονι ῥαιόντων,
> καὶ τοὺς ἤχους παραμὶξ πρὸς ἀλλήλους συμφωνοτάτους
> ἀποδιδόντων, πλὴν μιᾶς συζυγίας.
> ἐπεγίνωσκε δ' ἐν αὐτοῖς
> τὴν δὲ διὰ πασῶν καὶ τὴν διὰ πέντε καὶ τὴν διὰ τεσσάρων συνῳδίαν.
> τὴν δὲ μεταξύτητα τῆς τε διὰ τεσσάρων καὶ τῆς διὰ πέντε
> ἀσύμφωνον μὲν ἑώρα αὐτὴν καθ' ἑαυτήν,
> συμπληρωτικὴν δὲ ἄλλως τῆς ἐν αὐτοῖς μείζονος.
> Nicomachus, *Harmonicum enchiridion* 6.1, 11–20
>
> and happening by some heaven–sent chance (*daimoniou suntuchias*) to walk by a blacksmith's workshop (*chalkotupeion*),
> *he heard (epekouse) the hammers (rhaisteron) beating (rhaionton) iron (sideron)* on the anvil (*ep' akmoni*)
> and giving out sounds (*echous*) fully concordant (*sumphonotatous*)
> in combination with one another,
> with the exception of one pairing (*suzugias*):

and *he recognized (epeginoske)* among them the consonance (*synoidian*) of *the octave* and those of *the fifth* and *the fourth*.
*He noticed* that what lies in *between (metaxuteta) the fourth and the fifth [the tone]* was in itself discordant (*asumphonon*), but was essential in filling out the greater of these intervals.
[Nicomachus, *Harmonicum enchiridion*, trans. Barker, 1989, p. 256, with modifications by the authors]

It is clear that the musical intervals divide the (good) dichords into four *disjoint* equivalence classes, and thus the property of Transitivity is empirically established. By the property of *Transitivity* we mean that if the dichords ($a$, $b$) and ($c$, $d$) produce the same musical interval, and the dichords ($c$, $d$) and ($e$, $f$) produce the same musical interval, then the dichords ($a$, $b$) and ($e$, $f$) produce the same musical interval.

**1.2.** *The empirical Most Equality property of the musical intervals*

Nicomachus' account in the *Harmonicum enchiridion* continues in the passage 6.1, 20–35 with the description of another fundamental property of musical intervals:

ἄσμενος δὴ ὡς κατὰ θεὸν ἀνυομένης αὐτῷ τῆς προθέσεως
εἰσέδραμεν εἰς τὸ χαλκεῖον, καὶ ποικίλαις πείραις
παρὰ τὸν *ἐν τοῖς ῥαιστήρσιν ὄγκον* εὑρὼν τὴν διαφορὰν τοῦ ἤχου,
ἀλλ' οὐ παρὰ τὴν τῶν ῥαιόντων βίαν
οὐδὲ παρὰ τὰ σχήματα τῶν σφυρῶν
οὐδὲ παρὰ τὴν τοῦ ἐλαυνομένου σιδήρου μετάθεσιν,
*σηκώματα ἀκριβῶς* ἐκλαβὼν καὶ *ῥοπὰς ἰσαιτάτας* τῶν ῥαιστήρων
πρὸς ἑαυτὸν ἀπηλλάγη.
καὶ ἀπὸ τινός ἑνὸς πασσάλου διὰ γώνων ἐμπεπηγότος τοῖς τοίχοις,
ἵνα μὴ καὶ ἐκ τούτου διαφορά τις ὑποφαίνηται
ἢ ὅλως ὑπονοῆται πασσάλων ἰδιαζόντων παραλλαγή,
ἀπαρτήσας τέσσαρας χορδάς
*ὁμοΐλους καὶ ἰσοκώλους, ἰσοπαχεῖς τε καὶ ἰσόστροφους,*
ἑκάστην ἐφ' ἑκάστης ἐξήρτησεν,
*ὁλκὴν* προσδήσας ἐκ τοῦ κάτωθεν μέρους.
τὰ δὲ μήκη τῶν χορδῶν μηχανησάμενος ἐκ παντὸς ἰσαίτατα.
Nicomachus, *Harmonicum enchiridion* 6.1, 20–35

Overjoyed at the way his project had come, with god's help, to fulfillment,
he ran into the smithy (*chalkeion*),
and *through a great variety of experiments (peiras), he discovered that what stood in direct relation to the difference in the sound (diaphoran tou echou) was the weight of the hammers (ton en tois rhaistersin ogkon),* *not* the force (*bian*) of the strikers
*or* the shapes of the hammer–heads (*sphuron*)
*or* the alteration (*metathesin*) of the iron, which was being beaten. He *weighed* them *accurately (akribos),*
and took away for his own use
pieces of metal *exactly equal (isaitatas) in weight (rhopas) to the hammers.*
Then he fixed a single rod from corner to corner under his roof,
so that *no variation (me diaphora tis)* should arise
or even be suspected of arising
from the *peculiarities (parallage) of different rods*,
and hung from it four strings,
each of the *same material (homoulous),*
and consisting of an *equal number of strands (isokolous),*
and each of *equal thickness (isopacheis)*
and *twisted to the same extent (isostrophous)* as each of the others.
He then attached a *weight (holken)* to the lower part of each string.
And having so contrived it that the length of every string was in all respects *absolutely equal (isaitata).*
[Nicomachus, *Harmonicum enchiridion*, trans. Barker, 1989,
pp. 256–257, with modifications by the authors]

It is clear that Pythagoras' main concern, due to the numerous experiments he had performed, was to obtain weights *exactly (akribos), most equal (isaitatas)* to those of the hammers in the smithy that produced the desired musical intervals. Thus, he makes clear that two equivalent dichords must satisfy a rigorous *Most Equality* property: if ($a$, $b$) is a dichord producing some musical interval, then there are certainly many other ($c$, $d$) dichords, equivalent to it, producing the same musical interval, in ways yet unknown, but if we take a dichord of the form ($a$, $x$), then if we want ($a$, $x$) to be equivalent to ($a$, $b$), we must insist that $x$ is exactly equal (*isaitata*) to $b$, it will not do for $x$ to be merely close to $b$.

As we shall examine in Section 2, Pythagoras' claim that the Most Equality property applied to the *weights* of the dichords is seriously flawed, a fact that throws into doubt the genuineness of Pythagoras' experiments; the concern with the Most Equality property is valid not for the weights of the chords/hammers, but for their *lengths*.

### **1.3.** *The operation of Composition of dichords and musical intervals, and the empirical Commutativity property of Composition*

#### **1.3.1.** *The operation of the Composition of dichords and musical intervals*

Furthermore, the stringed musical instruments allowed the operation of *composition* for two dichords and their corresponding musical intervals to be formulated, as follows:
If ($a$, $b$, $c$) is a 3–chord, the musical interval produced by the dichord ($a$, $c$) is the *composition* of the musical interval produced by the dichord ($a$, $b$) with that produced by the dichord ($b$, $c$).
We denote this composition by $(a, c) = (a, b) * (b, c)$.
The operation is empirically *well defined* in the sense that equivalent dichords composed with equivalent dichords produce equivalent dichords.
Since the fundamental Pythagorean musical intervals were just four – the octave, the fifth, the fourth, and the tone – the possibilities for the operation of composition were limited to essentially just two:

a) an octave is the composition of a fifth and a fourth,
   octave = fifth * fourth, and
b) a fifth is the composition of a fourth and a tone,
   fifth = fourth * tone.

**1.3.2.** *The empirical commutativity property of the composition of intervals*

A fundamental empirical rule that was certainly observed since old times was the *commutativity of composition*; indeed, in music, composition is a well–defined operation not only on dichords but also on musical intervals. We can formulate the commutativity property as follows:

*Commutativity of Composition of Musical Intervals.*
If $a$, $b$, $c$, and $d$, $e$, $f$ are two trichords such that:
the dichords ($a$, $b$) and ($e$, $f$) produce the same musical interval, and
the dichords ($b$, $c$) and ($d$, $e$) produce the same musical interval,
then the composition $(a, c) = (a, b) * (b, c)$
and the composition $(d, f) = (d, e) * (e, f)$ produce the same musical interval.
We note that this property is analogous to the statement of Proposition V.23 in Euclid's *Elements* on the Perturbed Proportion (*tetaragmene analogia*) for ratios of magnitudes.

**1.4.** *The empirical equivalence of the dichords (a, b) and (ka, kb)*

Although there is no direct mention of this in ancient sources, it is highly probable that the ancient musicians had noticed that, if a dichord ($a$, $b$) produces one of the four musical intervals, then, for any natural number $k$, the dichord ($ka$, $kb$) also produces the same musical interval.

**2.** *The acoustical experiment, leading to the arithmetization of music and the construction of the 4–chord 6, 8, 9, 12, and attributed to Pythagoras,* *reported by Nicomachus in* Harmonicum Enchiridion*, Chapter 6, shown to be fictitious*

**2.1.** *The Arithmetization of Music according to Nicomachus'* Harmonicum Enchiridion*, Chapter 6*

The first part of Nicomachus' report on the acoustical experiments attributed to Pythagoras is given in Section 1.1, above. The report continues as follows:

τὴν μὲν γὰρ ὑπὸ τοῦ μεγίστου ἐξαρτήματος *τεινομένην*
πρὸς τὴν ὑπὸ τοῦ μικροτάτου διὰ πασῶν φθεγγομένην κατελάμβανεν.
ἦν δὲ ἡ μὲν δωδεκά τινων ὁλκῶν, ἡ δὲ ἕξ.
ἐν *διπλασίῳ* δὴ λόγῳ ἀπέφαινε τὴν *διὰ πασῶν*,
ὅπερ καὶ αὐτὰ τὰ βάρη ὑπέφαινε.
τὴν δ' αὖ μεγίστην πρὸς τὴν παρὰ τὴν μικροτάτην (οὖσαν ὀκτὼ ὁλκῶν)
*διὰ πέντε* συμφωνοῦσαν, ἔνθεν ταύτην ἀπέφαινεν ἐν *ἡμιολίῳ* λόγῳ,
ἐν ᾧπερ καὶ αἱ ὁλκαὶ ὑπῆρχον πρὸς ἀλλήλας.
πρὸς δὲ τὴν μεθ' ἑαυτὴν μὲν τῷ βάρει, τῶν δὲ λοιπῶν μείζονα,
ἐννέα σταθμῶν ὑπάρχουσαν, τὴν *διὰ τεσσάρων*, ἀναλόγως τοῖς βρίθεσι.
καὶ ταύτην δὴ *ἐπίτριτον* ἄντικρυς κατελαμβάνετο,
ἡμιολίαν τὴν αὐτὴν φύσει ὑπάρχουσαν τῆς μικροτάτης
(τὰ γὰρ ἐννέα πρὸς τὰ ἕξ οὕτως ἔχει),
ὅνπερ τρόπον ἡ παρὰ τὴν μικρὰν ἡ ὀκτώ
πρὸς μὲν τὴν τὰ ἕξ ἔχουσαν ἐν ἐπιτρίτῳ ἦν,
πρὸς δὲ τὴν τὰ δώδεκα ἐν ἡμιολίῳ.
τὸ ἄρα μεταξὺ τῆς διὰ πέντε καὶ τῆς διὰ τεσσάρων,
τουτέστιν ᾧ ὑπερέχει ἡ διὰ πέντε τῆς διὰ τεσσάρων,
ἐβεβαιοῦτο ἐν ἐπογδόῳ λόγῳ ὑπάρχειν,
ἐν ᾧπερ τὰ ἐννέα πρὸς τὰ ὀκτώ.
Nicomachus, *Harmonicum enchiridion* 6.1, 35–56

He then plucked (*teinomenen*) strings two at a time in turn,
and found the concords previously mentioned,
a different concord for each pairing.
He perceived that the string was under tension from the biggest object attached sounded as *octave*
in relation to the one under tension from the smallest.
The former was of *twelve* units of weight (*holkon*), the latter of *six.*
Hence, he showed that *the octave is in duple ratio,*
as the weights (*bare*) themselves implied.
He found that the *biggest* [12] sounded at a *fifth* in relation to
the smallest but one which had *eight units* [8] of weight (*holkon*),

and revealed from this that *the fifth is in a hemiolic ratio* [3/2],
the ratio in which these weights (*holkai*) stood to each other.
In relation to the one second in weight (*toi barei*) to itself
and greater than the others, which was of *nine* units (*stathmon*),
*it* [12] sounded at the interval of a *fourth*, in conformity with the relations of the weights (*brithesi*).
And he at once perceived that this ratio was *epitritic,* and
that *this same string* [12] was in a *hemiolic* ratio to *the smallest*
(since that is *the ratio of 9 to* 6):
and in the same way, the smallest but one, carrying eight units, stood
in the epitritic ratio to the one that carried six, and
in hemiolic ratio to that which carried twelve.
And hence he established that what lies *between the fourth and the fifth,* that is, that by which the fifth exceeds the fourth, is in *epogdoic* ratio,
that in which nine units stand to eight.
[Nicomachus, *Harmonicum enchiridion*, trans. Barker, 1989, p. 257, with modifications by the authors]

The discovery that Nicomachus' account assigns to Pythagoras consists in:
(a) the arithmetization of the four basic musical intervals, and their representation as *absolutely exact simple ratios of numbers in least terms,* namely the octave by the ratio 2/1, the fifth by the ratio 3/2, the fourth by the ratio 4/3, and the tone by the ratio 9/8, and
(b) the construction of the Pythagorean 4–chord, with chords equal to the numbers 6, 8, 9, 12 units of weight.

**2.2.** *The Validity of Early Pythagorean Acoustical Experiments*

**2.2.1.** *Summary of Early Pythagorean Acoustical Experiments*

We summarize Pythagoras' experiments.
Each of the dyads of the four hammers in the blacksmith's workshop acted like dichords. Activated by beating iron (rather than being beaten by it), they produced, by their nature, each of the four musical intervals: the octave, the fifth, the fourth, and the tone. Pythagoras is supposed to have determined that the only reason the dyads of hammers were producing musical intervals was the difference in their *weights*. Thus, in his experiments, the number assigned to a chord/note by Pythagoras was essentially a *weight* hung from it.
Thus,
the chord/note with the name *hypate* has weight 6,
the chord/note with the name *mese* has weight 8,
the chord/note with the name *paramese* has weight 9, and
the chord/note with the name *nete* has weight 12.

The arithmetization of the fundamental musical intervals attributed to Pythagoras by Nicomachus, which we examined in Section 2.1, provides the correct arithmetical values, but the validity of the experiments leading to these values has been put into serious question. As will be detailed in the remainder of Section 2.2, these Pythagorean accounts faced early skepticism from Ptolemy and were later conclusively refuted by the experiments of Vincenzo Galilei, leading modern scholars to dismiss the weight experiments as fictions. As an alternative to these questionable accounts, in Section 3 we will explore alternate, valid musical experiments attributed by Theon of Smyrna to Lasus and Hippasus, including the construction of the 4–chord 6, 8, 9, 12 attributed to Hippasus by a source going back to Aristoxenus.

**2.2.2.** *Ptolemy in his* Harmonics *expressed some doubts about the validity of Pythagoras' acoustical experiments*

> ἐπί τε τῶν ἐξαπτομένων ταῖς χορδαῖς βαρῶν
> μὴ διασῳζομένων ἀπαραλλάκτων ἀλλήλαις παντάπασι τῶν χορδῶν,
> ὁπότε καὶ πρὸς αὑτὴν ἑκάστην οὕτως ἔχουσαν εὑρεῖν ἔργον,
> οὐκέτι δυνατὸν ἔσται τοὺς τῶν βαρῶν λόγους ἐφαρμόσαι τοῖς γινομένοις δι' αὐτῶν ψόφοις τῷ καὶ τὰς πυκνοτέρας καὶ λεπτοτέρας ἐν ταῖς αὐταῖς τάσεσιν
> ὀξυτέρους φθόγγους ποιεῖν.
> πολὺ δὲ ἔτι πρότερον κἂν ταῦτά τις ὑπόθηται δυνατὰ
> καὶ ἔτι τὸ μῆκος τῶν χορδῶν ἴσον, τὸ μεῖζον βάρος τῇ πλείονι τάσει τὴν τῆς ἀρτώσης αὐτῷ χορδῆς διάστασιν αὐξήσει τε καὶ πυκνώσει μᾶλλον,
> ὥστε καὶ διὰ τοῦτο συμπίπτειν τινὰ παρὰ τὸν λόγον τῶν βαρῶν ἐν τοῖς ψόφοις ὑπεροχήν.
> Ptolemy, *Harmonics*, 1, 8, 11–20

> In the case of weights attached to strings, where the strings are not kept in all respects identical with one another – since it is a hard job to find strings of which each is in this condition, even with respect to itself – it will no longer be possible to fit the ratios of the weights to the sounds that arise through them, since denser and finer strings under the same tension make higher notes. Much more important even than that is the fact that even if one assumes that these things are possible, and again that the lengths of the strings are equal, the bigger weight by its greater tension will increase the length of the string attached to it, and will make it denser, so that from this too will arise a difference in the sounds that is not in accordance with the ratio of the weights.
> [translation Barker, 1989, p. 291]

**2.2.3.** *Vincenzo Galilei in 1581, 1589 was the first to prove experimentally that Pythagoras' experiments were definitely invalid*

As Papadopoulos (2025) notes:

> Of particular note among Galilei's treatises is his *Dialogo della musica antica et moderna*.(2) He is credited with being the first to observe that to obtain the octave of a note emitted by a stretched string, the tension of the string must be multiplied by four and not by two, contradicting the Pythagorean creed that the ratio of integers 1/2 is universally associated with the octave. Vincenzo Galilei describes such an experiment in his *Discorso intorno alla diversita delle forme del Diapason*, written around 1589, a treatise in which he examines the different ways of producing the octave (called, in ancient Greek terminology, the diapason, hence the title of the book).(3) His experiment consisted in attaching a system of weights to one end of a vibrating string and observing that in order to obtain the octave while keeping the length constant, the weights must be quadrupled, not doubled. Similarly, to obtain the fifth, one must multiply the weights by 4/9, and not by 2/3. Vincenzo Galilei's discovery is a consequence of the general law of string vibration stated by Marin Mersenne a few decades later, which says that the frequency of vibrations produced by a stretched string is inversely proportional to the square of the string's tension, and not to the tension itself. Mersenne's law is sometimes considered to be the first known non–linear law of physics.
>
> > (2) [Vincenzo Galilei, *Dialogo della musica antica e della moderna*, Florence, G. Marescotti, 1581, Eng. trans. by C.V. Palisca, Dialogue on ancient and modern music, in Humanism in Renaissance musical thought, New Haven–London, Yale University Press, 1985].
> >
> > (3) [See V. Galilei, *Discorso intorno alla diversità delle forme del Diapason*, Eng. trans. by C.V. Palisca, A special discourse concerning the diversity of the ratios of the diapason, in Claude V. Palisca, The Florentine Camerata: documentary studies and translations, New Haven, Yale University Press, 1988].
>
> (Papadopoulos, 2025, p. 3)

Thus, while Galilei's experiments refute the myth of Pythagoras' room with suspended weights, they equally invalidate the accounts of the experiments in the smithy.

**2.2.4.** *The discovery of Vincenzo Galilei was followed by Marin Mersenne, who discovered by experiments his famous «formule de Mersenne», and by Christiaan Huygens, who gave a mathematical proof of this formula*

Building upon Vincenzo Galilei's empirical observations, Marin Mersenne formalized this relationship into a general physical law. In his *Harmonie universelle* (1636), Mersenne stated what is often considered the first non–linear law of physics: the frequency of a vibrating string is proportional to the square root of its tension. However, as Papadopoulos (2025, p. 10) notes, Mersenne sought a rigorous mathematical explanation for this phenomenon. In a letter dated November 16, 1646, Mersenne challenged the seventeen–year–old Christiaan Huygens to prove mathematically why achieving an octave requires quadrupling the weight, while simply halving the string's length achieves the same result (De Waard, 1948).
Many years later, in 1673, Huygens successfully provided the definitive mathematical proof of this law in his work *Découverte de la théorie de l'isochronisme des vibrations*. In his writings, Huygens was unequivocally dismissive of the ancient Pythagorean myths concerning weights and strings. In his essay *Origine du chant*, he explicitly states that the accounts of Pythagoras' weight experiments are false, emphasizing the quadratic nature of the relationship:

> Quelques anciens auteurs de musique racontent qu'après cela il [Pythagore] attacha des poids suivant ces proportions trouvées à des cordes pour les tendre, et qu'il trouva que le poids double tendait la corde à l'octave plus haut, et le sesquialtère à la quinte, ce qui est faux: et si ces auteurs s'étaient donné la peine de faire l'expérience ils auraient trouvé qu'il faut le poids quadruple du premier pour faire monter une corde à l'octave, qu'il faut qu'il soit comme 9 à 4 pour faire la quinte, et qu'universellement la raison des poids doit être double de celle qui détermine les consonances par les parties d'une corde tendue. (Papadopoulos et al., in Huygens, 2021, p. 273)

To further clarify the physical laws underlying Huygens' critique, the editors of the modern edition of his musical writings provide two highly illuminating footnotes. They explain that Huygens' observation is deeply connected to what he called the "rule of bell founders" (*la règle des fondeurs*) and to "Mersenne's formula". Specifically, they point out that for three–dimensional resonant objects like bells, producing an octave requires doubling the linear dimensions, which consequently demands octupling the weight, proving once again that simple linear weight ratios cannot produce the fundamental musical consonances. This quadratic relationship between tension and frequency, which Mersenne had empirically observed, was ultimately proven mathematically, by Huygens, twenty–seven years after Mersenne's initial challenge. The editors' original annotations are as follows:

> – Il est en effet certain que ce n'est pas le poids double, mais le poids octuple d'un objet résonnant semblable du même métal qui donne l'octave. Or, ceci était fort bien connu. Dans ses *Harmonicorum libri* («Harmonic. Instrumentorum lib. 4 de campanis», Prop. VII) Mersenne dit que les dimensions linéaires d'une cloche doivent être doublées «ut... campana... habeatur quae facit Octauam cum prima». Dans le «Harmonia lib. IV» (p. 364 des *Cogitata physico–mathematica* de 1644) il dit expressément: «Praeter hoc in eo libro [manuscrit de Ioannes Faber] mihi placuit primo quod plurimis obseruationibus nitatur, quibus recte concludit 4 malleos in ea ratione, quam Pythagorae tribuunt, diapasonis divisionem in Quintam et Quartam minime facere, atque adeo falsum esse hinc illum rationes harmonicas desumpsisse». Huygens fait mention, dans une Pièce de 1672 [t. XIII, dernier alinéa de la p. 804], de «la règle des fondeurs», qui doublent le diamètre des cloches qu'ils veulent avoir à l'octave l'une de l'autre.
> – Cette règle [qu'universellement la raison des poids doit être double de celle qui détermine les consonances par les parties d'une corde tendue] fait partie de la fameuse «formule de Mersenne» qui énonce en particulier que la fréquence de la note produite par une corde vibrante est proportionnelle à la racine carrée de la tension de cette corde. Mersenne, qui avait observé cette règle, demanda à Christiaan Huygens, dans une lettre qu'il lui adressa le 16 novembre 1646, d'en trouver une explication rationnelle (le passage est cité au chapitre 1, §1.4 du présent volume). C'est vingt –sept ans plus tard que Huygens rédigea une démonstration mathématique de cette loi; voir la pièce *Découverte de la théorie de l'isochronisme des vibrations*, chap. 16, p. 389 de ce volume.
> (Papadopoulos et al., in Huygens, 2021, p. 273)

**2.2.5.** *The impact of Vincenzo Galilei's musical experiments on his son Galileo Galilei*

Vincenzo Galilei's debunking of the Pythagorean musical myths – as recounted by Nicomachus – by experimentation was a watershed turning point; its impact was not confined only to correcting musical theory, but extended much beyond, by leaving an indelible imprint on his son, Galileo, and instilling in him the conviction that rigorous experimentation is the cornerstone of modern physics.

As Papadopoulos (2025) notes:

> Vincenzo Galilei conducted his own experiments. In this context, as in others, he taught his son Galileo the importance of experimentation, and of not relying blindly on established theories, no matter how old they may be. Through his experiments, we get a glimpse of Galilei's "technician" side, a technician in the original sense of the Greek word *technê*, meaning both "art" and "craft". The following statement, extracted from the same treatise, has been considered as an indication of Galilei's adherence to the "experimental method" that was to take over science a few decades later: "There are few things that cannot be weighed, counted or measured".
>
> Galilei's experiments led him to distrust dogmas, or at least not to blindly accept them, even those of well–established authorities, as he wrote in the *Dialogo della musica antica et della moderna*:
>
> In connection with [the theories of Pythagoras] I wish to point out two false opinions of which men have been persuaded by various writings, and which I myself shared until I ascertained the truth by means of experiment, the teacher of all things.
>
> (Papadopoulos, 2025, pp. 3–4)

**2.2.6.** *Modern scholars, including van der Waerden, 1943, Szabó, 1978, Barker, 1989, Levin, 1994 have noted the fictitious nature of Pythagoras' experiment*

Modern scholars correctly note that Pythagoras' alleged weight experiments are physically invalid. As is well established today, the strictly one–dimensional length of the vibrating string, rather than the suspended weight it carries (which is governed by a quadratic proportion), is the correct measure of the musical interval produced. However, it is a striking historiographical oversight that these prominent scholars point out the physical and mathematical impossibility of the Pythagorean myth without ever mentioning that it was Vincenzo Galilei who first demonstrated this fallacy experimentally in the late sixteenth century.

Van der Waerden, 1943; Szabó, 1978; Barker, 1989; Levin, 1994 correctly note that this experiment is most likely invalid, since, among others, as we know today, the one–dimensional length of the chord, and not the three–dimensional weight that the chord carries, is the correct measure of the musical interval produced.

B. L. van der Waerden (1943) points out the mathematical error in the accounts of the weights:

> Falls man tatsächlich versucht haben sollte, die Zahlenverhältnisse mit Hilfe angehängter Gewichte exakt zu messen, wie es für

> Pythagoras überliefert ist, so hat man dabei einen neuen Mißerfolg buchen müssen, da die Gewichte bei der Oktave sich nicht wie 1:2, sondern wie 1: $\sqrt{2}$ verhält
> (Van der Waerden, 1943, p. 173)

A. Szabó (1978) notes:

> If one had actually tried to carry out this experiment and measure the numerical ratios (of the consonances) by means of suspended weights, as Pythagoras is said to have done, then the attempt would inevitably have ended in failure. For the weights corresponding to the octave are in the ratio 1:$\sqrt{2}$ not 1: 2. (Szabó, 1978, p. 122, Section 2.5)

A. Barker (1989) emphasizes the legendary character of these accounts:

> The tradition that Pythagoras discovered the fundamental intervallic ratios is very persistent… though others attributed the discovery to his disciples….
>
> Many of the stories told in this connection (including the present one …) are plainly fictions. We cannot be sure even if the discovery was first made by Pythagoreans, though they made notable use of it. It may have been due to practical musicians or instrument–makers, who would have found it useful in the construction of harps (these became popular in Greece during the sixth century B.C., but were widespread in eastern Mediterranean cultures from much earlier times). Most scholars are prepared to give Pythagoras the benefit of the doubt, at least so far as the dissemination of these ideas among the Greeks is concerned. (Barker, 1989, p. 256)
>
> This betrays the legendary character of the tale of this harmonious blacksmith. It is not true that the ratios between the pitches will correspond to those between the weights of the hammer–heads. The desired results will again fail to be produced. The pitch ratios are not directly related to the ratios of the weights, but to those of their square roots.
> (Barker, 1989, p. 257)

F. R. Levin (1994) provides a definitive conclusion on the matter:

> In sum, then, the smith's hammers, the weights attached to the strings and their supposed equivalents – tensions on the strings – have no place in the world of acoustical fact. They are, in every respect, fabulous elements in a long–standing legend. (Levin, 1994, p. 93)

**3.** *The acoustical experiments with vessels by Lasus and Hippasus, leading to the arithmetization of three musical intervals, reported by Theon of Smyrna 59,7–21 and the construction of the 4–chord 6, 8, 9, 12 by Hippasus, reported in the anonymous Scholion to Plato's* Phaedo *108d4 ("skill of Glaucus"), a Scholion whose origin goes back to Aristoxenus*

**3.1.** *The acoustical experiments with vessels by Lasus and Hippasus*

Lasus of Hermione and the Pythagorean Hippasus of Metapontium performed some epoch–making simple experiments that revealed the relation between music – until then an ancient but strictly empirical art – and arithmetical ratios, then in a primitive stage, a discovery that eventually marked, as this paper demonstrates, the birth of number theory.

*Λᾶσος δὲ ὁ Ἑρμιονεύς*, ὥς φασι,
καὶ οἱ περὶ τὸν *Μεταποντῖνον Ἵππασον*, Πυθαγορικὸν ἄνδρα,
συνέπεσθαι τῶν κινήσεων τὰ τάχη καὶ τὰς βραδυτῆτας,
δι’ ὧν αἱ συμφωνίαι…
ἐν ἀριθμοῖς ἡγούμενος λόγους τοιούτους ἐλάμβανεν ἐπ’ ἀγγείων.
ἴσων γὰρ ὄντων καὶ ὁμοίων πάντων τῶν ἀγγείων,
τὸ μὲν κενὸν ἐᾶσας, τὸ δὲ ἥμισυ ὑγροῦ πληρώσας,
ἐψόφει ἑκατέρῳ, καὶ αὐτῷ ἡ *διὰ πασῶν* ἀπεδίδοτο συμφωνία·
θάτερον δὲ πάλιν τῶν ἀγγείων κενὸν
ἐῶν, εἰς θάτερον τῶν τεσσάρων μερῶν τὸ ἕν ἐνέχεε,
καὶ κρούσαντι αὐτῷ ἡ *διὰ τεσσάρων* συμφωνία ἀπεδίδοτο·
ἡ δὲ *διὰ πέντε*, ὅτε ἓν μέρος τῶν τριῶν συνεπλήρου,
οὔσης τῆς κενώσεως πρὸς τὴν ἑτέραν
*ἐν μὲν τῇ διὰ πασῶν ὡς β' πρὸς ἕν,*
*ἐν δὲ τῷ διὰ πέντε ὡς γ' πρὸς β',*
*ἐν δὲ τῷ διὰ τεσσάρων ὡς δ' πρὸς γ'.*
Theon Smyrneus, *Expositio rerum mathematicarum ad legendum Platonem utilium, [De utilitate mathematicae]* 59, 7–21

*Lasus of Hermione*, so they say, and
the followers of *Hippasus of Metapontum*,
a Pythagorean (*Puthagorikon andra*),
pursued the speeds and slownesses of the movements,
through which the concords arise
...
Thinking that... in numbers (*en arithmois*),
he constructed ratios (*logous*) of these sorts in vessels.
All the vessels were equal (*ison*) and similar (*homoion*).
[*Note*: here again the Most Equality property]
Leaving one empty (*kenon*)
and filling (*plerosas*) the other up to *halfway* with liquid,
he made a sound (*epsophei*) on each, and the concord of *the octave (diapason)* was given out for him.
Then, again, leaving one of the vessels empty, he poured into (*enechee*) the other one part out of the four,

and when he struck (*krousanti*) it the concord of *the fourth (dia tessaron)* was given out for him,
as was *the fifth (dia pente)*
when he filled up one part out of the three.
The one empty space stood to the other
in the octave as 2 to 1, in the fifth as 3 to 2, and in the fourth as 4 to 3.
[Translation by Barker, 1989, p. 218, with modifications by the authors]

Lasus of Hermione was a musician and music theorist of the late 6th century B.C., renowned for his transformative contributions to musical practice. According to the 10th–century Byzantine encyclopedic lexicon, the *Suda*, his birth is traditionally placed during the 58th Olympiad, spanning the years 548 to 544 B.C. In his *Histories* (VII.6), Herodotus records that Lasus of Hermione caught the diviner Onomacritus in the act of forging an oracle of Musaeus. This event occurred shortly before the assassination of Hipparchus in 514 B.C., confirming that Lasus was already a prominent and influential figure in the Athenian court during the late 6th century B.C. (VII.6, line 15)

Barker (1989) notes:

> Lasus of Hermione was a distinguished late sixth–century musician, said to have written the first book on music (see the entry under his name in the Suda, cf. Martianus Capella De Nupt. IX.936), credited with experiments in acoustics (see I.4 Theon Smyrn. 59.4) and innovations in musical practice, especially in Athens (ps.–Plut. De Mus.1141b–c, Herodotus Hist.VII.6, Athenaeus Deipn.455c, 624e–f).
> (Barker, 1989, p. 128)

As Iamblichus notes in his work *On the Pythagorean Life (chapter 18, section 81, line 5 – section 82, line 7)*, Hippasus of Metapontum was a key successor to Pythagoras who fundamentally shaped the *Mathematici* branch of the school. Moving beyond the oral, unproven traditions of the *Acousmatici*, Hippasus emphasized a rigorous mathematical and experimental approach. This distinction is crucial, as it identifies him not merely as a follower of dogma, but as a pioneer of the scientific method that bridged harmony and number theory. According to Huffman (1993, p. 8), "Hippasus of Metapontum lived approximately between 530 and 450 B.C.".
We may conjecture, on the basis of these dates and on the fact that Theon mentions Lasus first, that Lasus was the first to devise these acoustical experiments and that Hippasus, in a later time, verified them.

We note that the experiments conducted by Lasus and Hippasus utilize identical vessels; for them the only variable is the height of the empty part of the vessel, a strictly one–dimensional, linear variable presumably behaving like a string/line segment. In contrast to the undoubtedly fictitious experiments attributed by Nicomachus to Pythagoras, Theon's accounts regarding Lasus' and Hippasus' experiments are in need of careful assessment, which we will present in Sections 3.3 and 3.4.

**3.2.** *The construction of the 4–chord 6, 8, 9, 12 by Hippasus, reported in the anonymous Scholion to Plato's* Phaedo *108d4 ("skill of Glaucus"), a Scholion whose origin goes back to Aristoxenus*

> *Γλαύκου τέχνη· σημείωσαι παροιμίαν*
> *ἐπὶ τῶν μὴ ῥᾳδίως κατεργαζομένων ἤτοι ἐπὶ τῶν πάνυ ἐπιμελῶς καὶ ἐν τέχνῃ εἰργασμένων.*
> *Ἵππασος γάρ τις κατεσκεύασε χαλκοῦς τέτταρας δίσκους οὕτως,*
> *ὥστε τὰς μὲν διαμέτρους αὐτῶν ἴσας ὑπάρχειν,*
> *τὸ δὲ τοῦ πρώτου δίσκου πάχος*
> *ἐπίτριτον μὲν εἶναι τοῦ δευτέρου, ἡμιόλιον δὲ τοῦ τρίτου, διπλάσιον δὲ τοῦ τετάρτου·*
> *κρουομένους δὲ τούτους ἐπιτελεῖν συμφωνίαν τινά.*
> *Καὶ λέγεται Γλαῦκον ἰδόντα τοὺς ἐπὶ τῶν δίσκων φθόγγους*
> *πρῶτον ἐγχειρῆσαι δι' αὐτῶν χειρουργεῖν, καὶ ἀπὸ ταύτης τῆς πραγματείας*
> *ἔτι καὶ νῦν λέγεσθαι τὴν καλουμένην Γλαύκου τέχνην.*
> *Μέμνηται δὲ τούτων Ἀριστόξενος ἐν τῷ Περὶ τῆς μουσικῆς ἀκροάσεως*
> *καὶ Νικοκλῆς ἐν τῷ Περὶ θεωρίας.*
> *Scholion Platon Phaedo 108d, 1–12*
>
> It is said of things that are not accomplished easily,
> or of things that are made with great care and skill.
> For a certain *Hippasus* made four bronze discs in such a way that while their diameters were equal,
> the thickness of the first disc was
> *epitritic* [4/3] in relation to that of the second,
> *hemiolic* [3/2] in relation to that of the third,
> and *double* [2/1] that of the fourth,
> And it is said that when *Glaucus* noticed the notes made by the discs,
> he was the first to set himself to making music with them,
> and that it is as a result of this endeavour that people still speak of the "skill of Glaucus", as it is called.
> Aristoxenus and Nicocles recall this information in their respective works.
> [Translation by Barker, 1989, p. 31, with modifications by the authors]

The testimony regarding Hippasus' 4–chord goes back to Aristoxenus and carries considerable weight, as Barker notes:

> This Glaucus is Glaucus of Rhegium, a musical writer of the late fifth century… Here he is represented as a practicing musician rather than a historian and critic: hence his proverbial skill. The whole passage is based on a report by Aristoxenus (frag. 90), whose testimony carries weight.
> (Barker, 1989, p. 31)

**3.3.** *From Empirical Vessels to the Arithmetized 4–chord*

Regarding the acoustical experiments reported by Theon in Section 3.1, it is not explicitly recorded that *cylindrical* vessels were used; it is only noted that they were identical. There are some ancient cylindrical vessels, but Theon does not specify that the experiment used vessels that were cylindrical. However, Hippasus' 4–chord described in Section 3.2 is the direct evolution of these experiments. Specifically, an evolved form of Hippasus' vessel experiments can be identified within the 4–chord by the first, third and fourth chords. Thus,
the cylinder of height 6 corresponds to the vessel half–filled with liquid,
the cylinder of height 9 corresponds to the vessel filled by one–fourth, and
the cylinder of height 12 corresponds to the empty vessel.

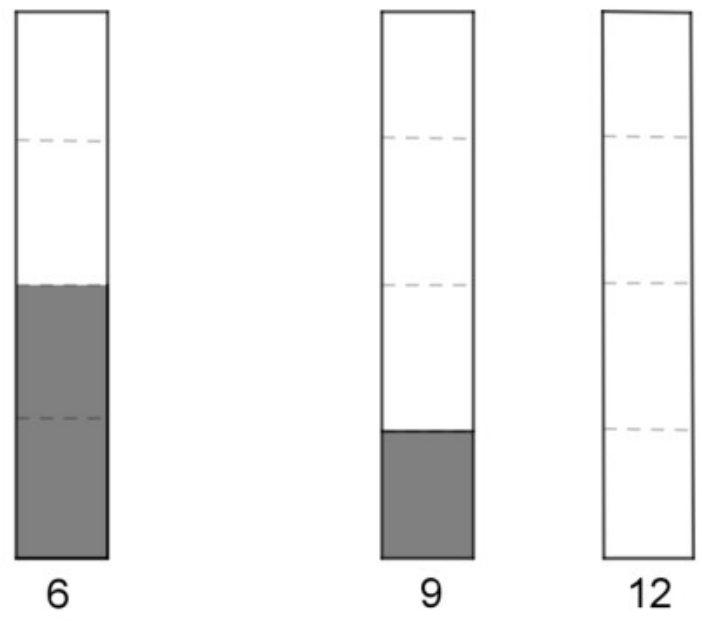

By comparing the accounts in Sections 3.1 and 3.2, a significant conclusion emerges. While Theon's account (Section 3.1) is vague about the vessels' exact shape, the description of the 4–chord (Section 3.2) explicitly specifies cylindrical forms, namely bronze discs with equal diameters.
This subsequent use of identical cylinders by Hippasus for the construction of his 4–chord leads us to the conclusion that the vessels employed in his earlier acoustical experiments must have also been cylindrical. Consequently, the empty air columns inside these identical cylindrical vessels acted essentially as musical pipes. Therefore, according to Euler's law, these air columns behaved mathematically like a linear vibrating string, where the pitch is determined by a strictly one–dimensional linear variable (length) rather than volume.
Our conclusion is that Hippasus could not have constructed the 4–chord unless he had eventually performed these acoustical experiments intended to arithmetize the three fundamental musical intervals (octave, fifth, and fourth), employing cylinders with identical bases. As we shall explore in Section 4, this construction had profound mathematical implications.

### **3.4.** *The validity of the acoustical experiments with vessels by Hippasus*

#### **3.4.1.** *Vincenzo Galilei's pipe experiments: The Walker–Palisca Debate and the genius of Hippasus*

Walker (1978) considers a question, set by V. Galilei:

> "In his unpublished *Discorso intorno alla diversità delle forme del diapason,* [Galilei] asks the question: *what interval would be given by two pipes of the same diameter but one of which is double the length of the other?*" (Walker, 1978, p. 24)

Palisca (1989) referring to the same passage in the *Discourse Concerning the Diapason,* [translated by Palisca, 1989, pp. 189–190], notes:

> "Galilei… asked the question: what sort of interval would two pipes make that have the same diameter but duple length?"
> (Palisca, 1992, p. 150, Footnote 19)

Granted that Hippasus' cylinders behave like musical pipes, the question Vincenzo Galilei asks here bears precisely on whether Hippasus' claim – that the two cylinders/pipes of height 2 to 1 give the musical interval of the octave – is true and correct.

To this question, Galilei answers that they would produce a major third (namely the musical interval 5/4). As Walker points out, Galilei's assumption is incorrect, since, according to the physical law mathematically formulated by Leonhard Euler in his *Dissertatio physica de sono* (1727), the pitch of a pipe is a function of its length, not its cubic capacity.

This answer sparked a significant debate among modern scholars, notably between Claude V. Palisca and D. P. Walker, as Capecchi and Capecchi (2026) also note.

Palisca (1992), argues that V. Galilei *"apparently made a slip"* in confusing diameter with length:

> A major third of the intense tuning of Aristoxenus, which, in fact, is the third part of the octave, whereas if they are in quadruple proportion, they would sound a minor sixth of the same intense tuning. *The question should have read: "have the same length but duple diameter."* Mersenne's table corroborates Galilei's findings with respect to volume, since twice the diameter will produce twice the volume. For the major third in the 5:4 proportion, the table gives the proportion of volumes as 125:64 (1.95); increasing the length slightly would make it the bigger third of the Aristoxenian system. Similarly, his ratio for the minor sixth, 8:5 (incorrectly

> printed as 6:5), is 512:125 (4.09). Cutting the length slightly gives Galilei's result of 4:1 for the Aristoxenian minor sixth.
> (Palisca, 1992, p. 150, Footnote 19)

Palisca (1992, p. 150) then points out that V. Galilei is guilty of a similar slip a little later in the same work:

> It is insufficient for two pipes to be in duple proportion, Galilei asserted (p. 195); this will give the major third of Aristoxenus (smaller than true); nor is quadruple proportion sufficient when comparing volumes, for this will sound a minor sixth. Galilei here did not specify which dimension was in duple proportion. *But he must have meant diameter or width of pipes, rather than length. Because two pipes otherwise equal that are in duple proportion of length will sound approximately an octave.* Mersenne also investigated the influence of diameter on pitch when other dimensions were equal. He found that *a pipe twice another in diameter* sounded *a minor third below*, though he noted that blowing the shorter pipe harder made that a *major third, which agrees with Galilei's result*. A pipe four times another produced a diminished seventh or major sixth. This again is not far from what Galilei observed. Vincenzo did not speak in this passage of the ratio of lengths, perhaps because he set out to prove that the duple ratio was not the only one that produced the octave. (Palisca, 1992, p. 150, Footnote 19)

*Our note*: if Palisca is correct and Galilei meant to vary the diameter rather than the length, then Galilei's question has no relation to the experiments of Lasus and Hippasus and the construction of the 4–chord, which strictly relied on varying the lengths while keeping the base diameter constant.

Walker (1978), on the contrary, believes that Galilei means exactly what he states, namely, he compares two pipes with the same diameter, one having double the length of the other. He argues that Galilei's answer is mistaken:

> In his unpublished *Discorso intorno alla diversità delle forme del diapason*, he asks the question: what interval would be given by two pipes of the same diameter but one of which is double the length of the other? and answers that it would be an equally tempered major third, which, by his own erroneous rule is correct.
> (Walker, 1978, p. 24).

*Our note*: if Walker is correct and Galilei's formulation was intentional, then the experiment described by Galilei is the same as the one performed by Hippasus with his bronze cylinders. In this scenario, Galilei is fundamentally wrong in his physical and musical conclusions, while Hippasus is correct. Here Walker refers to Euler's law since the pitch of a pipe is a function of its length and not of its cubic capacity.

In conclusion, Vincenzo Galilei does not present any experimental or theoretical reason to reject Hippasus' claim that his 4–chord, consisting of bronze cylinders of identical bases and heights 6, 8, 9, 12, was in any way musically deficient.

### **3.4.2.** *Theoretical vindication: From Hippasus to Euler's Dissertation De Sono*

In 1727, Leonhard Euler wrote his *Dissertatio physica de sono*, providing the first rigorous mathematical expression to show that the fundamental frequency of an organ pipe is proportional to the inverse of its length. He argued that the column of air inside a pipe does not act as a three–dimensional volume, but rather behaves in direct physical analogy to a vibrating string, establishing length – and not cubic volume – as the primary determinant of pitch. According to Euler, just as a regular musical string requires a physical weight to stretch it and create tension, these "air strings" are put under tension by the weight of the atmosphere. He acknowledged a natural difference: while a hanging weight pulls a regular string apart, atmospheric pressure does the opposite, compressing the air column and making it narrower. Yet, he argued that the physical effect is exactly the same, making the mathematical analogy completely legitimate. The key difference in a flute or pipe is that instead of vibrating from a single plucked point like a string, the entire column of air compresses and expands together along its whole length.

To turn this physical insight into a mathematical model, Euler established a ratio of tension $p$ to weight $q$, where $a$ represents the length of the pipe (corresponding to the length of a string), $p$ represents the atmospheric tension measured by the height of mercury in a barometer $k$, and $q$ is the actual weight of the air trapped inside the tube. By factoring in the relative density of mercury to air $n$:1, Euler proved that the ratio of tension to weight is given by the formula: $p/q = n \cdot k/a$.

This simple formula carries a profound physical conclusion: the frequency of the sound waves depends directly on the one–dimensional linear length ($a$) of the pipe, rather than its volume. By showing that the air column behaves mathematically like a vibrating string, Euler's law provides a solid theoretical foundation that vindicates Hippasus' ancient experiments with cylindrical vessels.

Thus, Euler theoretically proved what Hippasus had already demonstrated experimentally: that in cylindrical vessels acting as musical pipes, frequency is inversely proportional to the linear length of the empty space, definitively vindicating the one–dimensional approach to acoustics.

While this early work (Euler, 1727) established the inverse relationship between the length of a pipe and its frequency, Euler's later and more mature researches on the physics of sound are presented in his three landmark memoirs, *De la propagation du son* (Euler, 1766a)*, Supplément aux recherches sur la propagation du son* (Euler, 1766b), and *Continuation* *des recherches sur la propagation du son* (Euler, 1766c), as well as his final synthesis on the topic, *Éclaircissemens plus détaillés sur la génération et la propagation du son et sur la formation de l'écho* (Euler, 1767).

**3.4.3.** *Modern scholars, including Koestler (1959), Levin (1994), Huffman (1993), Ferguson (2011), Lloyd (2014), and Primavesi (2014), have noted the experimental validity of the acoustical experiments attributed to Hippasus, yet completely overlooked their mathematical consequence: the birth of anthyphairesis*

In stark contrast to the "fabulous" myths of Pythagoras' hammers and weights, modern scholars unanimously recognize the scientific rigor and physical validity inherent in the experiments with vessels attributed to Lasus and Hippasus. However, there is a profound historiographical blind spot in their analyses. While these prominent scholars fully understand and validate the acoustics of these experiments, they fail to grasp their ultimate historical and mathematical significance: they do not realize that it was precisely from these experimentally established intervals and the subsequent construction of the 4–chord that the concept of anthyphairesis emerged. The following excerpts demonstrate how modern scholars successfully validate the physical reality of the experiments, yet entirely miss their anthyphairetic consequence.

Barker (1989, p. 31) notes:

> In fact, the procedures described in those passages give the wrong results [attributed to Pythagoras], whereas that of Hippasus' works: if the discs have equal diameters, the frequencies are proportional to their thicknesses. The present report may well be historically sound.

D. Creese (2010, p. 94) writes on Hippasus' 4–chord:

> Unlike any of the other procedures involving metal objects (either sounded by percussion or as weights for strings, whether attributed to Pythagoras or not), this one will work as described.
>
> If the discs have an equal diameter, varying their thicknesses in the ratios described (12: 9: 8: 6 in lowest terms) will produce exactly the opposite result as varying the speaking length of a true string or a cylindrical pipe: the disc that is twice as thick will sound an octave above the first, and so on.
>
> The fact that the procedure will work as described does not in itself guarantee that Hippasus actually carried it out. We have already seen that not all fragments of Aristoxenus are credible. Paul Tannery, noting the lack of an earlier version of the story and the fact that Hippasus appears to have left no writings himself, treated the account with some doubt. Burkert was less dubious, and Barker treats Aristoxenus in this context as a "reliable

> authority". M. L. West notes that as far back as the eighth century B.C., there was a tradition of disc – chime – making in southern Italy, where Hippasus lived. This suggests that the materials and methods of manufacturing such discs would have been available to Hippasus.

D. Creese (2010, p. 96) continues:

> Our account of it, however, is devoid of the language either of empiricism or of demonstration. Hippasus merely "constructed" (*κατασκεύασε*) four bronze discs in certain dimensions, and a certain acoustic result followed from the manner of their construction. If Aristoxenus' choice of vocabulary tells us anything about the scientific context of the procedure (and it may well not), his use of the verb *kataskeuazein* could suggest that he thought of Hippasus as engaged in a mathematical project, in which a construction (*kataskeué*) can be the solution of a problem (*probléma*) or a stage leading toward a demonstration/proof (*apodeixis*). While it would be anachronistic to see Hippasus as a sort of harmonic Archimedes of the early fifth century, it is certainly possible that Aristoxenus thought of Hippasus' construction as a component of the sort of rigorous mathematical argument which had been developed in the fourth century.

Koestler (1959) emphasizes the general scientific importance of the discovery but misses the specific mathematical mechanism it spawned:

> The Pythagorean discovery that the pitch of a note depends on the length of the string which produces it, and that concordant intervals in the scale are produced by simple numerical ratios (2:1 octave, 3:2 fifth, 4:3 fourth, etc.), was epoch–making: it was the first successful reduction of quality to quantity, the first step towards the mathematization of human experience, and therefore the beginning of Science. (Koestler, 1959, p. 28)

Levin (1994) acknowledges the validity of Hippasus' experiments in contrast to the legends:

> Thus, if, as Nicomachus has it, Pythagoras extended his tests to metal plates, he would have discovered this most basic fact about the role of percussion in the production of sound. Interestingly enough, there are reports of similar experiments conducted by a follower of Pythagoras, namely, Hippasus of Metapontum. Using four bronze discs of identical diameter, their thicknesses being of the proportions 4:3, 3:2, and 2:1, Hippasus is said to have produced the consonances fourth, fifth, and octave, respectively. What is more, when Hippasus presented his four well–tuned discs to Glaukos of Rhegium, a musician of great repute, Glaukos could actually play a tune on them. (Levin, 1994, pp. 92–93)

Huffman (1993) explicitly separates the impossible stories from Hippasus' sound science:

> The second crucial point is to recognize that the stories that the tradition tells about the discovery of the ratios for the most part describe observations that are impossible (e.g. the story that Pythagoras detected the concords in sounds he heard as he passed a smithy, which falsely presupposes that the pitch of sounds emitted from hammers as they strike the anvil is proportional to their weight). The only observation that is scientifically correct is the one assigned to Hippasus, that in the case of bronze disks of equal diameter the pitch emitted when struck will vary with their thickness. Thus, we can have some confidence that Hippasus (early fifth century) had knowledge of the ratios...
> (Huffman, 1993, pp. 147–148)

K. Ferguson (2011) also corroborates the physical accuracy of the disks:

> Hippasus, himself a contemporary of Pythagoras, made four bronze disks, all equal in diameter but of different thicknesses. The thickness of one "was 4/3 that of the second, 3/2 that of the third, and 2/1 that of the fourth". Hippasus suspended the disks to swing freely. Then he struck them, and the disks produced consonant intervals. This experiment is correct in terms of the physical principles involved, for the vibration frequency of a free–swinging disk is directly proportional to its thickness. Whoever designed and executed this experiment understood the basic harmonic ratios... (Ferguson, 2011, p. 69)

Lloyd (2014) highlights the distinction between the myths and the actual tests:

> ...none of those stories can be true, for the simple reason that they do not in fact produce the results claimed. The fact that no fewer than eight ancient authors repeat one or other version of these fictions is a shocking indication of the way one writer repeats what he has found in another quite uncritically... The two types of tests that could reveal the relations are those with bronze disks (associated with Hippasus) and with lengths of pipe or string. (Lloyd, 2014, p. 37)

Finally, Primavesi (2014) confirms the generation of perfect intervals:

> An experimental proof was provided, according to Aristotle's pupil Aristoxenus, already in the early fifth century B.C. by the Pythagorean Hippasus Of Metapontum: Hippasus is reported to have found out that striking four brazen discs of the same diameter the thicknesses of which stand in a proportion of 12: 9: 8: 6 brings about a "certain consonance" (*συμφωνίαν τινά*). Since "with free–swinging circular metal plates of the same diameter, the vibration frequencies are directly proportional to their thickness", such discs can indeed produce the perfect intervals... (Primavesi, 2014, p. 245)

**3.4.4.** *Conclusion*

Galilei's musical investigations do not seem to bear directly on Hippasus' experiments and construction. Euler's law definitely validates Hippasus' construction of the 4–chord. The fundamental role of Hippasus' 4–chord for the birth of anthyphairesis has not been noted by modern scholars.

**4.** *Hippasus' 4–chord and discovery of the concept of musical anthyphairesis*

**4.1.** *The acoustical experiments reveal the multiplicative nature of the operation of musical interval composition*

By establishing through the experiments of Lasus and Hippasus that fundamental intervals correspond to specific arithmetical ratios, it became possible to explore how these intervals combine. This leads us to the next crucial stage: the discovery that the composition of musical intervals is a multiplication between the ratios that correspond to these intervals.
We recall that the empirical rule of musical composition in a trichord with chords $a$, $b$, $c$, is the operation, which we denote by *, of generating the interval of the dichord $(a, c)$ from the musical intervals of the dichords $(a, b)$ and $(b, c)$, in symbols $(a, c) = (a, b) * (b, c)$. The fundamental application of the empirical rule of composition was no doubt the rule

octave = fifth * fourth.

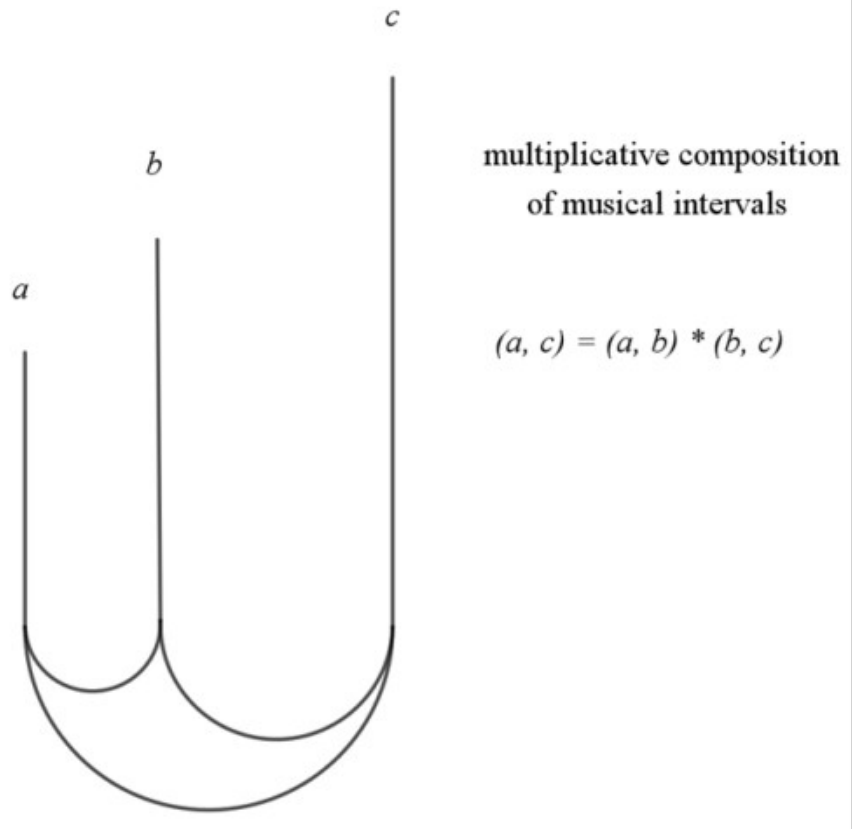

There is no doubt that, following the acoustical experiments and the arithmetization of musical intervals, it would be immediately understood that this fundamental rule could be realized by the trichord with chords of length 2, 3, 4.

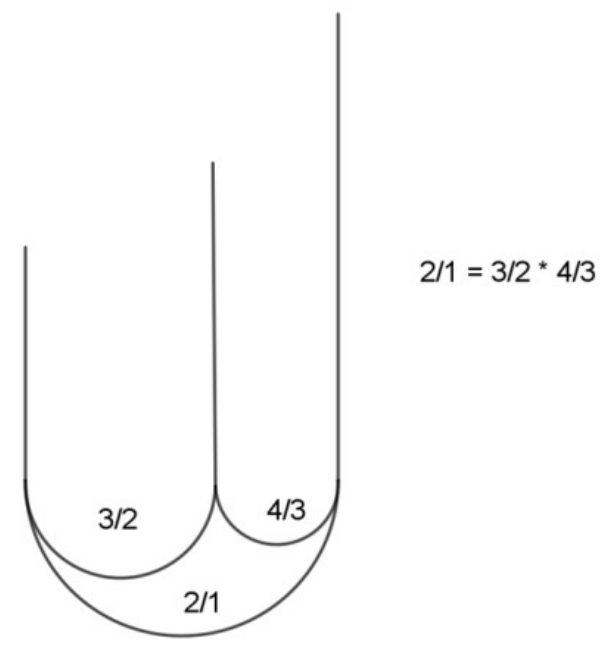

**4.2.** *The construction of Hippasus' 4–chord leads to the discovery of the concept of musical anthyphairesis*

How did Hippasus construct the 4–chord 6, 8, 9, 12? More specifically: how should the transition achieved by Hippasus from the arithmetization of the three basic musical intervals, achieved by Lasus and by himself, to the next step of constructing the 4–chord be conceived?
We recall that according to Koestler (1959), the performance of the acoustical experiments, resulting in the arithmetization of the three principal musical intervals, marked *the beginning of Science*. But the subsequent construction of the 4–chord with chords consisting of four cylinders of equal bases and heights 6, 8, 9, 12, reported by an *Anonymous Scholion to Plato's Phaedo,* whose eventual source is Aristoxenus, is even more momentous, as it marked *the birth of anthyphairesis*, in the musical/multiplicative form of musical intervals.
Indeed, Hippasus wanted to construct the simplest possible musical instrument able to produce some music. The trichord 2, 3, 4 consisting of the three basic intervals in the form 2/1 = 3/2 * 4/3 does not have such power. For this purpose, Hippasus divided the fifth by the fourth, producing a fourth interval, a theoretical *metaxuteta* between the fifth and the fourth, as Nicomachus had described it (in *Harmonicum enchiridion* 6.1, 17–18). Now it would not be necessary to make another acoustical experiment; it would be found by following the already discovered multiplicative nature of the composition: 3/2 = 4/3 * $x$ means 3/2 = 4/3 times $x$, thus $x$ = 9/8, thus the new musical interval, *the tone* is arithmetized by the ratio 9/8, and we have
3 /2 = 4/3 * 9/8.

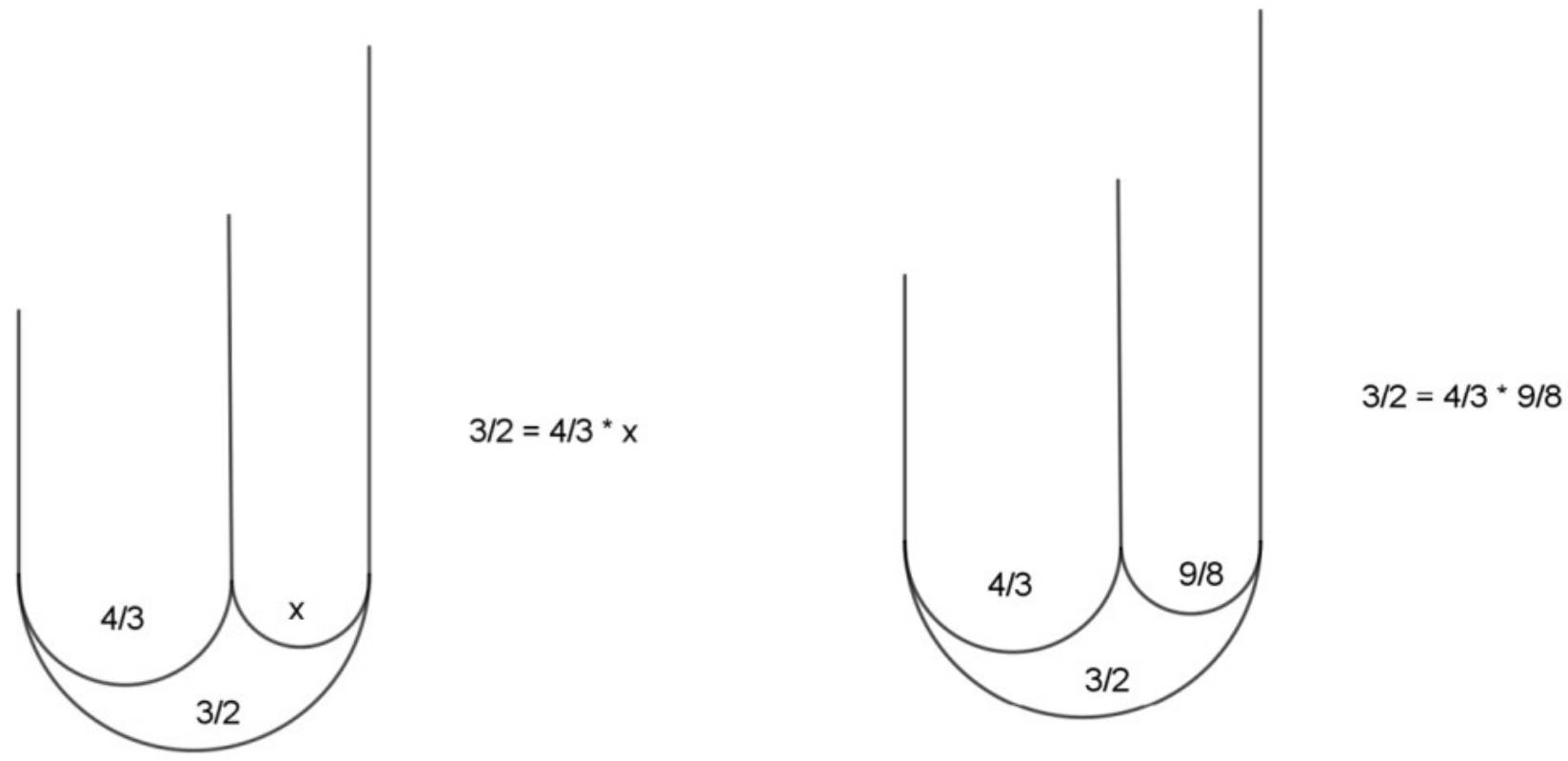

Thus, the construction of the 4–chord was achieved by *reversing* the role of the fifth:

in the first relation 2/1 = 3/2 * 4/3 the fifth is *active/dividing* with respect to the passive/divided octave,

while in the second relation 3/2 = 4/3 * 9/8 the same fifth is now *passive/divided*, with respect to the active fourth.

This reversal from **active** and **dividing** to **passive** and **being divided** is the essence of anthyphairesis. It was an easy matter to note that the least whole numbers realizing the 4–chord are the numbers 6, 8, 9, 12. Thus, Hippasus' construction of the 4–chord marked *the birth of the concept of anthyphairesis in Pythagorean music*.

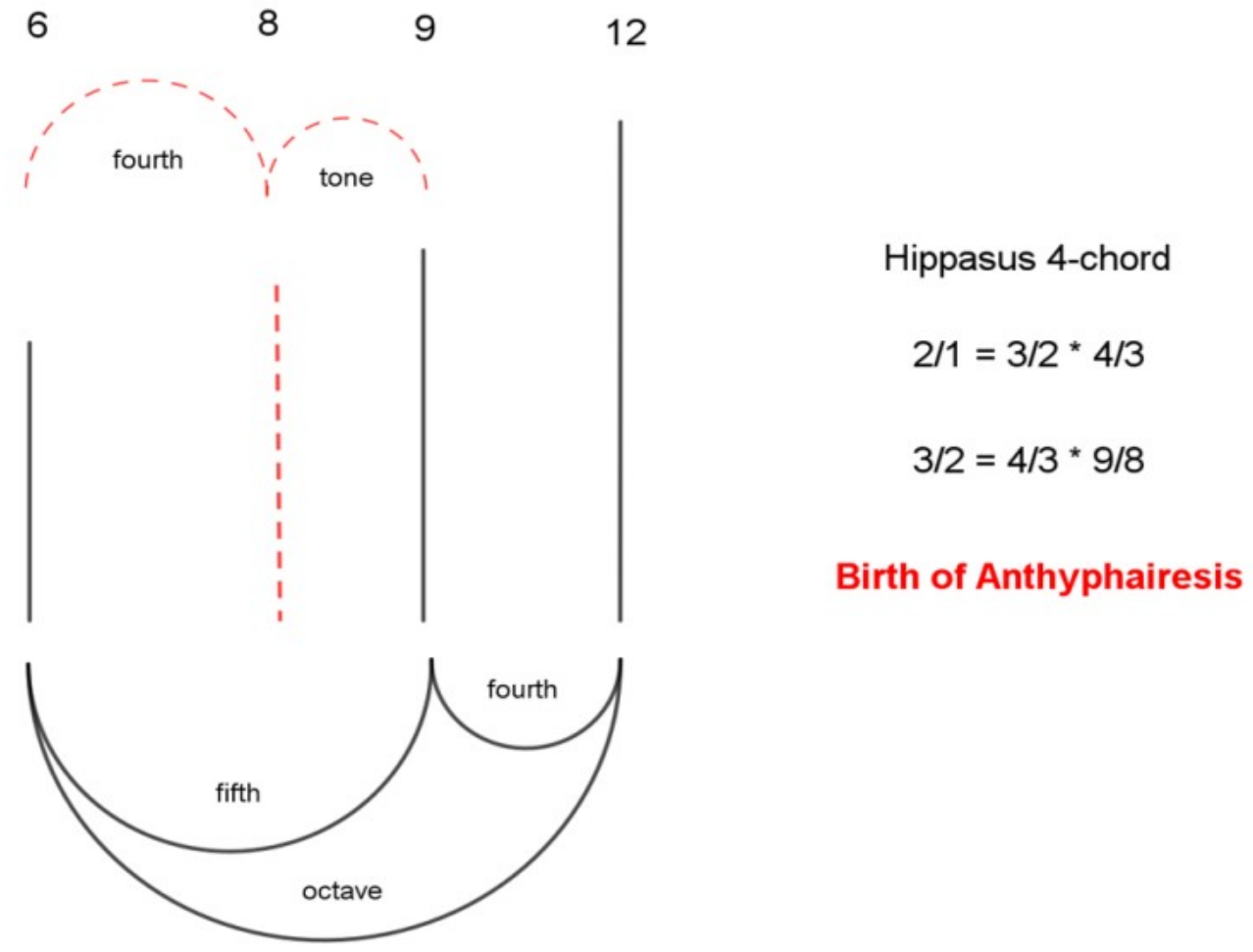

This attempt proved successful, since, according to the *Anonymous Scholion*, many years later, the famed musician Glaukos showed a proverbial skill in producing real music from Hippasus' 4–chord.

Hippasus' 4–chord consists of metallic cylinders, with identical bases and variable heights 6, 8, 9, 12, thus again a strictly one–dimensional linear variable, hence valid. It is crucial to emphasize that this 4–chord organically contains the acoustical experiments that preceded it. The mathematical structure of this 4–chord is governed by two fundamental relations, which both rightfully belong to Hippasus, even though they were famously preserved by later Pythagoreans (as will be discussed in Section 5 regarding Philolaus' *Fragment 6*).
The first relation (octave = fifth * fourth, or 2/1 = 3/2 * 4/3) represents the exact mathematical translation of the initial acoustical tests. In essence, the first relation *is* the experiments.
The second relation (fifth = fourth * tone, or 3/2 = 4/3 * 9/8) introduces the tone and marks the first anthyphairetic step, where the fifth reverses its role from dividing to being divided. The combination of these two relations completely defines and constructs the 4–chord. Therefore, Hippasus' 4–chord is the direct theoretical and physical culmination of his experimental discoveries.

**4.3.** *Babylonian origin of the Pythagorean 4–chord?*

According to Iamblichus in *Commentary to Nicomachus' Arithmetic (Nicomachi arithmeticam introductionem*) 118, 19–24, the musical proportion in four terms, *the musical 4–chord*, is a discovery of the Babylonians, and brought to Greece by Pythagoras.

> Τὰ νῦν δὲ περὶ τῆς τελειοτάτης ἀναλογίας ῥητέον
> ἐν τέσσαρσιν ὅροις ὑπαρχούσης, καὶ ἰδίως μουσικῆς ἐπικληθείσης,
> διὰ τὸ τοὺς μουσικοὺς λόγους τῶν καθ' ἁρμονίαν συμφωνιῶν
> τρανότατα ἐν αὐτῇ περιέχεσθαι.
> Εὕρημα δ' αὐτήν φασιν εἶναι Βαβυλωνίων, καὶ
> διὰ Πυθαγόρου πρώτου εἰς ῞Ελληνας ἐλθεῖν.
> Iamblichus, *Commentary to Nicomachus' Arithmetic 118,19–24* (1)
> Il faut maintenant parler de la proportion la plus parfaite, qui est en
> quatre termes et qu'on appelle en propre «musicale», du fait qu'elle
> renferme en elle, de la façon la plus claire, les rapports musicaux
> des consonances harmoniques. On dit que la découverte remonte
> aux Babyloniens et que Pythagore fut le premier à l'introduire en
> Grèce. [Iamblichus, trans. Vinel, 2014]

But Iamblichus, a little later, in 121, 13–16, follows Nicomachus' account and attributes the arithmetization of the 4–chord to Pythagoras himself.

> τὰς δὲ ἐπιτάσεις καὶ ἀνέσεις τῶν χορδῶν
> κατὰ τοὺς εἰρημένους λόγους γινομένας
> πρῶτον Πυθαγόραν ἱστοροῦσι συμμετρήσασθαι·
> παριόντα γὰρ εἴς τι χαλκοτυπεῖον…
> Iamblichus, *Commentary to Nicomachus' Arithmetic* 121, 13–16 (2)

(1) It is now necessary to speak of the most perfect proportion, consisting of four terms, which is specifically called "musical", because it most clearly contains the musical ratios of harmonic consonances. They say that it is a discovery of the Babylonians, and that it was first introduced to the Greeks by Pythagoras. (Translated by the authors from the Greek text)

(2) As for the tensions and relaxations of the strings produced according to the aforementioned ratios, it is recounted that Pythagoras was the first to harmonize them: as he was walking past a forge... (Translated by the authors from the Greek text)

Quant aux tensions et relâchements des cordes produits selon les rapports susdits, on raconte que Pythagore fut le premier à les harmoniser: comme il passait devant une forge… [Iamblichus, trans. Vinel, 2014]

We were not able to locate in Babylonian music something akin to the 4–chord that played a crucial role in Pythagorean music. The Babylonian origin of Hippasus' 4–chord is doubtful.

**4.4.** *Hippasus' 4–chord in Raphael's The School of Athens*

At this point, let us look at Raphael's painting, *The School of Athens*, created in 1509–11. Plato and Aristotle dominate the center. Pythagoras is positioned at the bottom left (with Euclid at the bottom right). Looking more closely, next to Pythagoras stands a young man holding a tablet and showing it to him. Looking even closer, the tablet held by the young man is precisely Hippasus' 4–chord. It is inscribed in Greek: *δια πασών* (octave), *δια πέντε* (fifth), *δια τεσσάρων* (fourth), a somewhat misspelled *επόγδοον* (whole tone), and the lengths of the strings are written in Roman numerals (with the number 8 incorrectly rendered as VII). It is highly improbable that Raphael intended the young man to represent Hippasus; he would have had to be aware of the anonymous commentary on Plato's dialogue *Phaedo* (cf. Section 3.2).

It is most likely that he had in mind Nicomachus' account, in which Pythagoras himself is presented as having discovered the numerical ratios of musical intervals through experiments and procedures that are, however, scientifically invalid. It is an exceptionally fortunate inspiration by Raphael, which – in hindsight – makes it appear as though he intended the young man to represent Hippasus, based on what we know today!

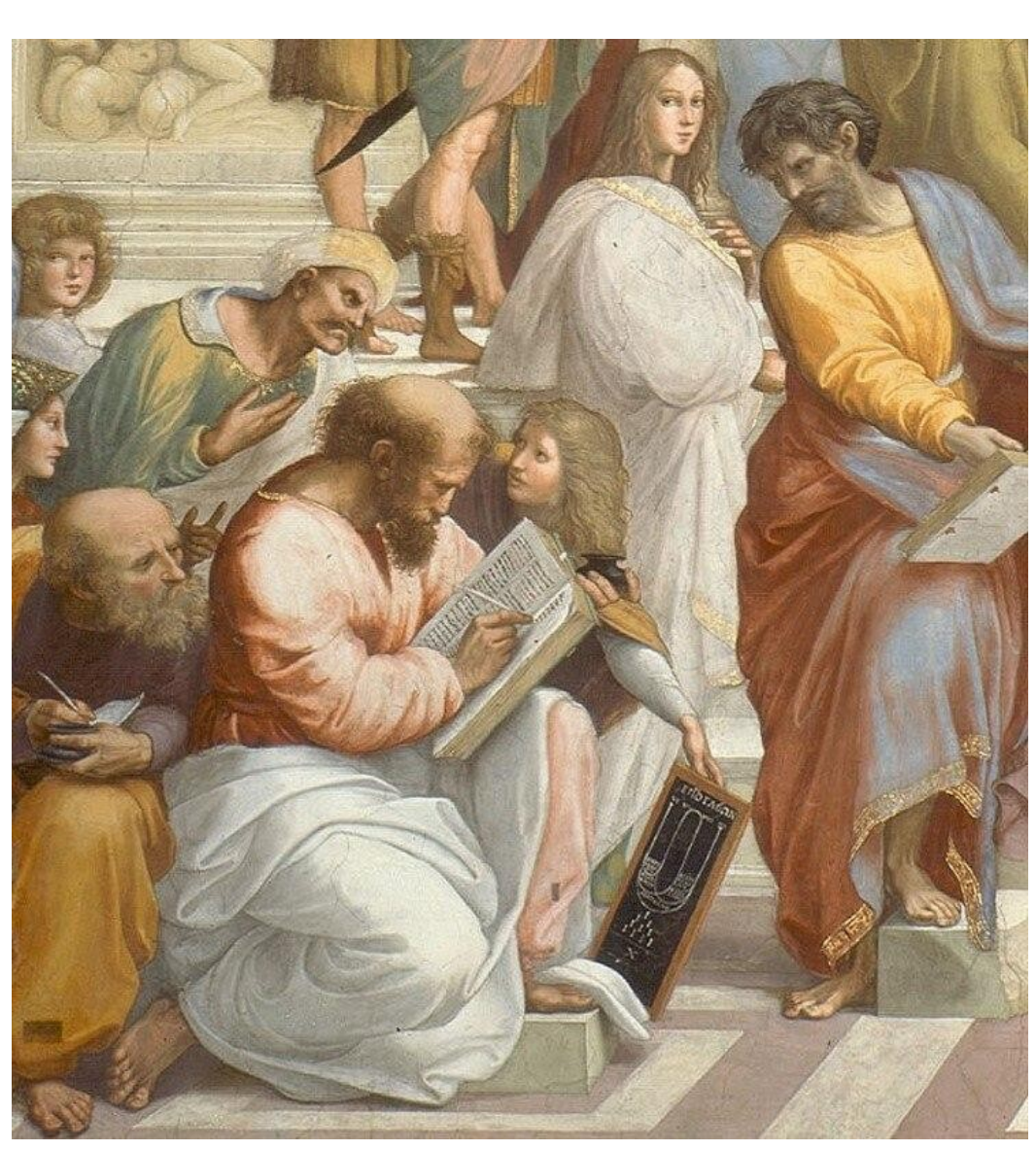

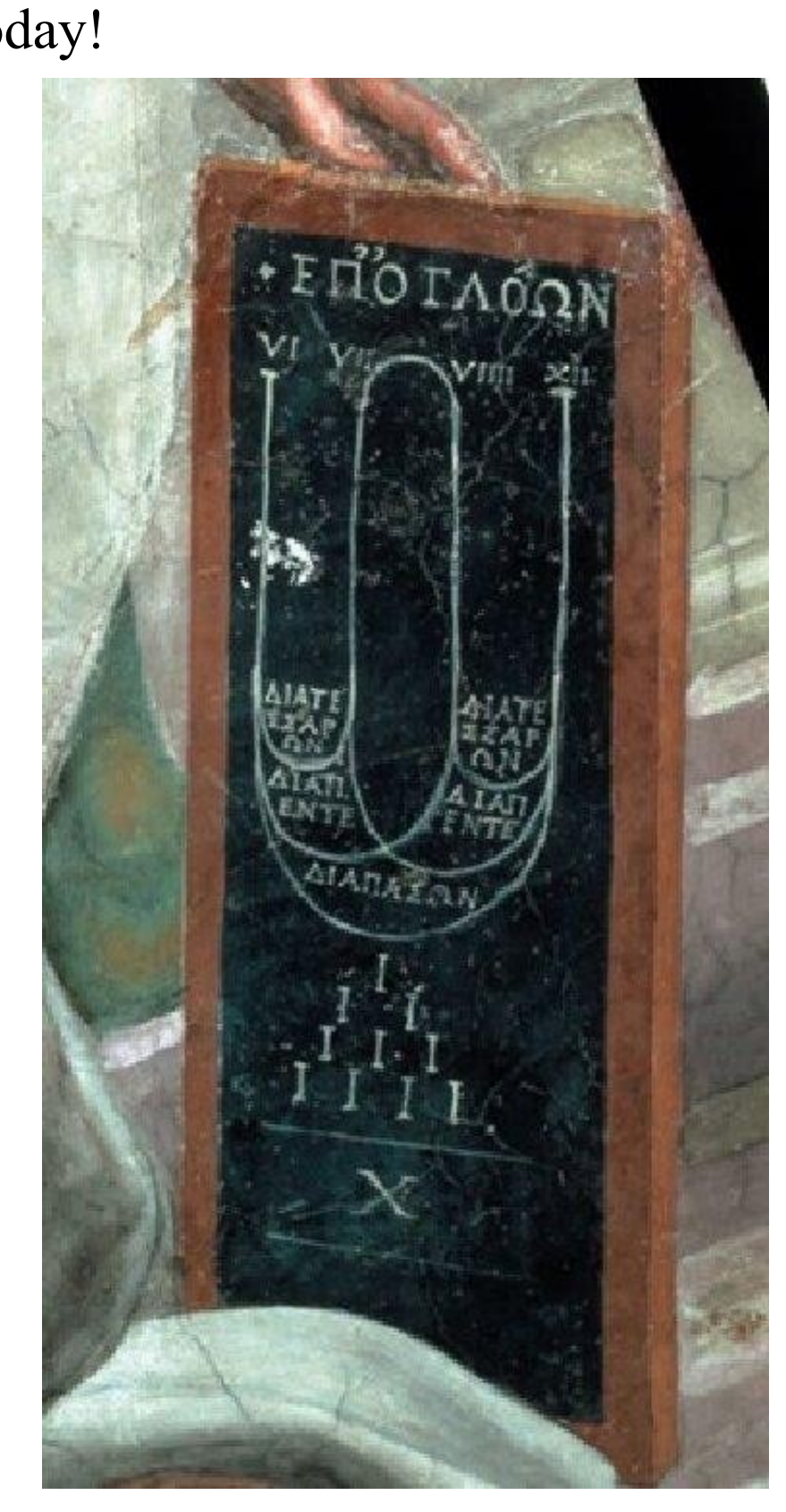

**5.** *Philolaus'* Fragment 6 *confirms this discovery*

The musical theory of the Pythagoreans, which arose from the original musical experiments and the construction of the 4–chord leading to the concept of musical anthyphairesis by Hippasus, was later extended in the important *Fragment 6* of Philolaus. Thus, Philolaus' *Fragment* 6 provides confirmation of the anthyphairetic nature of this theory.

**5.1.** *Philolaus, Fragment 6, 2– 15*

περὶ δὲ φύσιος καὶ ἁρμονίας ὧδε ἔχει·
ἁ μὲν ἐστὼ τῶν πραγμάτων ἀίδιος ἔσσα
καὶ αὐτὰ μὲν ἁ φύσις θείαν γα καὶ οὐκ ἀνθρωπίνην ἐνδέχεται γνῶσιν
πλέον γα ἢ ὅτι οὐχ οἷόν τ' ἦν οὐθὲν τῶν ἐόντων καὶ γιγνωσκόμενον ὑφ'
ἁμῶν γαγενέσθαι μὴ ὑπαρχούσας τᾶς ἐστοῦς τῶν πραγμάτων,
ἐξ ὧν συνέστα ὁ κόσμος, καὶ τῶν περαινόντων καὶ τῶν ἀπείρων.
ἐπεὶ δὲ ταὶ ἀρχαὶ ὑπᾶρχον οὐχ ὁμοῖαι οὐδ' ὁμόφυλοι ἔσσαι,
ἤδη ἀδύνατον ἦς κα αὐταῖς κοσμηθῆναι,
εἰ μὴ ἁρμονία ἐπεγένετο ὡιτινιῶν ἅδε τρόπωι ἐγένετο.
τὰ μὲν ὦν ὁμοῖα καὶ ὁμόφυλα ἁρμονίας οὐδὲν ἐπεδέοντο,
τὰ δὲ ἀνόμοια μηδὲ ὁμόφυλα μηδὲ ἰσοταγῆ
ἀνάγκα τᾶιτοι αύται ἁρμονίαι συγκεκλεῖσθαι,
οἵαι μέλλοντι ἐν κόσμωι κατέχεσθαι.
(Philolaus, *Fragment 6,* 2– 15)

Concerning *nature* and *harmonia,* it is like this.
The being of the things, which is eternal, and *nature* itself,
admit of divine and not human *knowledge,*
except that it was not possible
for any of the things that exist and are *known* by us
to have come into being,
if it were not for the existence of the beings of the things from which
the universe is composed of both *the finitizers and the infinites.*
And since there existed these principles, being neither alike nor of the same race,
it would then have been impossible for them to be organized together,
if *harmonia* had not come upon them, in whatever way it arose.
The things that were alike and of the same race had no need of harmonia as well, but things that were unlike and not of the same race nor equal in rank, for such things, it was necessary to have been *locked together by harmonia*, if they were to be held together in a cosmos.
(Huffman, 1993, pp. 123–124)

[1] ἁρμονίας δὲ μέγεθός ἐστι συλλαβὰ καὶ δι' ὀξειᾶν·
[2] τὸ δὲ δι' ὀξειᾶν μεῖζον τᾶς συλλαβᾶς ἐπογδόωι.
ἔστι γὰρ
ἀπὸ ὑπάτας ἐπὶ μέσσαν συλλαβά,
ἀπὸ δὲ μέσσας ἐπὶ νεάταν δι' ὀξειᾶν,
ἀπὸ δὲ νεάτας ἐς τρίταν συλλαβά,
ἀπὸ δὲ τρίτας ἐς ὑπάταν δι' ὀξειᾶν·
τὸ δ' ἐν μέσω μέσσας καὶ τρίτας ἐπόγδοον·
ἁ δὲ συλλαβὰ ἐπίτριτον, τὸ δὲ δι' ὀξειᾶν ἡμιόλιον, τὸ διὰ πασᾶν δὲ διπλόον.
(Philolaus, *Fragment* 6, 16–22)

[1] The magnitude of *harmonia* is *syllaba* and *di' oxeian.*
Thus, 2/1 = 3/2 * 4/3, with 3/2 > 4/3.
The first relation of the intervals in the 4–chord, octave = fifth * fourth, when arithmetized takes the form 2/1 = 3/2 * 4/3. Thus, the first step is a consequence of the acoustical experiments by Lasus and Hippasus.

[2] The *di' oxeian* is greater than the *syllaba* in the *epogdoic* ratio.
fifth = fourth * tone, 3/2 = 4/3 * 9/8.
The second step is already seen to be the crucial step that reveals the anthyphairetic nature of the 4–chord. The musical interval fifth, *active* and *dividing* under the operation of composition in the first relation is turned into *passive* and *being divided* under the operation of composition in the second relation. Since composition is really a homonymous multiplication of ratios, the 4–chord comprehends the first two steps of the musical/compositional/multiplicative *anthyphairesis* of the initial octave 2/1 to fifth 3 /2. Thus, 3/2 = 4/3 * 9/8, with 4/3 > 9/8.
In [1] and [2], we recognize the relations by means of which the 4–chord is constructed.
Philolaus proceeds immediately to describe Hippasus' 4–chord, thus directly correlating this construction with the relations [1] and [2].

From hypate to mese is a syllaba,
from mese to nete is a di'oxeian,
from nete to trite is a syllaba, and
from trite to hypate is a di' oxeian.
The interval between trite and mese is epogdoic,
the syllaba is epitritic, the di'oxeian hemiolic, and the diapason is duple. (Huffman, 1993, pp. 146–147)

[1] The first relation octave 2/1 = fifth 3/2 * fourth 4/3 produces the first chord

<table>
<tr><td>Hypate</td><td>Mese</td><td></td><td>Nete</td></tr>
<tr><td colspan="2">fourth 4/3</td><td colspan="2"></td></tr>
<tr><td></td><td colspan="3">fifth 3/2</td></tr>
</table>

[2] The second relation fifth 3/2 = fourth 4/3 * tone 9/8 produces *Hippasus' tetrachord*

<table>
<tr><td>Hypate<br>6</td><td>Mese<br>8</td><td>Trite<br>9</td><td>Nete<br>12</td></tr>
<tr><td colspan="2">fourth 4/3</td><td colspan="2"></td></tr>
<tr><td></td><td colspan="3">fifth 3/2</td></tr>
<tr><td></td><td colspan="2">tone 9/8</td><td></td></tr>
<tr><td colspan="2"></td><td colspan="2">fourth 4/3</td></tr>
</table>

But we realize that relations [1] and [2] are anthyphairetic.

**5.2.** *The construction of the octachord lyra from the first three steps of the musical anthyphairesis*

> Thus, the only reasonably solid conclusion that can be drawn is that in the generation before Philolaus, Hippasus, at least, knew of the ratios that corresponded to the concordant intervals of the octave, fourth, and fifth, but there is no evidence of knowledge of the ratios that correspond to the tone and the "remainder" (diesis or leimma) which are used to fill out the rest of the diatonic scale.
> (Huffman, 1993, p. 148)

In order to confirm both the anthyphairetic nature of Philolaus *Fragment 6* and its relevance for the construction of the 8–chord lyra and its intermediates, we consider the next part of *Fragment 6.*

> οὕτως
> [3b] ἁρμονία πέντε ἐπόγδοα καὶ δύο διέσιες,
> [3a] δι' ὀξειᾶν δὲ τρία ἐπόγδοα καὶ δίεσις,
> [3] συλλαβὰ δὲ δύ' ἐπόγδοα καὶ δίεσις.
> (Philolaus, *Fragment 6*, 22–24)

> Thus,
> [3] *syllaba* is two *epogdoics* and a *diesis;*
> [3a] *di'oxeian* is three epogdoics and a diesis; and
> [3b] *harmonia* consists of five epogdoics and two dieses.
> (Huffman, 1993, pp. 123–124)

According to Boethius' *De Institutione Musica* III.5:

> Philolaus, a Pythagorean, tried to divide the tone in another manner,
> postulating that the tone had its origin in the number that constitutes the first cube of the first odd number, for that number was highly revered among the Pythagoreans.
> Since 3 is the first odd number,
> if you multiply 3 by 3, then this by 3, 27 necessarily arises,
> which stands at the distance of a tone from the number 24,
> the same 3 being the difference.
> For 3 is an eighth part of the quantity 24, and, added to the same, it gives the first cube of 3, 27.
> From this number, 27, Philolaus made two parts,
> one that is more than half, which he called the *"apotome"*, and
> the remainder, which is less than half, which he called the *"diesis"*.
> (The diesis later came to be called the *"minor semitone"*.)
> The difference between these he called the *"comma"*.
> (Boethius, trans. Bower, 1989, p. 96)

Thus, Philolaus divided first the tone into a diesis and an apotome; in numbers, 9/8 = 256/243 * 32/27.
This division explains well the construction of a 6–chord lyra.

According to Boethius, *De Institutione Musica* III.8:

> Philolaus igitur haec atque his minora spatialibus definitionibus includit: [3] Diesis, inquit, est spatium, quo maior est sesquitertia proportio duobus tonis. [4] Comma vero est spatium, quo maior est sesquioctava proportio duabus diesibus.
> Philolaus incorporates these and intervals smaller than these in the following definitions.
> [3] The *diesis* is the interval by which a *sesquitertian* [= fourth] ratio is larger than two *tones*.
> fourth = two tones plus a diesis
> [4] The *comma* is the interval by which the *sesquioctave* ratio [= tone] is larger than two *diesis*
> tone = two dieses plus a comma.
> (Boethius, trans. Bower, 1989, p. 97)

According to Burkert:

> Philolaus makes two parts, one, which is larger than half, and which he calls apotome, and another, which is smaller than half, which he calls diesis; later it was called "the smaller semitone". The difference between these parts he calls [comma].
> (Burkert, 1972, p. 395)

So when the Pythagoreans realized the third step of their theory of music, reported in the important Philolaus' *Fragment 6*, namely the equation in musical intervals fourth = two tones * diesis, in ratios of numbers
$4/3 = (9/8) * (32/27) = (9/8) * (9/8) * (256/243)$,
it was clear on the one hand, that the third step is clearly anthyphairetic under the composition of intervals/ homonymous multiplication of ratios, and, on the other hand, that the third step produces the intermediate 7–chord octave and mainly the final 8–chord octave; thus:

octave 2/1 (dichord)

<table>
<tr><td>H</td><td></td><td></td><td></td><td></td><td></td><td></td><td>N</td></tr>
<tr><td colspan="8">octave 2/1</td></tr>
</table>

[1] First anthyphairetic relation
octave 2/1 = fifth 3/2 * fourth 4/3 (3–*chord*)

<table>
<tr><td>H</td><td></td><td></td><td>M</td><td></td><td></td><td></td><td>N</td></tr>
<tr><td colspan="4">fourth 4/3</td><td></td><td></td><td></td><td></td></tr>
<tr><td></td><td></td><td></td><td colspan="5">fifth 3/2</td></tr>
</table>

[2] The second anthyphairetic relation fifth 3/2 = fourth 4/3 * tone 9/8
produces Hippasus' *4–chord*

<table>
<tr><td>H</td><td></td><td></td><td>M</td><td>T</td><td></td><td></td><td>N</td></tr>
<tr><td colspan="4">fourth 4/3</td><td></td><td></td><td></td><td></td></tr>
<tr><td></td><td></td><td></td><td colspan="2">tone 9/8</td><td></td><td></td><td></td></tr>
<tr><td></td><td></td><td></td><td></td><td colspan="4">fourth 4/3</td></tr>
</table>

[3a] The initial form of the third anthyphairetic relation
fourth 4/3 = tone 9/8 * apotome 32/27 produces the *6–chord*

<table>
<tr><td>H</td><td>PH</td><td></td><td>M</td><td>T</td><td>PN</td><td></td><td>N</td></tr>
<tr><td colspan="2">tone 9/8</td><td></td><td></td><td></td><td></td><td></td><td></td></tr>
<tr><td></td><td colspan="3">32/27</td><td></td><td></td><td></td><td></td></tr>
<tr><td></td><td></td><td></td><td colspan="2">tone 9/8</td><td></td><td></td><td></td></tr>
<tr><td></td><td></td><td></td><td></td><td colspan="2">tone 9/8</td><td></td><td></td></tr>
<tr><td></td><td></td><td></td><td></td><td></td><td colspan="3">32/27</td></tr>
</table>

[3b] An intermediate form of the third anthyphairetic relation, involving both the interval 32/27 of the initial form and the diesis of the final form of the third relation

$4/3 = 9/8 * (32/27) = (9/8)^2 * (256/243)$, produces the *7–chord* [Philolaus]

| H | PH | L | M | T | PN | | N |
|---|---|---|---|---|---|---|---|
| tone 9/8 | | | | | | | |
| | tone 9/8 | | | | | | |
| | | diesis 256/243 | | | | | |
| | | | tone 9/8 | | | | |
| | | | | tone 9/8 | | | |
| | | | | | 32/27 | | |

[3c] The final form of the third anthyphairetic relation

fourth $4/3$ = two tones $(9/8)^2$ * diesis $(256/243)$, final stage produces the *8–chord*

| H | PH | L | M | PM | T | PN | N |
|---|---|---|---|---|---|---|---|
| tone | | | | | | | |
| | tone | | | | | | |
| | | diesis | | | | | |
| | | | tone | | | | |
| | | | | tone | | | |
| | | | | | tone | | |
| | | | | | | diesis | |

[H hypate, PH Parhypate, L Lichanos, M Mese, PM Paramese, T Trite, PN Paranete, N Nete]

| *t* | *t* | *d* | *t* | *t* | *t* | *d* |
|---|---|---|---|---|---|---|

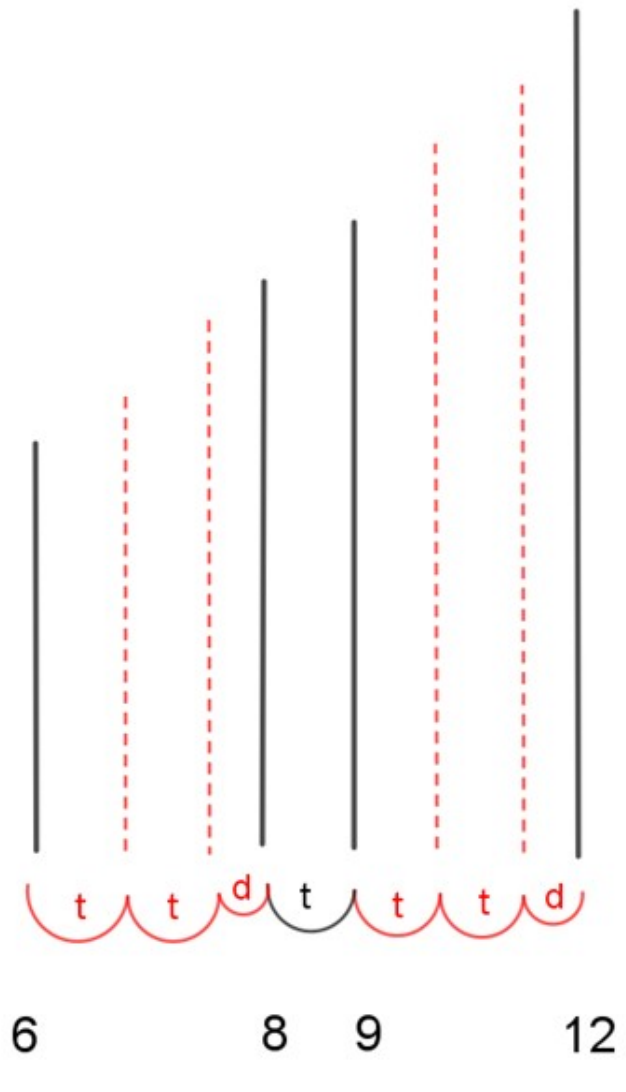

Generation of the Pythagorean 8-chord from the third Philolaus anthyphairetic relation, applied to the Hippasus 4-chord.

The smallest numbers for the length of the chords of an octachord lyra.

| H | PH | L | M | PM | T | PN | N |
|---|---|---|---|---|---|---|---|
| $2^7\,3$ | $2^4\,3^3$ | $2\,3^5$ | $2^9$ | $2^6\,3^2$ | $2^3\,3^4$ | $3^6$ | $2^8\,3$ |
| 384 | 432 | 486 | 512 | 576 | 648 | 729 | 768 |

<table>
<tr><td colspan="2">tone<br>$3 \cdot 2^4 \cdot 9/3 \cdot 2^4 \cdot 8 = 9/8$</td><td></td><td></td><td></td><td></td><td></td><td></td></tr>
<tr><td></td><td colspan="2">tone<br>$3^3 \cdot 2 \cdot 9/3^3 \cdot 2 \cdot 8 = 9/8$</td><td></td><td></td><td></td><td></td><td></td></tr>
<tr><td></td><td></td><td colspan="2">diesis<br>$2 \cdot 256/2 \cdot 243 = 256/243$</td><td></td><td></td><td></td><td></td></tr>
<tr><td colspan="4">fourth<br>$512/384 = 2^9/2^7 3 = 2^7 \cdot 4/2^7 \cdot 3 = 128 \cdot 4/128 \cdot 3 = 4/3$</td><td></td><td></td><td></td><td></td></tr>
</table>

<table>
<tr><td></td><td></td><td></td><td colspan="2">tone<br>$2^6 \cdot 9/2^6 \cdot 8 = 9/8$</td><td></td><td></td><td></td></tr>
<tr><td colspan="5">fifth<br>$576/384 = 2^6 3^2/2^7 3 = (2^6 \cdot 3) \cdot 3/(2^6 \cdot 3) \cdot 2 = 192 \cdot 3/192 \cdot 2 = 3/2$</td><td></td><td></td><td></td></tr>
</table>

<table>
<tr><td></td><td></td><td></td><td></td><td colspan="2">tone<br>$3^2 \cdot 2^3 \cdot 9/3^2 \cdot 2^3 \cdot 8 = 9/8$</td><td></td><td></td></tr>
<tr><td></td><td></td><td></td><td></td><td></td><td colspan="2">tone<br>$3^4 \cdot 9/3^4 \cdot 8 = 9/8$</td><td></td></tr>
<tr><td></td><td></td><td></td><td></td><td></td><td></td><td colspan="2">diesis<br>$3 \cdot 256/3 \cdot 243 = 256/243$</td></tr>
<tr><td></td><td></td><td></td><td></td><td colspan="4">fourth<br>$2^8 3/2^6 3^2 = (2^6 \cdot 3) \cdot 4/(2^6 \cdot 3) \cdot 3 = 4/3$</td></tr>
</table>

<table>
<tr><td></td><td></td><td></td><td colspan="5">fifth<br>$2^8 3/2^9 = 2^8 \cdot 3/2^8 \cdot 2 = 256 \cdot 3/256 \cdot 2 = 3/2$,</td></tr>
</table>

<table>
<tr><td colspan="8">octave<br>$768/384 = 2^8 3/2^7 3 = (2^7 \cdot 3) \cdot 2/(2^7 \cdot 3) \cdot 1 = 384 \cdot 2/384 \cdot 1 = 2/1$</td></tr>
</table>

From the fact that one chord, Mese, is represented as a power of 2 ($2^9$), and another chord, Paranete, is represented as a power of 3 ($3^6$), it follows that the numbers for the strings are the smallest possible.
The fourth step of the musical anthyphairesis is:
tone = two dieses * Pythagorean comma. Denoting the Pythagorean comma by *c*, we note that this fourth step divides the tone into two almost equal intervals, namely into the diesis *d* and into the composition *dc* of the diesis with the comma. From this division, we obtain a further division of the octave by dividing each tone of the octachord into *d* and *dc*, resulting in a *12 –tone almost equal tempered scale:*

| d | dc | d | dc | d | d | dc | d | dc | d | dc | d |
|---|---|---|---|---|---|---|---|---|---|---|---|

**5.3.** *Burkert (1972) on the musical anthyphairesis*

Burkert's (1972) dismissal of these proportions as "mathematical absurdities" or "pure frivolity" stems from a modern arithmetic bias. By failing to recognize the inherently anthyphairetic nature of Philolaus' divisions, he overlooks their mathematical and structural necessity:

> We still have to consider the "mathematical absurdities" reported by Boethius…. Philolaus divides the whole tone into two unequal parts, *diesis* 256:243, and *apotome*; the difference between the two is called *comma.* According to this: The *apotome* would be 2187:2048, and the c*omma* 531441:524288 – pure frivolity. (Burkert, 1972, p. 394–395)

and also

> The musical *apotome* has still less to do with irrationality, either when correctly calculated or in the incorrect version of Philolaus. (Burkert, 1972, p. 395)

Thus, the persistent (symbolic or literal) connection of incommensurability with Infinity and infinite anthyphairesis, which appears both in the narratives of Pappus (1930), Plutarch (1967), Iamblichus (1989), as well as in the passage of Proclus (*In Euclidem 6,19–23;* trans. Morrow, 1992), weakens significantly the possibility that the reconstruction of the discovery of the proof of the incommensurability of the diameter to the side of a square would be realized with a variant of the following form of indirect argument, which Aristotle knew and often mentioned (cf. *Analytica Priora* 41a26–31), and which appears to be a later, but still ancient addition to Book X of the *Elements* (as Proposition X.117).

**5.4.** *The infinity of the musical anthyphairesis*

The significant conclusion of Philolaus' *Fragment* 6, lines 2–24 together with the addition, preserved only by Boethius, *De Institutione Musica* III 8, but attributed to Philolaus, is that the *Pythagorean theory of music* is based on a musical/ multiplicative anthyphairesis, corresponding to the operation of synthesis/ composition of musical intervals, denoted *, of the two basic musical intervals harmonia 2/1 and fifth 3/2, of which the first four steps are the following:

[1] $2/1 = 3/2 * 4/3$, with $3/2 > 4/3$
(harmonia/octave = di' oxeian/fifth * syllaba/fourth),
[2] $3/2 = 4/3 * 9/8$, with $4/3 > 9/8$
(di' oxeian = syllaba * epogdoon),
[3] $4/3 = (9/8)^2 * 256/243$, with $9/8 > 256/243$
(syllaba = two epogdoa * diesis),
[4] $9/8 = (256/243)^2 * 3^{12}/2^{19}$, with $256/243 > 3^{12}/2^{19}$
(epogdoon = two dieseis * Pythagorean comma)

*Notes.* The Pythagorean comma is examined in detail by Proclus (cf. trans. Taylor, 1820) in his *Commentary to Plato's Timaeus* 2,174, 11–190,31. The anthyphairetic nature of the Pythagorean theory is clearly noted by Fowler (1999, Section 4.5, p. 148).

*Proposition.* The Pythagorean musical anthyphairesis of the octave (2/1) to the fifth (3/2) is infinite.
*Proof.* The process follows a sequence of successive divisions where each musical interval is divided by the one preceding it. By mathematical induction, it is shown that every remainder in this sequence is a ratio of the form alternately a power of 2 to a power of 3 or conversely a power of 3 to a power of 2. The process of this anthyphairesis, since it is multiplicative, terminates only if a remainder becomes equal to the ratio 1/1, something clearly impossible. So, the anthyphairesis is infinite.

**5.5.** *Question. Did the Pythagoreans know that the musical anthyphairesis is infinite?*

The infinity of the anthyphairesis of harmony is not proved or explicitly stated in Philolaus' *Fragment 6*, and, therefore, we cannot be sure that the Pythagoreans had proved it. But the general sense of *Fragment 6* is that harmony brings "kosmos" (i.e., order), and that "kosmos" consists of infinites and finitizers. The only form of infinity that we can conceive for the Pythagorean harmony is the anthyphairetic one; in this sense, *Fragment 6* itself can be seen as a formulation of the infinity of harmony's anthyphairesis.
The following is a consequence of the infinity of the musical anthyphairesis:

*Corollary.* There is no common measure of the octave and the fifth, of the fifth and the fourth, of tone and diesis.
*Proof.* Analogous to the proof of Proposition X.2.
Thus, the presence of two measures in the Pythagorean scale appears to be indirectly connected with the infinity of the musical anthyphairesis.
In fact, it is reasonable to conjecture that the infinity of the multiplicative anthyphairesis of the octave vs. the fifth was known to the Pythagoreans, and that, since the Pythagoreans adopted the diatonic Pythagorean musical scale, tone tone diesis tone tone tone diesis, i.e., a musical scale with two measures (tone and diesis), and not by a common measure, as they would certainly prefer and as they hoped they might find, they were admitting the impossibility of finding and having a musical scale with a common measure.

But in fact, Aristotle, in *Metaphysics* 1053a12–18, refers to the necessity for two musical measures, correlating the musical incommensurability with the geometrical incommensurability of the diameter to the side of a square and with incommensurability in general.

> καὶ *ἐν μουσικῇ δίεσις*, ὅτι ἐλάχιστον [*μετρον*],
> καὶ *ἐν φωνῇ* στοιχεῖον.
> καὶ ταῦτα πάντα ἕν τι οὕτως,
> οὐχ ὡς κοινόν τι τὸ ἓν
> ἀλλ' ὥσπερ εἴρηται.
> οὐκ ἀεὶ δὲ τῷ ἀριθμῷ *ἓν τὸ μέτρον*
> ἀλλ' ἐνίοτε πλείω, οἷον
> *αἱ διέσεις δύο*, αἱ μὴ κατὰ τὴν ἀκοὴν ἀλλ' ἐν τοῖς λόγοις,
> καὶ αἱ φωναὶ πλείους αἷς μετροῦμεν,
> καὶ *ἡ διάμετρος δυσὶ μετρεῖται* καὶ ἡ *πλευρά*,
> καὶ τὰ *μεγέθη πάντα*.
> Aristotle, *Metaphysics* 1053a12–18

> Hence too in music, the diesis is the smallest [measure] and in speech, the letter. Not as common [measure] to them all, but as a way of speaking. However, arithmetically, we do not always have one measure, but sometimes more than one measure.
> *Diesis,* for instance, *is two,* as established not by hearing but through
> their ratios, and there is *more than one* sound by which speech is measured, *and the diagonal, too, and the side* of the square are measured by two magnitudes, *and* the same for all [incommensurable] magnitudes.
> [Aristotle, 1924. *Aristotle's Metaphysics*. Translation by W. D. Ross, Oxford: Clarendon Press. With modifications by the authors]

Aristotle explicitly confirms that the impossibility of finding a single common measure in the musical scale (hence the need for both tone and diesis) is the exact empirical manifestation of geometric incommensurability.

It is therefore reasonable to interpret the Pythagoreans' adoption of a musical scale with two measures (tone, diesis), instead of the single measure they would undoubtedly have preferred, as a strong indication that at some point, certainly at the time of Philolaus, they became aware of the infinity of multiplicative anthyphairesis of harmony and incommensurability.

**6.** *The arithmetized Pythagorean theory of music, following the musical experiments, the construction of the 4–chord, and the discovery of anthyphairesis*

**6.1.** *The empirical pre–Pythagorean music*

According to Sections 1, 3, and 4, we may regard empirical music as the system consisting of an inexact, non–mathematical definition of the term musical interval, and the three properties as Axioms/postulates:

the empirical/acoustical *definitions* of the four musical intervals
*octave, fifth, fourth,* and *tone*, and
the class of *dichords* generating each of these musical intervals; and
the three empirical *properties* to which these classes are subject:
*musical transitivity,*
*musical most equality, and the*
*commutativity of musical composition*.

We will now describe the way in which the Pythagoreans transformed the heretofore empirical art of music into a rigorous mathematical system by introducing a definition of the musical intervals.
To understand this transformation clearly, we consider the ancient empirical art of music as consisting of empirical rules without a clear definition or proof. The epoch–making experiments by Lasus and Hippasus, the discovery of musical anthyphairesis, and the construction of the 4–chord by Hippasus make it possible to introduce exact arithmetical definitions for the first three musical intervals.

**6.2.** *The definition of the four musical intervals*

*Definition*
The *octave* is the class of all dichords/dyads of line segments ($a$, $b$), $a > b$, such that $a = 2b$.

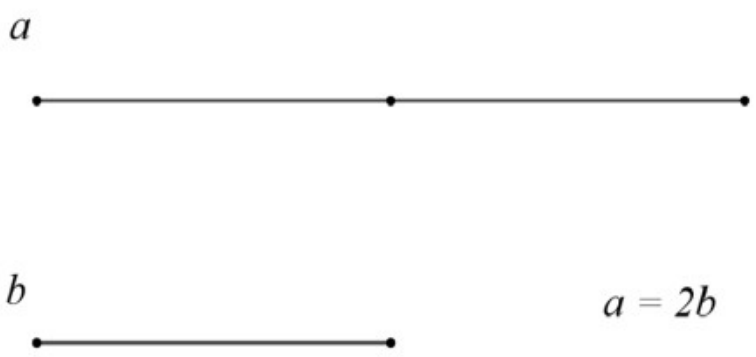

The *fifth* is the class of all dichords/dyads of line segments ($a$, $b$), $a > b$, such that $a = 3(a – b)$, $b = 2(a – b)$.

$a$ $c = a - b$ $a = 3(a-b)$

$b$ $b = 2(a-b)$

c

The *fourth* is the class of all dichords/dyads of line segments ($a$, $b$), $a > b$, such that $a = 4(a - b)$, $b = 3(a - b)$.
The *tone* is the class of all dichords/dyads of line segments ($a$, $b$), $a > b$, such that $a = 9(a - b)$, $b = 8(a - b)$.
The representation of a dichord as two unequal *line segments* helps in arriving at the definition.

**6.3.** *The transfer from the musical to the arithmetical anthyphairesis*

This is precisely the point where the concept of *musical anthyphairesis*, already discovered by the construction of the 4–chord by Hippasus, is naturally *transferred* to the concept of *arithmetical anthyphairesis.* In fact, as in the case of the musical anthyphairesis and as further suggested by the representation of numbers as chords, it would be natural to divide the greater number $a$ by the smaller $b$, for the case of the fifth, to obtain the arithmetical anthyphairetic relation:
$a = b + c$, $c < b$,
and then to reverse the role of the smaller $b$ in relation to the still smaller $c$, and divide $b$ by $c$:
$b = 2c$.
Thus, a musical interval is a fifth, (equivalently a ratio, $a$ to $b$ is hemiolic) if it satisfies the anthyphairetic description
$a = b + c$, $c < b$,
$b = 2c$.
This is explicitly given by Theon Smyrneus 78,1–5, as we will discuss in detail in Section 11.
Analogous steps are given for the fourth and for the tone.

*Proposition*
A dyad $a > b$ is a *fifth* if and only if, setting $c = a - b$, we have $a = b + c$, $b = 2c$ (Namely, *Anth* ($a$, $b$) = [1, 2]).
A dyad $a > b$ is a *fourth* if and only if, setting $c = a - b$, we have $a = b + c$, $b = 3c$ (Namely, *Anth* ($a$, $b$) = [1, 3]),
A dyad $a > b$ is a *tone* if and only if, setting $c = a - b$, we have $a = b + c$, $b = 8c$ (Namely, *Anth* ($a$, $b$) = [1, 8]).
Thus, a dichord $a > b$ is in the fifth, fourth, or tone if and only if the anthyphairesis of $a$ to $b$ has length 2 and sequence of quotients the dyad 1, 2, or 1, 3, or 1, 8, respectively.

**6.4.** *The arithmetical definition of musical intervals*

The arithmetical definition of the musical intervals has the power to prove the three empirical properties of musical intervals rigorously, but it is further able to prove two more properties of the musical intervals.

*Proposition (Musical Transitivity)*
If $(a, b)$, $(c, d)$ are in the same musical interval, and $(c, d)$, $(e, f)$ are in the same musical interval, then $(a, b)$, $(e, f)$ are also in the same musical interval.
*Proof.* We must verify each of the four cases separately.
E.g. let $a > b$ and $c > d$ be in the fourth:
$(a = 4(a – b), b = 3(a – b)$, and $c = 4(c – d), d = 3(c – d))$,
and $c > d$ and $e > f$ be in the fourth:
$(c = 4(c – d), d = 3(c – d)$ and $e = 4(e – f), f = 3(e – f))$.
Then $a > b$ and $e > f$ are in the fourth.

*Proposition (Musical Most Equality)*
If $(a, b)$, $(a, b')$ are dichords in the same musical interval, then $b = b'$.
*Proof.* Again, we must verify each of the four cases separately.
E.g., let $(a, b)$, $(a, b')$ be in the fifth.
Then, by the definition, $a = 3(a – b)$, $b = 2(a – b)$, and
$a = 3(a – b')$, $b' = 2(a – b')$;
then $3(a – b) = 3(a – b')$, then $b = b'$.

*Proposition (Synthesis, Diairesis)*
(a) (*Synthesis of equivalent dichords*)
If $(a, b)$ and $(c, d)$ are equivalent dichords, then $(a + c, b + d)$ and $(a, b)$ are equivalent dichords.
In particular, $(a, b)$ and $(ka, kb)$ are equivalent dichords for every natural number $k$.
(b) (*Diairesis of equivalent dichords*)
If $(a, b)$, and $(c, d)$, are equivalent dichords, and $a > b$, $c > d$, then
$(a – b, c – d)$ and $(a, b)$ are equivalent dichords.
*Proof.* E.g., for (a), let $(a, b)$, and $(c, d)$ be in the fifth.
$a = 3(a – b)$, $b = 2(a – b)$, and $c = 3(c – d)$, $d = 2(c – d)$.
Then $a + c = 3(a + c – b – d)$, and $b + d = 2(a + c – b – d)$.
Then, by definition, $(a + c, b + d)$ is in the fifth.

*Proposition (Commutativity of Musical Composition)*
If $a$, $b$, $c$, and $d$, $e$, $f$ are 3–chords, such that $(a, b)$, $(e, f)$ are equivalent, $(b, c)$, $(d, e)$ are equivalent, and $(a, c)$, $(d, f)$ are dichords in some musical interval, then $(a, c)$, $(d, f)$ are equivalent.
*Proof.* It can be seen that there are essentially only two possible cases.

*Case 1*. $(a, c)$ is in the octave, $(a, b)$ in the fifth, $(b, c)$ in the fourth.
For *Case 1*, by definition $a = 2(b - a)$, $b = 3(b - a)$, $b = 3(c - b)$, $c = 4(c - b)$, $e = 2(f - e)$, $f = 3(f - e)$, $d = 3(e - d)$, $e = 4(e - d)$;
$3a = 2b$, $4b = 3c$; hence $c = 2a$, it follows that $(a, c)$ is an octave.
$4d = 3e$, $3e = 2f$; hence $f = 2d$, it follows that $(d, f)$ is an octave.
Thus, $(a, c)$, $(d, f)$ are equivalent dichords.
Similarly with *Case 2.*

*Note.* Nicomachus' account, attributing the discovery of the 4–chord to Pythagoras himself, describes lucidly the commutativity of the fifth with the fourth in 6.1, 57–62:

> ἑκατέρως τε
> ἡ διὰ πασῶν σύστημα ἠλέγχετο
> τῆς διὰ πέντε καὶ διὰ τεσσάρων *ἐν συναφῇ*,
> ὡς ὁ διπλάσιος λόγος ἤτοι ἡμιολίου τε καὶ ἐπιτρίτου,
> οἷον δώδεκα ὀκτὼ ἕξ,
> ἢ ἀναστρόφως τῆς διὰ τεσσάρων καὶ διὰ πέντε,
> ὡς τὸ διπλάσιον ἐπιτρίτου τε καὶ ἡμιολίου,
> οἷον δώδεκα ἐννέα ἕξ ἐν τάξει τοιαύτῃ.
> Nicomachus, *Harmonicum enchiridion* 6.1, 57–62
>
> It was also proved that the octave can be constructed in each of two ways,
> either as the composition (*en synaphei*) of the fifth and the fourth,
> since the duple ratio consists of the composition of hemiolic and epitritic – as in the numbers 12, 8, 6 – or the other way round, as the composition of the fourth and the fifth,
> since the duple ratio consists of the composition of epitritic and hemiolic – as in the numbers 12, 9, 6, which are ordered in that sort of way.
> [Nicomachus, *Harmonicum enchiridion*, trans. Barker, 1989, p. 257, with modifications by the authors]

The following table provides a comparison between the ancient empirical music and the arithmetized Pythagorean theory of music, following Hippasus' 4–chord.

| *Ancient empirical music*<br>(Section 1) | | *Pythagorean music arithmetized by Lasus–Hippasus' experiments and Hippasus' 4–chord* (Section 6) |
|---|---|---|
| Empirical knowledge<br>of the four musical intervals | | Mathematical definition<br>of the four musical intervals |
| *Empirical properties<br>without proof* | | *Propositions with proof<br>based on the definition* |
| Empirical Transitivity<br>Pythagoras' recognition<br>of basic intervals by hearing<br>Nicomachus, *Harmonicum enchiridion*<br>6,1,11–20 (Section 1.1) | | Transitivity<br>proved for the 4 musical intervals<br>in the 4–chord 6, 8, 9, 12 |
| Empirical Most equality<br>Pythagoras' care for constructing<br>exact, most equal copies<br>Nicomachus, *Harmonicum enchiridion*<br>6,1,26–35 (Section 1.2) | → | Most equality<br>proved for the 4 musical intervals<br>in the 4–chord 6, 8, 9, 12 |
| Empirical Commutativity<br>of musical composition | | Commutativity<br>of musical composition<br>proved for the two cases:<br>fourth & fifth<br>8/6 * 12/8 = 9/6 * 12/9,<br>fourth & tone 8/6* 9/8 = 9/8 * 8/6<br>in the 4–chord 6, 8, 9, 12.<br>Also, syllaba * di'oxeian =<br>di'oxeian * syllaba<br>(Philolaus, *Fragment* 6, 18–20). |
| Empirical equivalence of the dichord<br>$a > b$ with any equi–multiple $ka > kb$ | | Equivalence of the dichord<br>$a > b$ with any equi–multiple<br>$ka > kb$<br>12/6 = 2/1, 12/8 = 9/6 = 3/2,<br>12/9 = 8/6 = 4/3 |

**7.** *Summary of Book VII of Euclid's* Elements *and related Comments and Questions*

We will now have a careful look at Book VII with the ultimate purpose of correlating it with Pythagorean music (in Section 11). We will provide a summary of the fundamental definitions and propositions of Book VII of Euclid's *Elements*. The definitions, propositions, and their proofs contained in Book VII present some striking structural and conceptual peculiarities, mostly in the form of unnatural and seemingly unexplained omissions. Tannery recognized that some of these features, such as ratio multiplication, originated from music, but he did not further develop a rigorous analysis of them. Despite certain omissions, the book maintains a highly rigorous and deep level of mathematical content leading to the Fundamental Theorem of Arithmetic; however, its structure presents several peculiarities that are at odds with modern Number Theory. We point out these issues in the form of *Questions*, and provide explanations and answers in subsequent Sections.

**7.1.** *Basic Definitions of Book VII*

*Definition VII.1.* A unit (*monas*) is that by virtue of which each of the beings (*hekaston ton onton*) is called one (*hen*).

*Definition VII.2.* A number (*arithmos*) is a multitude (*plethos*) of units.
*Notes*. Euclid implicitly means that a number is a *finite* multitude of units.
The definition of a unit is philosophical/Platonic; its main mathematical significance is that the units in a number are equal.
This definition would suggest a discrete representation of numbers as a collection of units/dots (*psephides*).

*Question 1*. Why, in the *Elements,* the number, although defined, in Definition VII.2, as a [finite] multitude of units, is nevertheless *represented,* not as a collection of units/dots (*psephides*) but *as a linear segment*?

Before defining parts and multiples, it is necessary to consider the definition of multiplication.

*Definition VII.16*. A number $a$ is said to *multiply* another number $b$ when the multiplied number $b$ is added to itself as many times as there are units in the multiplying number $a$, thus defining the product $a \cdot b$.
*Note.* It does not follow from this definition that the commutative property $a \cdot b = b \cdot a$ holds for multiplication of numbers. In fact, this is proved in Proposition VII.16.

*Definition VII.3.* For numbers $a$, $b$, with $a < b$, the number $a$ is a *part* (*meros*) of the number $b$, when $a$ measures (*katametrei*) $b$.
*Note.* According to the definition of multiplication (VII.16), $a$ is a part of $b$ if $a < b$ and there exists a number $n$ such that $b = n \cdot a$.

*Definition VII.4.* For numbers $a$, $b$, with $a < b$, the number $a$ is *parts* (*mere*) of the number $b$, when $a$ does not measure (*katametrei*) $b$.

*Definition VII.5.* For numbers $a$, $b$, with $a > b$, the number $a$ is a *multiple* (*pollaplasios*) of the number $b$, when $a$ is measured by $b$.
*Note.* A number $a$ is a multiple of a number $b$ if and only if $b$ is a part of $a$, if and only if $a = nb$ for some number $n$.

*Definition VII.21.* Numbers $a$, $b$, $c$, and $d$ are *proportional* (*analogon*), (in symbols $a/b = c/d$), when the first $a$ is of the second $b$
either the equal (*isakis*) multiple,
or the same (*to auto*) part,
or the same (*ta auta*) parts that (as) the third $c$ is of the fourth $d$.
*Notes*. The definition of proportionality VII.21 is the central concept of Book VII. Nevertheless, Definition VII.21 as given is not immediately clear, since Euclid does not define anywhere what is meant by *equal multiple*, *same part*, and especially by *same parts*. While equal multiple and same part are easily guessed, the condition of the same *parts* is quite opaque and in need of clarification. This clarification is indeed supplied by Euclid in the first place, where these notions are applied, namely in Propositions VII.5 and VII.6. As a result of these clarifications, *Definition VII.21* assumes the following equivalent form:

$a/b = c/d$, if and only if there are numbers $n$, $m$, $k$, $l$, with
either $a = nb$ and $c = nd$ for some number $n > 1$,
or $b = na$ and $d = nc$ for some number $n > 1$,
or $b = nk$, $a = mk$, and $d = nl$, $c = ml$, for some numbers $n > m > 1$, $k$, $l$.

*Question 2.* Why does Definition VII.21 of proportion exclude the simple but important proportion $a/a = b/b$, and is given only for unequal ratios?

*Question 3.* Why is Definition VII.21 of proportion $a/b = c/d$ given fully only for the case $a < b$, $c < d$ in terms of "*a* part of *b*" and "*a* parts of *b*", but quite inadequately for the case $a > b$, $c > d$ in terms only of "*a* multiple of *b*"? For example, 2/4 = 3/6 (*a* part of *b*), 4/10 = 12/30 (*a* parts of *b*), and 4/2 = 6/3 (*a* multiple of *b*) are defined, while 10/4=30/12 is not defined according to Definition VII.21.
*Note*. In any case, Definition VII.21 is incomplete and must be completed in two ways; the complete mathematical definition in the spirit of Definition VII.21 assumes the following form: $a/b = c/d$ if and only if there are numbers $m$, $n$, $k$, $l$, such that $a = mk$, $b = nk$, $c = ml$, $d = nl$.
After the Definitions we would normally expect the appearance of Postulates. But there are no explicit Postulates for Numbers, unlike the situation of Geometry. But in fact, all of Book VII is based on the *Principle of the Least*, according to which

"every strictly decreasing sequence of natural numbers is finite",

a Principle that, as we know today, is logically equivalent to the Principle of Mathematical Induction, which serves, by means of the Peano Axioms, as the modern axiomatic foundation of natural numbers. The Principle of the Least is used twice in Book VII, for the proof of Propositions VII.2 and VII.31, in the latter formulated explicitly. It is not clear why this Principle is not elevated to the status of a Postulate, in analogy to what happens with the Postulates of Geometry. We may conjecture that the definition of number, employed in the *Elements* and in fact standard in Greek Mathematics, as a (finite) multitude of units, is so strong that it makes the Principle of the Least appear evident.

### **7.2.** *Propositions VII.1–4 and Definition VII.21**

The first two propositions contain a fundamental discovery of the Pythagoreans and form the base of Book VII.
Propositions VII.1 and VII.2 introduce the process of anthyphairesis (*the Euclidean Algorithm*) of two natural numbers $a$, $b$, with $a > b$, in order to prove that the last step of the anthyphairesis is the Greatest Common Measure/Divisor (*GCM*) of $a$ and $b$. (Proposition VII.1 for relatively prime numbers, Proposition VII.2 for not relatively prime numbers, employing the Principle of the Least to ensure that the process ends after a finite number of steps).

*Proposition VII.1.*
Two unequal numbers being set out, and the less being anthyphairated/mutually subtracted (*anthyphairoumenou*) always from the greater, if the number which is left never measures the one before it until a unit is left, the original numbers will be prime to one another.

*Proposition VII. 2.*
Given two numbers not prime to one another, to find their greatest common measure/divisor.
*Proof.* Let $a$, $b$ be two numbers, $a > b$. We construct the anthyphairesis of $a$ to $b$ exactly as in Proposition VII.1.

$a = k_0 b + c_1$, with $c_1 < b$,
$b = k_1 c_1 + c_2$, with $c_2 < c_1$,
$c_1 = k_2 c_2 + c_3$, with $c_3 < c_2$,
…
$c_{i-1} = k_i c_i + c_{i+1}$, with $c_{i+1} < c_i$,
…

Continuing, the anthyphairesis sequence of remainders is constructed, a sequence of numbers forming a strictly decreasing sequence:

$$a > b > c_1 > c_2 > \ldots > c_i > c_{i+1} > \ldots$$

At this point, the proof in the *Elements* states:

> the less of the numbers AB [*a*], CD [*b*] being repeatedly anthyphairated from the greater, some number will be left which will measure the one before it.

Euclid here makes an implicit appeal to the Principle of the Least:
*Every non–empty subset of the natural numbers has a least element.*
*Note:* Perhaps the careful reader will not be fully convinced that Euclid indeed has the Principle of the Least in mind. However, in the proof of

Proposition VII.31 below, the appeal to the Principle of the least is unmistakable:

> For, if it is not found, *an infinite sequence of numbers* will measure the number *a, each of which is less than the other:* which is impossible in numbers.

Thus, every anthyphairesis of a dyad of numbers necessarily ends after a finite number of steps, and the only way to end is when a remainder measures the immediately previous remainder. Thus, the two final steps will be as follows:

$$c_{n-2} = k_{n-1}c_{n-1} + c_n$$
$$c_{n-1} = k_n c_n$$

We finally claim that the last remainder $c_n$ is the greatest common measure of $a$ and $b$.

*Claim 1. $c_n$ is a common measure of a and b.*

For the proof, we proceed by finite recursion along the steps of the anthyphairesis, going from the end to the beginning.

From the very last relation, $c_n$ is a measure of $c_{n-1}$.

Hence, from the next–to–last relation, $c_n$ is a measure of $c_{n-2}$.

Hence, continuing in this way, from the second relation, $c_n$ is a measure of $b$, and hence, from the first relation, $c_n$ is a measure of $a$, thus $c_n$ is a common measure of $a$ and $b$.

*Claim 2. $c_n$ is the greatest common measure of a and b.*

For the proof, we proceed by finite recursion along the steps of the anthyphairesis, going from the beginning to the end, essentially imitating the recursion in the proof of Proposition VII.1.

Let $c$ be a measure/divisor of both $a$ and $b$.

The first anthyphairetic relation implies that $c$ is a measure of the first remainder $c_1$.

The second anthyphairetic relation implies that $c$ is a measure of the second remainder $c_2$. Proceeding in this way to the end, we conclude that $c$ must be a measure of the last anthyphairetic remainder $c_n$.

Thus, EVERY common measure of $a$ and $b$ is a measure of $c_n$. In particular, $c_n$ is the greatest common measure of $a$ and $b$.

*Definition of the process of anthyphairesis*

Let $a$, $b$ be two given unequal numbers, say $b < a$.

The *first step* of the anthyphairesis is to find a natural number, the quotient $k_0$, such that,

$$a = k_0 b + c_1, \text{ with } c_1 < b.$$

In the ancient process, this is achieved by repeated mutual subtraction:

Since $a > b$, we may subtract $b$ from $a$ and form $a - b$.

If $a - b < b$ then we set $c_1 = a - b$;

If $a - b \geq b$, then we subtract once more $b$ from $a$ and form $a - 2b$.

If $c_1 = 0$, then we have $a = k_0 b$, and the process comes to an end.

But if $c_1 > 0$, then we come to the decisive next *second step* of anthyphairesis, in which the role of the number $b$ is *reversed* from subtracting/dividing/active to being subtracted/divided/passive in relation to the active/dividing *first remainder* $c_1$ (hence the Greek term *ant–h*uphairesis, *reverse* subtraction). We find $k_1$, $c_2$, such that

$b = k_1c_1 + c_2$, with $c_2 < c_1$.

If $c_2 = 0$, then we have $b = k_1c_1$ and the process comes to the end.
The process of anthyphairesis can be continued with reversals as follows

$c_1 = k_2c_2 + c_3$, with $c_3 < c_2$,

…

$c_{i-1} = k_ic_i + c_{i+1}$, with $c_{i+1} < c_i$

…

The hypotheses of the proposition state that the last two steps in the anthyphairesis are the following:

$c_{n-2} = k_{n-1}c_{n-1} + 1$

$c_{n-1} = k_n.$

*Question 4*. Why, in Propositions VII.1 & VII.2 of the *Elements,* is anthyphairesis defined ONLY for a dyad of unequal numbers, and not for equals?

Propositions VII.1, VII.2 are the most fundamental in all three arithmetical Books VII, VIII, and IX of the *Elements*. The process of finding the Greatest Common Divisor of two numbers by anthyphairesis is usually called the *Euclidean algorithm,* but this is certainly not because the method is due to Euclid himself, but because it appears in Euclid's *Elements*. So, it is reasonable to ask:

*Question 5.* How did the Pythagoreans arrive at the discovery of the arithmetical anthyphairesis/Euclidean algorithm?

*Proposition VII.4.* Any number is either a part or parts of any number, the less of the greater.
First of all, Proposition VII.4 is stated in a way that makes it appear a tautology, since already by Definitions VII.3, 4, every number is either part or parts, the less of the greater. But the real meaning of Proposition VII.4 can be deduced from its proof. We understand that the meaning of the terms "part", "parts" is anthyphairetic, namely (a) if $a > b$, then $b$ is part of $a$ if and only if the anthyphairesis of $a$ to $b$ consists of one step, and $b$ is parts of $a$ if their anthyphairesis has length greater than one, and (b) in both cases $a$ and $b$ are measures of their Greatest Common Divisor. But at the same time, Proposition VII.4 shows that the distinction in part and parts is quite unnecessary, since both cases are replaced by the single description $a = mk$, $b = nk$ for some numbers $m$, $n$, where $k = GCD\ (a, b)$.

Thus, we are led to another question:

*Question 6.* Why did Euclid enunciate Definition VII.21 of proportion of numbers by introducing solely for this purpose the distinction between "same part" and "same parts", even though this distinction, as clarified by Proposition VII.4, is not fundamental for the definition of proportion, and it is not used after Proposition VII.13?

As suggested by Proposition VII.4, the ultimate importance of Propositions VII.1 and 2 is in strengthening Definition VII.21 so as to compare ratios of numbers by considering only the *GCM* of the numbers in each ratio:

*Definition VII.21**
$a/b = c/d$ if, setting $k = GCM\ (a, b)$, $l = GCM\ (c, d)$,
there are numbers $m$, $n$ such that
$a = m{\cdot}k$, $b = n{\cdot}k$, and $c = m{\cdot}l$, $d = n{\cdot}l$.
In fact, it may be proved that Definition VII.21, completed so as to take care both of ratios $a/a$ and of all cases of $a > b$, can be shown to be equivalent to Definition VII.21*.
The role of Propositions VII.1 and 2 is precisely in order to provide the crucial strengthening of Definition VII.21 to Definition VII.21* and they are never used directly after Proposition VII.4.

**7.3.** *Propositions VII.A, B, 11, 12, 13, 14, 16, 20*

The following proofs are based on the arithmetical definition of proportion (Definition VII.21*) derived from the Pythagorean musical–anthyphairetic tradition.

*Proposition VII.A.* (*Transitivity property of proportion*). If $a$, $b$, $a_1$, $b_1$ and $a_2$, $b_2$ are numbers such that $a/b = a_1/b_1$ and $a_1/b_1 = a_2/b_2$, then $a/b = a_2/b_2$.
*Proof.* Let $k = GCD\ (a, b)$, $k_1 = GCD\ (a_1, b_1)$, $k_2 = GCD\ (a_2, b_2)$.
By Definition VII.21*, since $a/b = a_1/b_1$, there are numbers $m$, $n$, such that
$a = mk$, $b = nk$, and $a_1 = mk_1$, $b_1 = nk_1$, and
since $a_1/b_1 = a_2/b_2$, there are numbers $m_1$, $n_1$, such that
$a_1 = m_1k_1$, $b_1 = n_1k_1$, and $a_2 = m_1k_2$, $b_2 = n_1k_2$.
Thus, $m = m_1$, $n = n_1$.
Hence $a = mk$, $b = nk$, and $a_2 = mk_2$, $b_2 = nk_2$,
namely $a/b = a_2/b_2$.
*Note.* Since reflexivity ($a/b = a/b$) and symmetry (if $a/b = c/d$, then $c/d = a/b$) are trivially true, it follows from Proposition VII.A that proportion is an equivalence relation on arithmetical dyads.

*Proposition VII.B.* (*"most equal", isaitata*)
Let $a$, $b$, $b_1$ be numbers. If $a/b = a/b_1$, then $b = b_1$.
*Proof.* By Definition VII.21, there are $m$, $n$, $k$, $k_1$ such that
$a = mk$, $b = nk$, and $a = mk_1$, $b_1 = nk_1$.
Then $mk = a = mk_1$, hence $k = k_1$;
hence $b_1 = nk_1 = nk = b$.
*Note.* Propositions VII.A and VII.B are not stated explicitly, but are necessary for the proofs in Book VII.

*Proposition VII.11.* If $a/b = c/d$ and $a < c$, $b < d$, then $a/b = (a - c) / (b - d)$.
*Proof.* By Definition VII.21, there are $k$, $l$, $m$, $n$, such that $a = mk$, $b = nk$, $c = ml$,
$d = nl$. Then $a - c = m(k - l)$, $b - d = n(k - l)$.
By Definition VII.21, $a/b = (a - c) / (b - d)$.

*Proposition VII.12.*
If $a_1/b_1 = a_2/b_2 = \ldots = a_n/b_n$, then $(a_1 + a_2 + \ldots + a_n)/(b_1 + b_2 + \ldots + b_n) = a_1/b_1$.
*Proof.* By Definition VII.21, there are numbers $m$, $n$, $k_1$, $k_2, \ldots, k_n$, $l_1$, $l_2, \ldots, l_n$, such that $a_1 = mk_1$, $a_2 = mk_2$, $\ldots$, $a_n = mk_n$, $b_1 = nl_1$, $b_2 = nl_2$, $\ldots$, $b_n = nl_n$.
We set $k = k_1 + k_2 + \ldots + k_n$ and $l = l_1 + l_2 + \ldots + l_n$.
Then $a_1 + a_2 + \ldots + a_n = mk$, $b_1 + b_2 + \ldots + b_n = nl$.
By Definition VII.21, the conclusion follows.

An immediate *Corollary* to Proposition VII.12 is:
If $a$, $b$, and $n$ are numbers, then $na/nb = a/b$.

*Proposition VII.13.* (*Alternando/Enallax*)
If $a/b = c/d$, then $a/c = b/d$.
*Proof.* By Definition VII.21, there are numbers $m$, $n$, $k$, $l$
such that $a = mk$, $b = nk$, $c = ml$, $d = nl$.
By the Corollary to Proposition VII.12,
$a/c = mk/ml = k/l$, and $b/d = nk/nl = k/l$.
By Proposition VII.A, $a/c = b/d$.

*Proposition VII.14.* (*di' isou, ex aequali*)
If $a/b = d/e$, $b/c = e/f$, then $a/c = d/f$
*Proof.* By VII.13, $a/d = b/e$, $b/e = c/f$
By VII.A $a/d = c/f$
By VII.13 $a/c = d/f$
Proposition VII.14 shows that the operation of multiplication of ratios (positive fractions, in modern terminology), defined implicitly therein, is well defined for the equivalence classes of the ratios (positive rationals, in modern terminology), as well.

*Question 7.* Why is the operation * that corresponds to the multiplication of ratios restricted to *"homonymous"* ratios, namely only to multiplication of ratios of the form $(a_1/a_2) * (a_2/a_3) = a_1/a_3$, rather than being defined in general?

In Book VII of the *Elements,* there is no operation of addition for ratios of numbers.

*Question 8.* Why is the *addition of ratios* not defined at all in the *Elements* although it would be quite easy, especially so since the Least Common Multiple of two numbers is constructed in Proposition VII.34?

*Proposition VII.16. (Commutativity of multiplication of numbers).* For any two numbers $a$ and $b$, $ab = ba$.
*Note.* A straightforward proof of *Proposition VII.16* is easily available; there is no need to do anything more than observe the following:

```
                      •  •  •
•  •  •  •  •         •  •  •
•  •  •  •  •    =    •  •  •
•  •  •  •  •         •  •  •
                      •  •  •
```

Thus, the proof given by Euclid is most unexpected; it is based essentially on the Alternando/Enallax property (Proposition VII.13). It goes as follows:
$1/a = b/ba$, $1/b = a/ab$ by VII.12,
$1/a = b/ab$ by Alternando/Enallax VII.13,
$ab = ba$ by VII.B.
Thus, the Pythagoreans must have considered Proposition VII.13 as a commutativity property.

*Question 9*. Why did Euclid prove *the commutativity of multiplication* of numbers with the quite unexpected proof by using Proposition VII.13 *Alternando/Enallax*? In what sense would the Alternando property be considered *a commutativity property?*

*Proposition VII.20. (*Generation of each equivalence class by a unique ratio)
If $a$, $b$ are relatively prime numbers and $c$ and $d$ are numbers such that $a/b = c/d$, then there is a number $k$ such that $c = k \cdot a$, and $d = k \cdot b$.
*Proof.* Since $a$, $b$ are relatively prime, then $GCD\ (a, b) = 1$.
We set $k = GCD\ (c, d)$.
By Definition VII.21*, there are numbers $m$, $n$, such that $a = m \cdot 1$, $b = n \cdot 1$,
$c = m \cdot k$, $d = n \cdot k$. Thus, $m = a$ and $b = n$, and $c = k \cdot a$, $d = k \cdot b$.

*Note.* Proposition VII.20 states that every *equivalence class* under proportion is *generated* by the unique element of the class, the ratio consisting of two relatively prime numbers.
The propositions of Book VII after Proposition VII.20 are mostly concerned not with ratios of numbers but with numbers themselves (such as VII.27, VII.30 and VII.31). It would seem strange for a book on the theory of numbers to start with Propositions on ratios and rational numbers and continue with Propositions about natural numbers, and in fact using Propositions about proportion to prove basic number–theoretic theorems, such as VII.27 and VII.30. This is the content of our next question.

*Question 10*. Why does the arithmetical Book VII of the *Elements* begin with *Propositions about ratios and proportions of numbers,* not with numbers themselves, and continues with propositions about numbers themselves? It is like today one would first deal with rational numbers and then with natural numbers themselves.

**8.** *Aristotle's* Topics *158b24–29 Principle of a dynamic relation between Hypotheses/Postulates and Definitions*

At this point, we find it useful to introduce a Principle, enunciated by Aristotle, on the dynamic interaction between Definitions and Hypotheses/Postulates. We first describe this Principle.
According to Aristotle's Topics *158b24–29* Principle, a system of empirical properties cannot be proven before the discovery of a good definition, but these properties become propositions with rigorous proofs once a suitable definition is found.

> Πολλαῖς τε τῶν θέσεων
> μὴ καλῶς ἀποδιδομένου τοῦ ὁρισμοῦ
> οὐ ῥᾴδιον διαλέγεσθαι καὶ ἐπιχειρεῖν,
> οἷον πότερον ἓν ἑνὶ ἐναντίον ἢ πλείω·
> ὁρισθέντων δὲ τῶν ἐναντίων κατὰ τρόπον ῥᾴδιον συμβιβάσαι
> πότερον ἐνδέχεται πλείω τῷ αὐτῷ εἶναι ἐναντία ἢ οὔ.
> τὸν αὐτὸν δὲ τρόπον καὶ ἐπὶ τῶν ἄλλων τῶν ὁρισμοῦ δεομένων.
> Aristotle, *Topics Books I and VIII,* 158b24–29
>
> And with many theses (*theseon*),
> if the definition (*horismou*) is not given well,
> it is not easy (*radion*) to argue (*dialegesthai*) and deal (*epicheirein*) with them,
> e.g., whether one thing is contrary (*enantion*) to one, or several.
> But once "contraries" have been properly defined (*horisthenton*),
> it is easy to infer whether it is possible for there to be several contraries for the same thing or not.

> It is the same way also in other cases requiring definitions (*horismou*).
> [Aristotle, *Topics*, trans. Smith, 1997, p. 27]

Some of the "theses"/hypotheses are not easy, or are in fact impossible to establish or prove, if a suitable definition is lacking. But if a thesis that we believe to be true cannot be proved, because of a lack of a suitable definition, then we are obliged to take this thesis as an Axiom/Postulate without proof.

Thus, Aristotle describes a dynamic relation between Definitions and Postulates. If we have no suitable definition of a notion, then we will necessarily have as axioms all the Statements that employ this notion. But if we succeed, at some later time, in finding a suitable and good definition of the notion, then the Axioms/Postulates without proof will end up becoming Propositions with proof.

Aristotle describes an interesting occurrence of this Principle, as follows: At some early stage of Geometry, when there was no good definition of proportion of magnitudes, it was impossible, exactly because of the lack of the definition, to prove what is now Proposition VI.1 of the *Elements*, namely that the ratio of two lines $a/b$ is equal to the ratio of the rectangles $ac/bc$, and thus, this being a statement intuitively deemed as true and also useful for the theory, we would be obliged to take this proposition as a Postulate. But when Theaetetus provided a "good definition" of proportion of magnitudes, namely for magnitudes $a$, $b$, $c$, $d$, the proportion $a/b = c/d$ means that the sequence of the quotients of the anthyphairesis of $a$ to $b$ is equal to the sequence of the quotients of the anthyphairesis of $c$ to $d$, then Proposition VI.1, up to that point a Postulate without proof, now becomes a Proposition, with a (simple) proof. Schematically, we have the following dynamic interaction between the Definition of proportion and the Postulate/Proposition VI.1.

Aristotle's example of his *Topics* Principle

| *Empirical Proportion of magnitudes* | | *Theaetetus' theory of Proportion of magnitudes with a mathematical definition* |
|---|---|---|
| No mathematical definition of proportion | → | *Definition.* For lines $a/b = c/d$ if the sequences of the quotients of the anthyphairesis of $a$ to $b$ and of $c$ to $d$ are equal |
| The Basic property of proportion stating that "if $a$, $b$, $c$ are lines, then $a/b = ac/bc$" is a Postulate (without proof) | | The previous Postulate becomes Proposition VI.1 with proof |

We wish to comment further on the implications of Aristotle's Principle. Usually, in an axiomatic system, we distinguish among Definitions, Postulates and Propositions with proof. But Aristotle's Principle suggests a finer distinction: the propositions must be distinguished between those which are *one step removed* from the Postulates and Definitions, namely those that are proved by a direct appeal to the Definitions and Postulates (call them Class I), and those which are removed more than one step, namely that are proved solely from the Propositions of Class I, without any other appeal to Definitions and Postulates (call them Class II).
*Note:* We can regard the Pythagorean arithmetized musical theory based on the definition of the four musical intervals in relation to the previous empirical music as an instance of Aristotle's *Topics* Principle.

**9.** *Our idea is to apply Aristotle's* Topics *Principle to Book VII of the* Elements*, in order to obtain the state of the Pythagorean arithmetic, at the time after the development of the Pythagorean arithmetized theory of music, but before the discovery of Propositions VII.1 &2 and Definition VII.21**

This example mentioned by Aristotle on proportion in Geometry (in Section 8) is close to the situation with proportion in numbers. Surely there was an earlier time in the development of Pythagorean Number Theory, when, on the one hand it would be thought desirable to include the propositions in Class I of Book VII, but, on the other hand, Definition VII.21 was not yet available, so that, Aristotle's *Topics* Principle would be applicable, and all the propositions of Class I would have to function as Postulates (without proof) in order to logically sustain the proofs of Class II propositions.
Going now to Book VII, it is not difficult to realize, first, that these five or six properties correspond to Propositions VII.A (transitivity), VII.B (most equality), VII.11 (diairesis), VII.12 (synthesis), VII.13 (commutativity), and VII. 20 (generation of each equivalence class by a unique ratio), and secondly, these are exactly the propositions of Class I of Book VII, in the sense described in Section 8, namely, those propositions whose proofs require a direct appeal to the Definition VII.21*.
Schematically we have the following dynamic interaction between Definition VII.21* and the Postulates/Propositions VII.A, VII.B, VII.11, VII.12, VII.13, VII.20.

Backward application of Aristotle's *Topics* Principle

| _Pythagorean Arithmetic_<br>**_before_** _the discovery of_<br>_the "Euclidean algorithm"_<br>_and Definition VII.21*_<br>(Section 10) | | _Pythagorean Arithmetic (Book VII)_<br>**_based on_**<br>_the "Euclidean algorithm"_<br>_and Definition VII.21*_<br>(Section 8) |
|---|---|---|
| **No** mathematical definition<br>of proportion for arithmetical ratios | → | Mathematical definition<br>of proportion<br>Definition VII.21* |
| All immediate consequences<br>of the definition of proportion<br>VII.21*<br>**must be taken as Postulates:** | | All immediate consequences<br>of the definition of proportion<br>VII.21*<br>**now have a proof:** |
| VII.A (Transitivity) | | Proposition VII.A (Transitivity) |
| VII.B (Most equality) | | Proposition VII.B (Most equality) |
| VII.11, VII.12. | | Propositions VII.11, VII.12. |
| VII.13 (Alternando/Enallax) | | Proposition VII.13<br>(Alternando/Enallax) |
| VII.20<br>(generation of each equivalence<br>class by a unique ratio) | | Proposition VII.20<br>(generation of each equivalence<br>class by a unique ratio) |

What separates the two states described in the two columns is the discovery of arithmetical anthyphairesis (*the Euclidean algorithm*), leading to the method of finding the Greatest Common Measure and the subsequent "good definition" of proportion based on this method.
Our suggestion for how this discovery came about will be described in Sections 11 and 12. In the meantime, we will explain, in the next Section, how the Pythagoreans were motivated to consider plausible and to conjecture the first column.

**10.** *The Sweeping Pythagorean conjecture*

We cannot fail to observe that the list of the properties that have been proved to hold for the four basic musical intervals (in earlier Sections) is practically identical with the list of the propositions VII.A, VII.B, VII.11, VII.12, VII.13, VII.20 in Book VII that hold for all the arithmetical ratios. These specific Propositions are the immediate consequences of the Definition VII.21* and are thus turned into Postulates at an earlier stage, following the application of Aristotle's *Topics* Principle. (The only non–obvious correspondence between the commutativity of the composition of musical intervals and the Alternando Proposition VII.13 for arithmetical ratios, will be made clear in Section 12.8, below).
We are then led to conclude that the Pythagoreans,
**after** the development of the arithmetized music theory of the 4–chord, described in earlier Sections, achieved by the forward progress from empirical music to Hippasus' 4–chord to the Pythagorean mathematical theory of the four musical intervals of the 4–chord,
**but before** the discovery of the Euclidean algorithm (VII.1 & 2) and Definition VII.21* of proportion, conceived by applying the backward analytical step suggested by Aristotle's *Topics* Principle – moving from the structured Euclidean Book VII back to a state of primitive arithmetic – must have proceeded to *a sweeping conjecture* that what has been proved to hold for a minute part of arithmetical ratios, namely the multiple (Section 11.4.1) and the epimoric ratios (Section 11.4.2), must in fact hold for all arithmetical ratios without exception.
Schematically, we have:

The sweeping Pythagorean conjecture

| *Pythagorean arithmetized music theory* ***after*** *the construction of Hippasus' 4–chord* | | *Pythagorean Arithmetic* ***before*** *the discovery of the "Euclidean algorithm" and the Definition VII.21** |
|---|---|---|
| Definition of musical intervals $a$, $b$ octave if $a = 2b$, $a$, $b$ fifth if $a = 3(a–b)$, $b = 2(a–b)$, analogously for fourth, tone | → | **No** mathematical definition of proportion for arithmetical ratios |
| Propositions with proofs based on the definition for the four basic musical intervals: | | All immediate consequences of the definition of proportion VII.21* **must be taken as Postulates:** |
| Transitivity for dichords within the same musical interval | | VII.A (Transitivity) |
| Most equality proved for the 4 musical intervals in the 4–chord 6, 8, 9, 12 | | VII.B (Most equality) |
| Equivalence of the dichord $a > b$ with any equi–multiple $ka > kb$ 12/6 = 2/1, 12/8 = 9/6 = 3/2, 12/9 = 8/6 = 4/3 in the 4–chord 6, 8, 9, 12 | | VII.11, VII.12 |
| Commutativity of *fourth* with *fifth* 8/6 * 12/8 = 9/6 * 12/9, of *fourth* with *tone* 8/6 * 9/8 = 9/8 * 12/9 in the 4–chord 6, 8, 9, 12. (Nicomachus, *Harmonicum enchiridion* 6.1, 57–62, Section 1.3) | | VII.13 (Alternando/Enallax) |
| | | VII.20 (generation of each equivalence class by a unique ratio) |

**11.** *How did the Pythagoreans discover anthyphairesis/Euclidean algorithm and the canonical definition of proportion of numbers VII.21*?*

**11.1.** *The discovery of the arithmetical anthyphairesis/Euclidean algorithm*

*Question 5* (posed in Section 7.2). How did the Pythagoreans arrive at the discovery of the arithmetical anthyphairesis/Euclidean algorithm?

We will provide our answer to Question 5 in the present Section. We will argue now that the Pythagoreans' preoccupation with music was crucial in several essential ways towards this discovery. Szabó was the first to have the key insight that Pythagorean music was the origin of this Pythagorean discovery.

**11.2.** *Because the Pythagoreans wished to develop a general theory of ratios of natural numbers generalizing their arithmetized theory of music, they had a strong motive to prove that every arithmetical ratio is proportional to a ratio in least terms, equivalently they had a strong motive to find a method of constructing the greatest common divisor of any two numbers*

As mentioned in Section 10, the Pythagoreans had a *strong motive* to prove that every ratio of numbers is equivalent to one in least terms; equivalently they had a *strong motive* to find the Greatest Common Divisor of any two numbers, exactly because they wished to imitate and generalize their arithmetized musical theory, in which the fundamental musical intervals were described by arithmetical ratios in least relatively prime terms, to a theory of all arithmetical ratios.

**11.3.** *The representation of numbers as continuous chords/straight lines facilitates the discovery of the Euclidean algorithm*

The fact that the Pythagoreans did not represent natural numbers discretely by dots/pebbles but continuously by strings/line segments facilitated the discovery by means of the Euclidean algorithm. The question is how the Pythagoreans were helped in the discovery of the Euclidean algorithm by going from numbers to straight lines.
If the numbers had been given not as continuous quantities (*chords*), but in a discrete way (*pebbles*), then we would have had nothing but to count the pebbles, with each of these pebbles being an obvious unit. The existence of these discrete (*minimal*) units would obscure and impede the search for, and identification of, the non–obvious, but more importantly, the greatest common measure, and would not push them to discover an internal method for finding the greatest common measure.

However, the line segment is a continuous, and not discrete, quantity, and if there is some unit that measures these two straight segments, then this unit is hidden, to be found and not apparent. If, therefore, we seek an arithmetical relation between two quantities that are continuous (such as chords, straight lines, and do not consist of discrete pieces), measured internally, without external aids, then essentially only one method is available, the interaction between them. The continuous representation of two unequal numbers by two unequal chords/straight line segments hides the unit and, in fact, facilitates, even makes necessary, the comparison of the two strings directly between them, internally, bringing them into some interaction. The only conceivable way of comparing two line segments is to measure the greater by the smaller.
We cannot do anything other than divide the large interval by the small one, and we are necessarily led to anthyphairesis, precisely because there is no pre–given arithmetical unit to offer a trivial reliance on practical rules of divisibility.

**11.4.** *The generalization of the four irreducible ratios of the early Pythagorean theory of music leads to the anthyphairetic definitions of the multiple, the epimoric, and multi–epimoric ratios, namely to the ratios with anthyphairesis having length 1 and 2, respectively*

What is needed, according to the application of Aristotle's *Topics* Principle, is to discover Definition VII.21*. For the discovery of Definition VII.21*, the discovery of Propositions VII.1 and 2 on general arithmetic anthyphairesis are needed.
The steps of the discovery of the general arithmetic anthyphairesis will be described as an upward ladder starting from the multiple ratios (Section 11.4.1), proceeding to the epimoric and poly–epimoric ratios (Sections 11.4.2, 11.4.3), next, proceeding to the epimeric and poly–epimeric ratios (Section 11.4.4), next, formulating the general upward step, and finally, reversing the process and proceeding downwards, introducing the Principle of the Least, for the final proof of Propositions VII.1 and 2. We rely heavily on Theon's account *Expositio rerum mathematicarum ad legendum Platonem utilium, [De utilitate mathematicae]*.

**11.4.1.** *The musical interval of the octave, generated by the irreducible double ratio 2/1, was generalized to the multiple ratio m/1*

The octave is generated by a dichord with the greater chord being equal to double the smaller chord. This is generalized to a ratio with the greater number being a multiple of the smaller, or the smaller being a divisor of the greater one.

πολλαπλάσιος μὲν οὖν ἐστι λόγος,
ὅταν ὁ μείζων ὅρος πλεονάκις ἔχῃ τὸν ἐλάττονα
Theon 76, 8–9

XXIII. The relationship is *multiple*
when the larger term contains the smaller several times
[Theon of Smyrna, trans. Lawlor & Lawlor, 1979, with modifications by the authors]
*Note.* This corresponds to Definitions VII.3, 5 in the *Elements*.
ἔστι δὲ τῶν πολλαπλασίων λόγων
πρῶτος καὶ ἐλάχιστος ὁ διπλάσιος,
μετὰ δὲ τοῦτον ὁ τριπλάσιος,
εἶτα ὁ τετραπλάσιος,
καὶ οὕτως οἱ ἑξῆς ἐπ' ἄπειρον ἀεὶ οἱ μείζονες.
Theon 77,23–78,1

Among the *multiple* ratios,
the first and smallest is the *double*,
next comes the *triple*,
then the *quadruple*,
and so forth increasing indefinitely.
[Theon of Smyrna, trans. Lawlor & Lawlor, 1979, with modifications by the authors]

**11.4.2.** *The generalization of the ratios 3/2, 4/3, and 9/8, corresponding to the musical intervals fifth, fourth and tone, are generalized to the epimoric ratios of the form (n+1)/n*

The most crucial step was taken with the generalization of the three remaining ratios 3/2, 4/3, and 9/8 corresponding, according to the Pythagorean experiments, to the musical intervals of the fifth, fourth, and tone.

*The general epimoric ratio (n+1)/n*

τῶν δ' ἐπιμορίων λόγων πρῶτος καὶ μέγιστος ὁ ἡμιόλιος,
ὅτι δὴ καὶ τὸ ἥμισυ μέρος πρῶτον καὶ μέγιστον καὶ ἐγγυτάτω τῷ ὅλῳ,
μετὰ δὲ τοῦτον ὁ ἐπίτριτος, καὶ ὁ ἐπιτέταρτος,
καὶ οὕτω πάλιν ἐπ' ἄπειρον ἡ πρόοδος ἀεὶ ἐπ' ἐλάττονος.
Theon 78,1–5

Among the epimoric ratios, the first and largest is the *sesquialter* ratio (3/2), because the fraction 1/2 is the first, the largest and the one which most closely approaches the whole; then comes the *sesquitertian* ratio (4/3), then the *sesquiquartan* ratio (5/4) and so forth indefinitely, always proceeding by diminishing.
[Theon of Smyrna, trans. Lawlor & Lawlor, 1979, with modifications by the authors]

**11.4.3.** *The crucial step is the anthyphairetic definition of the multiple, epimoric, and poly–epimoric ratios given by Theon*

(i) Theon's definition of a *multiple* ratio as the ratio of $a$ to $b$ having an *anthyphairesis Anth* $(a, b) = [n]$

τουτέστιν ὅταν ὁ μείζων ὅρος καταμετρῆται ὑπὸ τοῦ ἐλάττονος
ἀπαρτιζόντως, ὡς μηδὲν ἔτι λείπεσθαι ἀπ' αὐτοῦ
Theon 76,9–11

that is, when the smaller term exactly measures the larger
without there remaining any part left over
[Theon of Smyrna, trans. Lawlor & Lawlor, 1979, with modifications by the authors]

Thus, $a$ is a multiple of $b$ if $a$ consists of $m$ copies of $b$, $a = mb$, with no remainder left.
Multiple (*pollaplasios*) is the ratio $a/b$ such that *Anth* $(a, b) = [k]$, with $k > 1$.

(ii) Theon's definition of an *epimoric* ratio as the ratio of $a$ to $b$ having an *anthyphairesis Anth* $(a, b) = [1, n]$

ἐπιμόριος δέ ἐστι λόγος,
ὅταν ὁ μείζων ὅρος ἅπαξ ἔχῃ τὸν ἐλάττονα καὶ μόριον ἕν τι τοῦ ἐλάττονος,
τουτέστιν ὅταν ὁ μείζων τοῦ ἐλάττονος ταύτην ἔχῃ τὴν ὑπεροχήν,
ἥτις τοῦ ἐλάττονος ἀριθμοῦ μέρος ἐστίν.
Theon 76,21–77,2

the greater term contains the smaller term plus a part of the smaller term one time, that is to say, when the larger term is greater than the smaller by a certain quantity which is a part of it.
[Theon of Smyrna, trans. Lawlor & Lawlor, 1979, with modifications by the authors]

Thus, if the ratio of two numbers $a$, $b$, $a > b$, is *epimoric*,
then the anthyphairesis of $a$ to $b$ has exactly two steps,

$a = b + c$, $c < b$,
$b = nc$

with sequence of quotients 1, n, namely *Anth* $(a, b) = [1, n]$,
hence $a = (n + 1)c$, hence $a/b = (n + 1)/n$

(iii) Theon's definition of a poly–epimoric ratio of $a$ to $b$ as the ratio having an *anthyphairesis Anth* $(a, b) = [k, n]$, with $k > 1$.

πολλαπλασιεπιμόριος δέ ἐστι λόγος,
ὅταν ὁ μείζων ὅρος δὶς ἢ πλεονάκις ἔχῃ τὸν ἐλάττονα καὶ ἔτι μέρος αὐτοῦ,
ὡς ὁ μὲν τῶν ζʹ δὶς ἔχει τὸν γʹ καὶ ἔτι τρίτον αὐτοῦ,
καὶ λέγεται αὐτοῦ διπλασιεπίτριτος,
ὁ δὲ τῶν θʹ δὶς ἔχει <τὸν> τῶν δʹ καὶ ἔτι τὸ τέταρτον αὐτοῦ,
λέγεται δὲ διπλασιεπιτέταρτος,
ὁ δὲ τῶν ιʹ τρὶς ἔχει τὸν τῶν γʹ καὶ τὸ τρίτον αὐτοῦ,
καὶ λέγεται τριπλασιεπίτριτος.
παραπλησίως δὲ θεωρείσθωσαν καὶ οἱ λοιποὶ πολλαπλασιεπιμόριοι.
τοῦτο δὲ συμβαίνει, ὅταν δυεῖν προτεθέντων ἀριθμῶν
ὁ ἐλάττων καταμετρῶν τὸν μείζονα
μὴ ἰσχύσῃ ὅλον καταμετρῆσαι,
ἀλλ' ἀπολείπῃ μέρος τοῦ μείζονος, ὅ ἐστιν αὐτοῦ τοῦ ἐλάσσονος μέρος·
οἷον ὁ τῶν κϛʹ τοῦ τῶν ηʹ πολλαπλασιεπιμόριος λέγεται,
ἐπειδή περ <ὁ> ηʹ τρὶς καταμετρήσας τὸν κϛʹ
οὐχ ὅλον ἀπήρτισεν,
ἀλλὰ μέχρι τῶν κδʹ ἐλθὼν δύο ἐκ τῶν κϛʹ ἀπέλιπεν, ὅ ἐστι τῶν ηʹ τέταρτον.
Theon 78,23–79,14

XXVI. The ratio is called *poly–epimoric* when the larger term contains the smaller two or more times plus a part of this smaller term.
7 contains in this way, 2 times 3 and in addition, a third of 3. Also, it is said that the relationship of 7 to 3 is *bi–epimoric*.
Likewise, 9 contains 2 times 4 and the fourth of 4 in addition; the ratio of 9 to 4 is *bi–epimoric.*
Again likewise, 10 contains 3 times 3, along with the third of 3, and the ratio is called *tri–epimoric*.
Other *poly–epimoric* ratios are recognized in the same manner. They occur in every case where, of the two proposed numbers, the smaller does not measure the larger exactly, but when the larger gives a remainder which is at the same time a remainder of the smaller.

Thus, the ratio 26 to 8 is *poly–epimoric* because 3 times 8 does not give 26 completely, but comes to 24 rather than 26, and there is a remainder of 2, which is a quarter.
[Theon of Smyrna, trans. Lawlor & Lawlor, 1979, with modifications by the authors]

Thus, a ratio of $a$ to $b$ has *anthyphairesis of length 2* if and only if the ratio of $a$ to $b$ is either *epimoric or poly–epimoric* if and only if $a = kb + c$, with $c < b$, $k$ any natural number, and $b$ is a multiple of $c$.

**11.4.4.** *Theon continues with the definition of the epimeric and poly–epimeric ratios a to b, which are precisely the ratios with length 3 of their anthyphairesis, while mentioning examples of ratios with length of anthyphairesis 4 and 5*

(i) epimeric ratios

ἐπιμερὴς δέ ἐστι λόγος,
ὅταν ὁ μείζων ὅρος ἅπαξ ἔχῃ τὸν ἐλάττονα
καὶ ἔτι πλείω μέρη αὐτοῦ [τοῦ ἐλάττονος],
εἴτε ταὐτὰ καὶ ὅμοια εἴτε ἕτερα καὶ διάφορα·
ταὐτὰ μὲν οἷον δύο τρίτα ἢ δύο πέμπτα καὶ εἴ τινα ἄλλα οὕτως·
ὁ μὲν γὰρ τῶν ε′ ἀριθμὸς τοῦ τῶν γ′ δὶς ἐπίτριτος,
ὁ δὲ τῶν ζ′ τοῦ τῶν ε′ δὶς ἐπίπεμπτος,
ὁ δὲ τῶν η′ τοῦ τῶν ε′ τρὶς ἐπίπεμπτος,
καὶ οἱ ἑξῆς ὁμοίως·
ἕτερα δὲ καὶ διάφορα
οἷον ὅταν ὁ μείζων αὐτόν τε ἔχῃ τὸν ἐλάττονα καὶ ἔτι ἥμισυ αὐτοῦ
καὶ τρίτον, οἷον ἔχει λόγον ὁ τῶν ια′ πρὸς τὸν τῶν ϛ′,
ἢ πάλιν ἥμισυ καὶ τέταρτον, ὅς ἐστι λόγος τῶν ζ′ πρὸς δ′,
ἢν ἡ Δία τρίτον καὶ τέταρτον, ὃν ἔχει λόγον τὰ ιθ′ πρὸς τὰ ιβ′.
Theon 78,6–17

XXV. A ratio is called *epimeric* when the larger term contains
the smaller one time plus several other parts of it, either similar parts or different parts; similar as two–thirds, two–fifths, etc.
Thus, the number 5 contains 3 plus two–thirds of 3;
[the ratio of 5 to 3 is epimeric,
since $5 = 1 \cdot 3 + 2$, and the ratio of 3 to 2 is epimoric]
the number 7 contains 5 plus two–fifths of 5
[the ratio of 7 to 5 is epimeric,
since $7 = 1 \cdot 5 + 2$, and the ratio of 5 to 2 is poly–epimoric];
the number 8 contains 5 and three–fifths of 5
[the ratio 8 to 5 has an anthyphairesis of length 4,
since $8 = 1 \cdot 5 + 3$, and 5 to 3 is epimeric];

and so on.
The parts are different when the largest term contains the smallest and in addition, its half and its third, as in the ratio of 11 to 6
[the ratio 11 to 6 is epimeric,
since 11 = 1·6 + 5, and the ratio 6 to 5 is epimoric],
or its half and its quarter, as in the ratio of 7 to 4
[the ratio 7 to 4 is epimeric,
since 7 = 1·4 + 3, and the ratio 4 to 3 is epimoric],
or again, by Zeus, the third and the quarter, as in the ratio of 19 to 12
[the ratio 19 to 12 has an anthyphairesis of length 5,
Since 19 = 1·12 + 7, 12 = 1·7 + 5, and 7 to 5 is epimeric].
[Theon of Smyrna, trans. Lawlor & Lawlor, 1979, with modifications by the authors]

*Note.* A ratio *a* to *b* is epimeric if the length of the anthyphairesis of *a* to *b* is 3, and the first quotient is 1.

(ii) poly–epimeric ratios

πολλαπλασιεπιμερὴς <δέ> ἐστι λόγος,
ὅταν ὁ μείζων ὅρος δὶς ἢ πλεονάκις ἔχῃ τὸν ἐλάττονα
καὶ δύο ἢ πλείω τινὰ μέρη αὐτοῦ εἴτε ὅμοια εἴτε διάφορα·
οἷον ὁ μὲν τῶν η′ δὶς ἔχει τὸν τῶν γ′ καὶ δύο τρίτα αὐτοῦ,
λέγεται δὲ διπλάσιος καὶ δὶς ἐπίτριτος,
ὁ δὲ τῶν ια′ τοῦ τῶν γ′ τριπλάσιος καὶ δὶς ἐπίτριτος,
ὁ δὲ τῶν ια′ τοῦ τῶν δ′ διπλάσιός τε καὶ ἡμιόλιος καὶ ἐπιτέταρτος
ἢ διπλάσιός τε καὶ τρὶς ἐπιτέταρτος.
καὶ τοὺς ἄλλους δὲ πολλαπλασιεπιμερεῖς
πολλοὺς καὶ ποικίλους ὄντας προχειρίζεσθαι ῥᾴδιον.
τοῦτο δὲ γίνεται, ὅταν ὁ ἐλάττων ἀριθμὸς καταμετρήσας τὸν μείζονα
μὴ ἰσχύσῃ ἀπαρτίσαι,
ἀλλ᾽ ἀπολείπῃ ἀριθμόν τινα, ἅ ἐστι μέρη αὐτοῦ,
ὡς ὁ τῶν ιδ′ τοῦ τῶν γ′·
ἡ γὰρ τριὰς καταμετρήσασα τὸν τῶν ιδ′ οὐκ ἴσχυσεν ἀπαρτίσαι,
ἀλλὰ προκόψασα τετράκις μέχρι τῶν ιβ′
τὴν λοιπὴν ἀπὸ τῶν ιδ′ ἀπέλιπε δυάδα,
ἥ τις ἐστὶ τῶν γ′ δίμοιρον, ἃ δὴ λέγεται δύο τρίτα.
ἀντίκειται δὲ καὶ τῷ πολλαπλασιεπιμερεῖ ὁ
ὑποπολλαπλασιεπιμερής.
Theon 79,15–80,6

XXVII. A ratio is called *poly–epimeric*
when the larger term contains the smaller two times or more,
along with two or several parts of the latter,

whether they be similar or different.
Thus, 8 contains 2 times 3 and in excess, two–thirds of 3,
and the ratio is called double with two–thirds in excess;
[The ratio 8 to 3 is poly–epimeric, since $8 = 2 \cdot 3 + 2$, and the ratio 3 to 2 is epimoric] [Theon of Smyrna, trans. Lawlor & Lawlor, 1979, with modifications by the authors]

Likewise, the ratio of 11 to 3 is triple with two–thirds in excess;
[The ratio 11 to 3 is poly–epimeric,
since $11 = 3 \cdot 3 + 2$, and the ratio 3 to 2 is epimoric]
The ratio of 11 to 4 is double with three–quarters in excess
[The ratio 11 to 4 is poly–epimeric,
since $11 = 2 \cdot 4 + 3$, and the ratio 4 to 3 is epimoric]
It is easy to find many other poly–epimeric ratios,
and this takes place each time that the smaller number does not exactly measure the larger, but there is a remainder formed of several parts of the smaller number, as in the ratio of 14 to 3, since 3 does not exactly measure 14, but 4 times 3 are 12, and of 14 there remains 2, which is two parts of three, and which is called two–thirds.
[The ratio 14 to 3 is poly–epimeric,
since $14 = 4 \cdot 3 + 2$, and the ratio 3 to 2 is epimoric]

*Note.* A ratio *a* to *b* is poly–epimeric if the length of the anthyphairesis of *a* to *b* is 3, and the first quotient is greater than 1.

*Note.* Iamblichus on whether the interval 8/3 [= 2/1 * 4/3] is consonant:
τὸ δὲ διὰ πασῶν ἅμα καὶ διὰ τεσσάρων λεγόμενον
οἱ Πυθαγορικοὶ μὲν σύμφωνον οὐκ ᾤοντο εἶναι,
διαφεῦγον πολλαπλάσιόν τε καὶ ἐπιμόριον λόγον καὶ ἔτι ἐπιμερῆ,
εἰς δὲ μικτὴν σχέσιν ἐκ πῖπτόν ἐστι·
καὶ γὰρ ὡς η′ πρὸς γ′,
διότι τὰ μὲν Ϛ′ τοῦ γ′ διπλάσια, τὰ δὲ η′ τοῦ Ϛ′ ἐπίτριτα·
εἰς δ’ οὖν τὸ παρὸν κατὰ τοὺς νεωτέρους νομιζέσθω καὶ αὐτὸ σύμφωνον,
σαφηνείας ἕνεκα τῶν ἑξῆς.
Iamblichus, *In Nicomachi Arithmeticam* 120,18–121,1

Les Pythagoriciens ne considéraient pas ce qu’on appelle l’ «octave –quarte» comme une consonance, parce qu’elle échappe aux rapports multiples, épimores et même épimères, et réside dans une relation mixte (par exemple 8 envers 3, car 6 est le double de 3 et 8 l’épitrite de 6); mais en l’occurrence, conformément aux modernes, il faut la considérer elle aussi comme une consonance, pour la clarté de la suite.
[Iamblichus, trans. Vinel, 2014]

*Note on the musical interval 8/3 by Creese, 2010:*

> The emphasis on reason and perception raises difficulties of two sorts. The first is that if reason is taken to dictate that "concord" is in fact a mathematical category, and that it includes only multiple and epimoric ratios, then while it will accept the fifth (3:2) and the octave plus fifth (3:1) as concords, it will reject the octave plus fourth (8:3) even though it will accept the fourth (4:3).
> This proved to be one of the most difficult and contentious issues in mathematical harmonics, and one of the authors of this period who makes the most determined attempt to include the 8:3 interval among the concords unaccountably omits it from a canonic division which includes all the other concords within a two–octave range (this is Adrastus).
> (Creese, 2010, p. 14)

(iii) Thus, a ratio of $a$ to $b$ has anthyphairesis of length 3 if and only if the ratio of $a$ to $b$ is epimeric or poly–epimeric if and only if $a = kb + c$, with $c < b$, $k$ natural number and the ratio $b$, $c$ is epimoric or poly–epimoric.
The ratio of 8 to 5 has anthyphairesis of length 4, and the ratio of 19 to 12 has anthyphairesis of length 5.

**11.4.5.** *The upward inductive ladder from the simple musical intervals (multiples, epimoric) to the general arithmetical ratios leads to Proposition VII.2 and the discovery of the Euclidean algorithm*

[1] It is clear that the steps described by Theon naturally lead to an implied *upward inductive ladder:*

*Definition*. The numbers $a$, $b$, $a > b$, have *anthyphairesis of length* $n + 1$ if for some natural numbers $k$ and $c < b$, $a = kb + c$, and $b$, $c$ have anthyphairesis of length $n$.
Inductively, this upward process can be continued *ad infinitum,* resulting in ratios with anthyphairesis of length any natural number.

[2] For any such ratio $a$, $b$, it can be proved, mimicking the upward process of Proposition VII.2, that the last remainder of the anthyphairesis is a common divisor of $a$, $b$, and it can be proved, mimicking the downward process of Proposition VII.2, that the last remainder of the anthyphairesis is the greatest common divisor of $a$, $b$, and, in fact, that every divisor of $a$, $b$ divides it.

[3] The question is whether every arithmetical ratio $a$, $b$ with $a > b$ appears at some level. It is at precisely this moment that Proposition VII.2 shows, with the use of the *Principle of the Least*, that this is indeed so.

[4] In conclusion, the music considerations, which start with the simple anthyphairetic definition of the epimoric musical intervals in terms of the multiple musical intervals, lead with the inductive upward process and the introduction of the Principle of the Least at the final stage to Proposition VII.2 and to the discovery of the Euclidean algorithm and anthyphairesis.

[5] One may question: are the first two or three steps enough to lead to an understanding of the general concept of anthyphairesis of numbers? According to what we mentioned before, after the definition of the second crucial step of the anthyphairesis of two numbers $a$, $b$ with $a > b$, through the epimoric ratio, ratios are adequately exhibited by anthyphairesis of exactly three steps. Was this enough to lead them to the discovery of the general Euclidean algorithm?
We note, by way of comparison, that Euclid, in order to prove the general case of the fundamental Proposition VII.2, establishing the construction by means of the Euclidean algorithm of the Greatest Common Divisor of any two non–relatively prime numbers, limits himself to exhibiting only the case of a pair of numbers with anthyphairesis consisting of exactly three steps:

> ὁ μὲν ΓΔ [$b$] τὸν ΒΕ [$a$] μετρῶν λειπέτω ἑαυτοῦ ἐλάσσονα τὸν ΕΑ [$c$],
> [$a = k_0 b + c$, $b > c$]
> ὁ δὲ ΕΑ [$c$] τὸν ΔΖ [$b$] μετρῶν λειπέτω ἑαυτοῦ ἐλάσσονα τὸν ΖΓ [$d$],
> [$b = k_1 c + d$, $c > d$]
> ὁ δὲ ΓΖ [$d$] τὸν ΑΕ [$c$] μετρείτω.
> [$c = k_2 d$]
> Euclid, *Elements,* Proposition VII.2

Thus, Euclid considers that describing the method of constructing the Greatest Common Divisor of two numbers, by an anthyphairesis of length three paradigmatically, provides a satisfactory exhibition of the general method for finding the Greatest Common Divisor of two numbers by an anthyphairesis of any (finite) length. Similarly, Plato in the second hypothesis of the *Parmenides* (142b1–143a3), in order to exhibit the infinity of the anthyphairesis of the dyad One and Being, is content to describe carefully just the first two steps; similarly, Proclus in his *Platonic Theology* (cf. Negrepontis, in press).
There are at least two ancient specimens in which an anthyphairesis is extended upwards: a mathematical example can be found in Proposition XIII.5 of the *Elements* concerning the mean and extreme ratio, and a philosophical one in the second hypothesis of Plato's *Parmenides* (142c7–d9), where the philosophical anthyphairesis One to Being is extended upwards by considering One + Being.

**11.5.** *The Evolutionary Stages from empirical music to the Pythagorean arithmetized music theory of Hippasus' 4–chord to the discovery of the Euclidean algorithm to the birth of Number theory (Book VII of the* Elements*)*

| Empirical Music | → Hippasus' experiments and construction of 4–chord. | Pythagorean Arithmetized Music theory of the 4–chord | → Sweeping Pythagorean conjecture | Book VII before arithmetical anthyphairesis and Definition VII.21* | → Discovery of Euclidean algorithm VII.1 & 2 | Birth of Number theory (Book VII) |
|---|---|---|---|---|---|---|

**12.** *The presence of numerous mathematical peculiarities in Book VII confirms the musical origin of Arithmetic*

The evolution of empirical pre–Pythagorean music into the mathematically rigorous arithmetized theory of Hippasus' 4–chord was achieved through the discovery of musical anthyphairesis and its subsequent transfer to arithmetic. This arithmetized system was later generalized into the deductive Number Theory of Book VII of the *Elements* by means of the general arithmetical anthyphairesis. The structural trajectory of this development strongly confirms that Greek arithmetic evolved directly from, rather than merely being applied to, the arithmetized theory of music.
There are two possibilities:
either there was a theory of Numbers, containing Book VII in some form, and this was applied to Music,
or the theory of Numbers evolved somehow from the arithmetized theory of Music.
This is the point to bring into the fore the considerable peculiarities that are discerned in Book VII of the *Elements*, if we compare it with a normal presentation of basic number theory, including the Fundamental Theorem of Number Theory and ratios/fractions and rationals in modern mathematics. We will take stock of these peculiarities, and we will provide explanations for them.

**12.1.** *Question 1* (posed in Section 7.1). Why, in the *Elements,* the number, although defined in Definition VII.2 as a [finite] multitude of units, is nevertheless *represented,* not as a collection of units/dots (*psephides*), but *as a linear segment*?

*Answer*. *Chords* are continuous in appearance, like *line segments*, and according to Nicomachus' account, chords were turned into *numbers* by Pythagoras' experiments. Thus, the representation of a number as a line in the *Elements* is in reality its representation as a musical chord.

**12.2.** *Question 2* (posed in Section 7.1). Why does *Definition VII.21* of proportion exclude the simple but important proportion $a/a = b/b$, and is given only for *unequal ratios*?

*Answer.* This is clear from Plato's passage in *Republic* 531 a4–b1:

> Νὴ τοὺς θεούς, ἔφη, καὶ γελοίως γε, πυκνώματ' ἄττα
> ὀνομάζοντες καὶ παραβάλλοντες τὰ ὦτα, οἷον ἐκ γειτόνων
> φωνὴν θηρευόμενοι, οἱ μέν φασιν ἔτι κατακούειν ἐν μέσῳ
> τινὰ ἠχὴν καὶ σμικρότατον εἶναι τοῦτο διάστημα, ᾧ μετρητέον,
> οἱ δὲ ἀμφισβητοῦντες ὡς ὅμοιον ἤδη φθεγγομένων, ἀμφότεροι
> ὦτα τοῦ νοῦ προστησάμενοι.
> (Plato, *Republic* 531 a4–b1)
>
> They talk of something they call minims (*puknomata*) and,
> laying their ears alongside,
> as if trying to catch a voice (*phonen*) from next door,
> some affirm that they can hear a note (*echen*) between
> and that this is *the least interval* (*diastema*) and the unit of measurement,
> while others insist that *the strings now render identical sounds*, that a musical interval, as the name *"diastema"* connotes, consists of unequal chords/numbers, and so a theory, a system for musical intervals, cannot but deal only with dichords with unequal chords.

In fact, Theon of Smyrna in 81,6–11 is even more explicit:

> διαφέρει δὲ διάστημα καὶ λόγος, ἐπειδὴ
> διάστημα μέν ἐστι τὸ μεταξὺ τῶν ὁμογενῶν τε καὶ ἀνίσων ὅρων,
> λόγος δὲ ἁπλῶς ἡ τῶν ὁμογενῶν ὅρων πρὸς ἀλλήλους σχέσις.
> διὸ καὶ τῶν ἴσων ὅρων
> διάστημα μὲν οὐδέν ἐστι μεταξύ,
> λόγος δὲ πρὸς ἀλλήλους εἷς καὶ ὁ αὐτὸς ὁ τῆς ἰσότητος
> Theon 81,6–11
>
> The *[musical] interval* and *the ratio logos* differ in that
> *the [musical] interval is* contained *between* homogeneous and *unequal* terms,
> while *the ratio logos* simply links homogeneous terms to one another.
> This is why there is *no [musical] interval between equal terms,*
> but there is *a ratio/logos* between them, which is that of *equality.*
> [Theon of Smyrna, trans. Lawlor & Lawlor, 1979, with modifications by the authors]

**12.3.** *Question 4* (posed in Section 7.2). Why, in Propositions VII.1 & VII.2 of the *Elements,* is anthyphairesis defined ONLY for a dyad of *unequal numbers*, and *not for equal ones*?

*Answer.* Similarly to our answer to Question 2, the anthyphairesis of two numbers, in Propositions VII.1&2, is explicitly considered only for a dyad of given *unequal* numbers (*duo arithmon anison ekkeimenon*). The equivalence class of a ratio of numbers is, as we saw, determined by its sequence of quotients, and since the equality of the terms of the ratio is excluded, for musical reasons, equality of the dyad of numbers in anthyphairesis must be excluded, as well.

**12.4.** *Question 3* (posed in Section 7.1). Why is Definition VII.21 of proportion $a/b = c/d$ *given fully* only for the case $a < b$, $c < d$ in terms of "*a part of b*" and "*a parts of b*", but *quite inadequately* for the case $a > b$, $c > d$ in terms only of "*a multiple of b*"? For example, 2/4 = 3/6, 4/10 = 12/30, and 4/2 = 6/3 are defined, while 10/4 = 30/12 is not defined according to Definition VII.21?

*Answer.* In music, the two chords $a$, $b$ of a dichord do not form an *ordered* dyad, as in the case of a mathematical ratio where the ratios $a/b$ and $b/a$ are different, but an *unordered* dyad; there is no differentiation between $(a, b)$ and the reciprocal $(b, a)$.
This is expressed vividly by Theon 81,11–16:

> τῶν δὲ ἀνίσων
> διάστημα μὲν ἓν καὶ τὸ αὐτὸ ἀφ’ ἑκατέρου <πρὸς> ἑκάτερον,
> λόγος δὲ ἕτερος καὶ ἐναντίος ἑκατέρου πρὸς ἑκάτερον·
> οἷον ἀπὸ τῶν β′ πρὸς τὸ ἓν καὶ ἀπὸ τοῦ ἑνὸς πρὸς τὰ β′
> διάστημα ἓν καὶ τὸ αὐτό,
> λόγος δὲ ἕτερος,
> τῶν μὲν δύο πρὸς τὸ ἓν διπλάσιος, τοῦ δὲ ἑνὸς πρὸς τὰ β′ ἥμισυς.
> Theon 81,11–16
>
> Between unequal terms,
> the interval between one and the other is unique and identical,
> while the ratios vary and are inverse from one term to the other;
> thus, 2 to 1 and 1 to 2 have only a single and identical interval,
> But there are two different ratios, the ratio of 2 to 1 being double, whereas the ratio of 1 to 2 is half.
> [Theon of Smyrna, trans. Lawlor & Lawlor, 1979, with modifications by the authors]

Thus, a musical interval differs from an arithmetical ratio in two ways:

A musical interval is defined *only for unequal numbers* $a$, $b$, but an arithmetical ratio is defined for both equal and unequal numbers $a$, $b$; and among unequal numbers $a$, $b$, the musical interval is an *unordered* pair of numbers $a$, $b$, but an arithmetical ratio is an ordered pair of numbers $a$, $b$. Definition VII.21 of Euclid's *Elements* ostensibly is a definition of the proportion of ratios, but in reality, it is a definition of equality of musical intervals, exactly because the definition of proportion $a/b = c/d$ is given only for unequal numbers (because equality does not make sense for musical intervals), and is complete only for $a < b$, $c < d$, but quite deficient for $a > b$, $c > d$ (because a musical interval is an unordered pair).

**12.5.** *Question 6* (posed in Section 7.2). Why did Euclid enunciate the Definition VII.21 of proportion of numbers by introducing solely for this purpose the distinction between "*same part*" and "*same parts*", even though this distinction, as clarified by Proposition VII.4, *is not fundamental* for the definition of proportion, and it is not used after Proposition VII.13?

*Answer.* It is made clear from Proposition VII.4 that the distinction between $a$ part of $b$ and $a$ parts of $b$ is exactly the distinction between the dyad $a$, $b$ having anthyphairesis of length one and of length greater than one. To base the definition of proportion of numbers on this distinction in two classes, grossly disproportionate to each other, does not illuminate the content of the definition, and in fact, this distinction is not used in Book VII, beyond Propositions VII.5–10. On the other hand, this distinction, restricted to the ratios of Hippasus' 4–chord, namely the distinction between multiple and epimoric ratios, makes good sense. Thus, we suggest that the peculiarities described by the three Questions 2, 3, and 6 on the definition VII.21 of proportion are explained by the musical origin of this definition.

**12.6.** *Question 7* (posed in Section 7.3). Why is the operation * that corresponds to the multiplication of ratios restricted to "*homonymous*" ratios, namely only to multiplication of ratios of the form $(a_1/a_2) * (a_2/a_3) = a_1/a_3$, rather than being defined in general?

*Answer.* Because this is the only multiplication that makes sense for musical intervals, first the composition of two dichords forming a 3–chord, and then, by Proposition VII.14, that this composition is well defined for musical intervals. General multiplication of arithmetical ratios $(a/b)\cdot(c/d)$ would most naturally be defined as $(a\cdot c)/(b\cdot d)$, and then proved well defined for equivalence classes with no mathematical difficulty, but such general multiplication makes no musical sense for musical intervals.

**12.7.** *Question 8* (posed in Section 7.3). Why is the *addition of ratios* not defined at all in the *Elements,* although it would be quite easy, especially so since the Least Common Multiple of two numbers is constructed in Proposition VII.34?

*Answer.* Again, there would be no mathematical difficulty for Euclid to define $a/b + c/d$, with or without the use of the Least Common Multiple, and then to prove that addition is well defined on equivalence classes, BUT addition of ratios makes no sense in musical intervals.

**12.8.** *Question 9* (posed in Section 7.3). Why did Euclid prove *the commutativity of multiplication of numbers* with the quite unexpected proof by using Proposition VII.13 *Alternando/Enallax*? In what sense would the Alternando property be considered a *commutativity property*?

*Answer.* This proof of Proposition VII.16 indicates that Euclid regards Proposition VII.13 as expressing some sort of commutativity. It was observed in Section 4 that in Hippasus' 4–chord the octave is produced as the composition of a fourth with a fifth, and of a fifth with a fourth, and also that the fifth is produced as the composition of a fourth with a tone and of a tone with a fourth. The general enunciation of the commutativity of the fourth with the tone in Hippasus' 4–chord is precisely the arithmetic analogue of the Perturbed Proportion (Proposition V.23) for magnitudes. Thus, we state and prove the following:

*Proposition. The commutativity property of composition* for musical intervals is equivalent to the Ex Aequali/*Alternando/Enallax property* for musical intervals.
*Proof.*
(⇒) Assume the Commutativity property of composition, let A, B, C, D be a 4–chord with A/B = C/D. We want to prove that A/C = B/D. We note that the 4–chord A, B, C, D has the same structure as Hippasus' 4–chord 6, 8, 9, 12. We compose the interval B/C first with the interval A/B from the right, and next with the equal interval C/D from the left. By the Commutativity property, A/C = A/B * B/C = B/C * C/D = B/D. This is the Ex Aequali/Alternando/Enallax Property.
(⇐) Assume the Ex Aequali/Alternando/Enallax property, and let A, B, C, and A′, B′, C′ be two trichords, such that A/B = B′/C′ and B/C = A′/B′.
We want to prove equality of composition A/C = A′/C′.
(This is the exact arithmetical analogue of Proposition V.23 of the *Elements* for the Perturbed proportion for magnitudes).
We form the trichords AA′, BA′, CA′ and BA′, BB′, BC′.
By the Ex Aequali/Alternando/Enallax property, it follows from the assumption B/C = A′/B′ that CA′ = BB′.
Thus, the two trichords have the form AA′, BA′, CA′ and BA′, CA′, BC′.
Thus, we can form the 4–chord AA′, BA′, CA′, BC′.

Furthermore,
AA′/BA′ = A/B = B′/C′ = BB′/BC′ = CA′/BC′, by the assumptions.
By the Ex Aequali/Alternando/Enallax property,
A/C = AA′/CA′ = BA′/BC′ = A′/C′.
This is the Commutativity property for Composition.

This Proposition shows clearly the musical origin of Proposition VII.13 (*Alternando*/*Enallax*), and, therefore, the musical origin of the proof of Proposition VII.16 (*commutativity of multiplication of numbers).*

**12.9.** *Question* 10 (posed in Section 7.3). Why does the arithmetical Book VII of the *Elements* begin with *Propositions about ratios and proportions of numbers,* not with numbers themselves, and continues with propositions about numbers themselves? It is like today one would first deal with rational numbers and then with natural numbers themselves.

*Answer.* Book VII, or to be more precise, the first part of Book VII, Propositions VII.1–22, is strictly about ratios of numbers, not about numbers themselves, and only the second part of Book VII deals with Propositions on numbers themselves, such as Propositions VII.27, 30, 31, 34.
According to Proclus, in his *Commentary to Euclid* 35,28–36,1:

> καὶ τὴν μὲν ἀριθμητικὴν τὸ καθ' αὑτὸ τὸ ποσὸν θεωρεῖν,
> τὴν δὲ μουσικὴν τὸ πρὸς ἄλλο
> Proclus, *Commentary to Euclid*, 35,28–36,1
>
> Arithmetic deals with numbers themselves (*to poson kath'auto*),
> Music with ratios (*to pros allo*) of numbers
> (Translation by Morrow, 1992)

Thus, the reason for this peculiarity is that
*(a) the first half* of Book VII, Propositions VII.1–22, mostly due to the early Pythagoreans, deals with ratios, because it has its origins in music, while,
*(b) the second half* of Book VII, Propositions VII.23–39, consists of pure number theory, containing Proposition VII.27, which is basic for Book VIII, generally considered to be due to Archytas, and Propositions VII.30, 31, which together with IX.14, a Corollary to VII.39, come quite close to the Fundamental Theorem of Arithmetic,
*(c)* the propositions of the first part about ratios are used for the proofs of the propositions of the second part about pure Number Theory.

**12.10.** *Conclusion*

At the end, it all boils down to a single fundamental observation: all the peculiarities of Book VII are well explained by realizing that they have a musical origin. Thus, the weight of these features strongly suggests that what happened in reality is the first alternative: Pythagorean music was arithmetized first, and was subsequently generalized to the birth of Arithmetic.
Ultimately, the long array of striking structural and conceptual peculiarities embedded within Book VII of the *Elements* cannot be dismissed as arbitrary stylistic choices. Every single unexpected anomaly – from the restriction of the definition of proportion to strictly unequal numbers, to the complete absence of addition of ratios, and the unexpected proof of the commutativity of multiplication via the *Alternando* property – finds a systematic and coherent explanation only when viewed through the lens of its musical origin. The cumulative weight of these peculiarities strongly reinforces our core thesis: early empirical Pythagorean music was arithmetized first through the construction of Hippasus' 4–chord and the discovery of musical anthyphairesis, providing the foundational framework that directly facilitated the birth of ancient Greek Number Theory.

**13.** *The key insight of Szabó (1978)*

Szabó (1978) was the first with the key insight that the arithmetical anthyphairesis (Euclidean algorithm) has its origin in Pythagorean music, yet relied on an anachronistic monochord model rather than the foundational evidence like Hippasus' 4–chord, Theon's account, and the peculiarities of Book VII. The main historians of Greek mathematics such as Heath (1921), van der Waerden (1943; 2012), and Knorr (1975) overlooked this connection.
According to Crocker (1963):

> The Pythagoreans based their musical theorems directly on their arithmetic. As a result, the kind of problem that could be taken up, as well as the kind of solution that could be found for it, depended upon the nature of pythagorean arithmetic. (p. 192)
>
> In this respect too, the Pythagoreans theory of consonance, dealing as it did with the simplest intervals, was the musical corollary of their arithmetic. (p. 195)

The thesis that "the Euclidean algorithm may well have originated in the Pythagorean theory of music" was first enunciated by Szabó (1978). This is how Szabó described the thesis:

> The sequence of steps that make up this method may be summarized as follows: (a) On discovering the shorter section of string which, together with the whole string, yielded the desired consonance, one was actually faced with two lines of different length. The one was the whole *monochord,* and the other was the section which produced the second note. (b) The shorter line was considered to be the unit and was subtracted from the longer one. This left the piece of string, which had been kept from vibrating, as a remainder. (c) To ascertain the length of the remainder, one tried to see *how many times it could be subtracted from the shorter of the two lines*.

This could be done twice in the case of the fifth and three times in the case of the fourth, because in the former, the piece of string which did not vibrate was *half* as long as the one which produced the second note, whereas in the latter, it was only a *third* as long. Thus, the shorter line was first subtracted from the longer one, and then the remainder was subtracted from the shorter, until in the end nothing was left over. The method described above is well known to readers of Euclid. It is the so–called *"Euclidean algorithm",* or method of successive subtraction, as it is termed in the more recent literature on the history of Greek mathematics.

(Szabó, 1978, pp. 134–135)

As far as I know, the origin of this method has never been elucidated. (Szabó, 1978, p. 135)

The explanation presented above of the musical terms *hemiolion* and *epitriton diastema* suggests that the Euclidean algorithm (successive subtraction) may well have originated in the Pythagorean theory of music. (Szabó, 1978, p. 136)

We make the following observations:

[1] Szabó bases his thesis on the assumption that Pythagoras possessed and utilized the monochord:

The invention of the canon is traditionally attributed to Pythagoras (and whatever one may think of the half–legendary figure of Pythagoras, this dating is not to be doubted). (Szabó, 1978, p. 138)

But according to Creese (2010):

One of the most problematic [arguments] was put forward by Arpad Szabó, who claimed that "the canon can be conclusively proved to have existed at least at the time of Plato" (p. 118).

The sole basis for this claim was a passage in the pseudo–Platonic *Epinomis,*

> …ἐν μέσῳ δὲ τοῦ ἓξ πρὸς τὰ δώδεκα συνέβη τότε ἡμιόλιον καὶ ἐπίτριτον…Plato, *Epinomis* 991a7–b1
>
> Between six and twelve comes the whole–and–a–half (9 = 6+3) and whole–and–a–third (8 = 6+2) (Translation by W.R.M. Lamb)

where it is observed that the hemiolic and epitritic ratios, by which the arithmetic and harmonic means are found, happen to fall within the ratio 6 to 12. By taking a passage of Gaudentius (Harm. 341.13–342.6: Pythagoras' first canonic division) as an indication that the division of the string of the monochord into twelve was more or less standard, Szabó concluded that the numbers 6 and 12 for the octave ratio in the *Epinomis* passage are only explicable if one assumes the existence of the monochord:

> I believe that these questions cannot be answered properly unless one bears in mind that the measuring instrument of the Pythagoreans (the canon), which was used to illustrate the proportional numbers of the consonances, was divided into twelve parts. In other words, Plato's remark proves convincingly that the canon existed at that time. Modern attempts to regard the canon as "an artificial piece of apparatus which was devised later" have not been successful. (Szabó, 1978, p. 119)

However Creese (2010) notes:

> Because he credited Gaudentius' account of the invention of the monochord with some historical value, accepting the date if not the attribution to Pythagoras, Szabó thought that the content of Gaudentius' division could be applied directly to the harmonic landscape of the fourth century B.C. simply because some of the numbers were the same. But Gaudentius himself evidently knew very little about the instrument's early history: his repetition of the Pythagoras–and–the–forge myth proves as much. Szabó deduced no other evidence from any of the undisputed dialogues of Plato, nor did he attempt to integrate the monochord into the harmonics of Plato's generation. He may not have been wrong about the date, but his argument fell far short of proving that he was right.
> (Creese, 2010, p. 102–103)

Thus, it is not clear at all that the monochord existed during the time of Pythagoras, not even during the time of Plato.

However, the problematic status of the monochord does not necessarily invalidate Szabó's argument. As we have noted, a number in music is represented by a chord/line, and a musical interval is represented by a dichord; this representation helps construct the Greatest Common Divisor of two numbers, especially for the simple multiple 2/1 and epimoric ratios 3/2, 4/3, 9/8 that arise in the arithmetized theory of the 4–chord. Thus, Szabó's argument can be modified so as to avoid the use of the monochord.

[2] Szabó mentions the two relations,

> octave = fifth composed with fourth, and fifth = fourth composed with tone,

but he does not correlate these with Hippasus' 4–chord and does not realize that these two relations mark the discovery of the concept of anthyphairesis. (Sections 3, 4 & 5)

> These ancient names for the three most important consonances are all mentioned in a fragment of Philolaus from which I would like to quote a few words here: The extent of an *octave* is *a fourth* and a *fifth…* (Szabó, 1978, p. 110)
>
> We have already quoted a fragment of Philolaus which states that in musical practice the octave consisted of the combination of a fourth and of a fifth. (Szabó, 1978, pp. 118–119)

Szabó (1978) mentions Hippasus only once:

> Although the concepts *diastema* and *horoi* originated in the course of experiments with the canon, Hippasus was already able to establish the most important consonances by means of bronze *diskoi…*
>
> (Szabó, 1978, p. 121)

Szabó, however, does not realize that these two equations,

2/1 = 3/2 * 4/3, and

3/2 = 4/3 * 9/8,

contained in Hippasus' 4–chord constitute musical anthyphairesis, in fact marking the discovery of the concept of anthyphairesis!

[3] Szabó does not mention the evidence in favor of the thesis that the arithmetical anthyphairesis (Euclidean algorithm) given by the careful description of the first three steps (from multiple ratios to epimoric, poly–epimoric, and to epimeric, poly–epimeric) of the upward ladder of constructing the general arithmetical anthyphairesis given by Theon of Smyrna in *Expositio rerum mathematicarum ad legendum Platonem utilium, [De utilitate mathematicae]*. This is surprising, considering that Theon's single work is cited and used for other questions in Szabó (1978). (Section 11)

[4] Szabó does not mention the peculiarities (Sections 7 & 12) in the arithmetical theory of ratios, as given in Book VII of the *Elements*, which, as we saw, are explained by its musical origin.

Thus, Szabó had the right intuition, the key insight about the musical origin of the arithmetic of Book VII, but did not really provide the required arguments in its favor.

Barker (1989) notes:

> Music enters the matter with the discovery that the relations between notes framing an organized melodic structure can themselves be expressed in very simple and neat numerical formulae. The lengths of two sections of a string giving notes an octave apart are in the ratio 2:1, while the ratio 3:2 gives a fifth and 4:3 a fourth. These fundamental harmonic relations thus correspond to what are evidently elegant and fundamental mathematical relations, and encourage the idea that all properly harmonic intervals gain their musical status because of their mathematical properties. […] The order found in music is a mathematical order; the principles of the coherence of a coordinated harmonic system are mathematical principles. And since these are principles that generate a perceptibly beautiful and satisfying system of organization, perhaps it is these same mathematical relations, or some extension of them, that underlie the admirable order of the cosmos, and the order to which the human soul can aspire.
These ideas fueled enthusiasm for investigations in mathematics, together with less rationalistic speculations about the symbolic import of individual numbers. Particular attention was focused on the mathematics of ratio and proportion, with constant reference back to the paradigmatic ratios governing basic musical structures. (Barker, 1989, p. 6)

**14.** *Conclusion: The theory of numbers in Book VII of the* Elements *has its origin and birth in the Pythagorean arithmetization of music*

Our conclusion that Number Theory, as presented in Book VII of the *Elements*, emerged from Pythagorean music rests on the following findings:
[1] Musical anthyphairesis was discovered by Hippasus in the course of the construction of his 4–chord. That Hippasus' 4–chord was seen as the first two steps of the musical anthyphairesis of the octave to the fifth is ascertained by the fact that these two steps were continued later by Philolaus in his *Fragment 6*.
[2] The transfer of the musical anthyphairesis of the octave to the fifth to the arithmetical anthyphairesis of the musical intervals themselves, the multiple 2/1 of length 1, the epimoric 3/2, 4/3, and 9/8, and of the fifth, fourth, tone of length 2.
[3] This transfer is helped

(a) by the additive language describing the multiplicative operation of composition of musical intervals [(3/2) * (4/3) = 2/1 is described as synthesis, (3/2) / (4/3) = 9/8 is described as subtraction/division of musical intervals], and

(b) by the representation of the numbers as chords/lines, a representation that illuminates the nature of the Greatest Common Divisor.

[4] Theon's detailed account (76.8–80,6) of the upward ladder from multiple ratios to epimoric and poly–epimoric ratios to epimeric and poly–epimeric ratios, a ladder that may be continued *ad infinitum.*
[5] The completion of the proof of VII.1 & 2 by reversing the upward process to obtain arithmetical anthyphairesis of two numbers, introducing the Principle of the Least.
[6] Propositions VII.A, B, 11, 12, 13, 20, the immediate consequences of the Principle of the Least and arithmetical anthyphairesis are natural generalizations, to all pairs of numbers, of the propositions that hold for the arithmetized theory of the four musical intervals of the 4–chord.
[7] A large number of striking peculiarities in Book VII (not mentioned by Szabó or other scholars on the mathematics of Pythagorean music), all explainable in terms of the musical origin of Book VII, from the musical theory of Hippasus' 4–chord.

Thus, the birth of Number Theory in Book VII, based on the Principle of the Least and the Euclidean algorithm, leading essentially to the Fundamental Theorem of Arithmetic, has its origin in the arithmetized Pythagorean theory of music, essentially in the theory of Hippasus' 4–chord.

## **15.** *The crucial role of Hippasus in the Pythagorean mathematical discoveries*

### **15.1.** *The role of Hippasus in arithmetizing empirical musical theory and in discovering the musical form of anthyphairesis*

This role has been extensively analyzed in Sections 3, 4, 5, and 6 of this paper.

### **15.2.** *The role of Hippasus in the discovery of incommensurability*

There are strong arguments, on the basis of persistent stories about the Pythagoreans, that Hippasus was the discoverer of the epoch–making proof of incommensurability of the diameter to the side of a square, by transferring the arithmetical concept of anthyphairesis to geometry, and by proving that the geometric anthyphairesis of the diameter to the side of a square is infinite, and by actually computing it. (Negrepontis & Farmaki, 2025)

### **15.3.** *The role of Hippasus in the formation of the two Pythagorean philosophical principles, Infinite and Finite*

Furthermore, according to ancient sources:

Ἵππασος Μεταποντῖνος καὶ αὐτὸς Πυθαγορικός. ἔφη δὲ
χρόνον ὡρισμένον εἶναι τῆς τοῦ κόσμου μεταβολῆς
καὶ πεπερασμένον εἶναι τὸ πᾶν καὶ ἀεικίνητον.
Diogenes Laërtius, *Lives of Eminent Philosophers* 8.84, 1–3

Hippasus of Metapontium, a Pythagorean, believed that there is a definite time that the changes in the universe take to complete, and the universe is finite and ever in motion.
[Diogenes Laërtius, *Lives of Eminent Philosophers*, 8.84.1–3, trans. Hicks, 1925, with modifications by the authors]

Ἵππασος δὲ ὁ Μεταποντῖνος καὶ Ἡράκλειτος ὁ Ἐφέσιος ἓν καὶ οὗτοι καὶ κινούμενον καὶ πεπερασμένον, ἀλλὰ πῦρ ἐποίησαν τὴν ἀρχὴν καὶ ἐκ πυρὸς ποιοῦσι τὰ ὄντα πυκνώσει καὶ μανώσει καὶ διαλύουσι πάλιν εἰς πῦρ, ὡς ταύτης μιᾶς οὔσης φύσεως τῆς ὑποκειμένης· πυρὸς γὰρ ἀμοιβὴν εἶναί φησιν Ἡράκλειτος πάντα.
Simplicius, *Commentary to Aristotle's Physics* 23, 33–24,4

Hippasus of Metapontum and Heraclitus of Ephesus too said that [the principle] was one and moved and finite, but they made the principle fire, and they make the things that are come from fire by condensation and rarefaction, and they resolve them again into fire, on the grounds that this one nature is what underlies: for Heraclitus says that all things are an exchange for fire. [Simplicius, *In Aristotelis Physicorum*, trans. Menn, 2022, with modifications made by the authors]

Ἡράκλειτος δὲ καὶ Ἵππασος τὸ πῦρ (τὴν γὰρ γῆν οὐδεὶς
ἠξίωσεν ὑποθέσθαι μόνην διὰ τὸ δυσαλλοίωτον), τινὲς δὲ ἄλλο τι τῶν τριῶν ὑπέθεντο, ὅ ἐστι πυρὸς μὲν πυκνότερον, ἀέρος δὲ λεπτότερον, ἢ ὥς ἐν ἄλλοις φησίν, ἀέρος μὲν πυκνότερον, ὕδατος δὲ λεπτότερον.
Simplicius, *Commentary to Aristotle's Physics* 149,8 – 11

Heraclitus and Hippasus [postulated] fire (for no one thought it appropriate to postulate earth alone because it is hard to alter), while some posited some other of the three, which is denser than fire but finer than air, or as he says elsewhere, denser than air but finer than water.
[Simplicius, *In Aristotelis Physicorum*, trans. Menn, 2022, with modifications made by the authors]

Ἵππασος δὲ ὁ Μεταποντῖνος…ἓν εἶναι τὸ πᾶν ἀεικίνητον καὶ πεπερασμένον. Aetius, De placitis reliquiae 292,1–5

> According to Hippasus of Metapontium, … the universe is one, ever in motion, and finite.[Aetius, *De placitis reliquiae*, 292.1–5, ed. Diels, 1879/1965, authors' translation]

Thus, according to Hippasus, the universe is *both ever in motion* (reflecting the principle of the Infinite) *and finite* (reflecting the principle of the Finite). This cosmological thesis is the direct counterpart to the standard Pythagorean mathematical thesis of the two principles. As demonstrated in Negrepontis & Farmaki (2025), these Pythagorean Principles have their exact origin in the proof of the incommensurability of the diameter to the side of a square: more specifically, the Infinite corresponds to the infinite anthyphairesis, and the Finite to the preservation of Gnomons. Since the anthyphairetic proof of incommensurability has been attributed to Hippasus, it logically follows that these corresponding philosophical principles owe their origin to him as well.

*Tannery on Hippasus' "fire" as the precursor of the Pythagorean "apeiron"*

> According to Tannery, Hippasus' theory is rather evidence for the early stage of the Pythagorean theory of principles, in which *peras* and *apeiron* were still conceived in material terms: the solid earth as *peras* and the fluid and refined fire as *apeiron*. (Betegh, 2011, pp. 370–371)

**15.4.** *The transfer of musical anthyphairesis to arithmetical anthyphairesis (Euclidean algorithm) and the birth of Number theory (Book VII of Euclid's* Elements*)*

We have attempted a reconstruction of the Pythagorean discovery of arithmetical anthyphairesis (*Euclidean algorithm*) in Section 11. There is no ancient testimony that Hippasus discovered the arithmetical anthyphairesis. But Hippasus was the key person in the early step of the arithmetization of music and the construction of the 4–chord, leading to the discovery of the concept of (musical) anthyphairesis; and quite possibly Hippasus was the key person in the later steps of discovering and proving the incommensurability by anthyphairesis, and in formulating an early form of the Pythagorean philosophy based on the two principles Infinite and Finite. Thus, in accordance with our suggestion that the arithmetical anthyphairesis was discovered by a gradual transfer of the musical anthyphairesis, we may conjecture that Hippasus might be the person who made the intermediate discovery of the Euclidean algorithm.

In conclusion,

*the initial step* of the discovery of the *musical anthyphairesis* resulting in the 4–chord,

*the third step,* consisting of the transfer to the *geometric anthyphairesis* of the diameter to the side of a square, resulting in the Pythagorean proof of incommensurability, and

*the fourth step* resulting in the Pythagorean principles of Infinite and Finite, both having anthyphairetic content (Negrepontis & Farmaki, 2025), are all attributed by ancient sources to Hippasus. It is a reasonable conjecture that the remaining *second step*, the discovery of *arithmetic anthyphairesis/Euclidean algorithm*, is also due to the genius of Hippasus.

**15.5.** From our foregoing analysis, Hippasus emerges as the key figure in Pythagorean Mathematics and Philosophy.

**Acknowledgment.** The authors wish to express our deep thanks to Athanase Papadopoulos for his critical reading of a preliminary version of this manuscript. His highly constructive and illuminating comments have essentially improved the final text. Of course, all possible shortcomings still remaining are strictly our own.

## References

### Primary Sources

**Secondary Sources**